\documentclass[preprint,12pt]{article}

\usepackage{arxiv}
\usepackage[english]{babel}
\usepackage[utf8x]{inputenc}
\usepackage[T1]{fontenc}
\usepackage{amsmath,amsfonts,amssymb,systeme} 
\usepackage{graphicx} 
\usepackage{lineno} 
\usepackage{epstopdf} 
\usepackage{subcaption} 
\usepackage{hyperref}
\usepackage{cleveref}
\usepackage{xcolor}
\usepackage{makecell}

\graphicspath{{img/}}

\DeclareMathOperator{\atantwo}{atan2}

\title{Phase resonance nonlinear modes of mechanical systems}
\author{Martin Volvert\\
	Space Structures and Systems Laboratory,\\
	Aerospace and Mechanical Engineering Department,\\
	University of Liège, Belgium\\
	\texttt{m.volvert@uliege.be}
	\and Gaëtan Kerschen\\
	Space Structures and Systems Laboratory,\\
	Aerospace and Mechanical Engineering Department,\\
	University of Liège, Belgium\\
	\texttt{g.kerschen@uliege.be}}
\begin{document}
	\maketitle
    
    
    \begin{abstract}
	The resonances of forced dynamical systems occur when either the amplitude of the frequency response undergoes a local maximum (amplitude resonance) or phase lag quadrature takes places (phase resonance). This study focuses on the phase resonance of nonlinear mechanical systems subjected to single-point, single-harmonic excitation. In this context, the main contribution of this paper is to develop a computational framework which can predict the mode shapes and oscillation frequencies at phase resonance. The resulting nonlinear modes are termed phase resonance nonlinear modes (PRNMs). A key property of PRNMs is that, besides primary resonances, they can accurately characterize superharmonic, subharmonic and ultra-subharmonic resonances for which, as shall be shown in this paper, phase lags at resonance may be different from $\pi/2$. The proposed developments are demonstrated using one- and two-degree-of-freedom systems featuring a cubic nonlinearity. 
\end{abstract}


\section{Introduction}
\label{SECTION:INTRODUCTION}
Modal analysis has been, and continues to be, the dominant dynamical method used in structural design. The goal of modal analysis is to find the vibration modes, resonance frequencies and damping ratios of the considered system \cite{Ewins}. One key assumption of modal analysis is linearity, which, however, real-world structures violate because they may feature advanced materials, friction and contact \cite{MSSPREVIEW}. 

The theory of nonlinear normal modes (NNMs) was developed to generalize the concept of a vibration mode to nonlinear systems \cite{VAKAKISBOOK}. There exist two main definitions of a NNM based on either periodic motions or invariant manifolds. In direct analogy to a linear mode, Rosenberg defined a NNM as a synchronous vibration of the undamped, unforced system for which all points reach their extreme values or pass through zero simultaneously \cite{ROSENBERG3,ROSENBERG}. This definition was later generalized to non-necessarily synchronous periodic oscillations of the system to encompass modal interactions \cite{LEE+KERSCHEN}. The concept of a periodic solution was then used to calculate NNMs of dissipative systems through an additional damping term which is large enough to compensate for the nonconservative forces \cite{KRACK}. Based on the theory of invariant manifolds, Shaw and Pierre proposed an elegant generalization of a nonlinear mode to damped systems \cite{SHAW,SHAW2}. The manifold defining the NNM is invariant under the flow, i.e., the orbits that start out in the manifold remain in it for all times, which extends the invariance property of linear normal modes to nonlinear systems. Recently, Haller and Ponsioen provided exact mathematical existence and uniqueness conditions for NNMs defined as invariant manifolds \cite{HALLER}. 

The usefulness of the NNM theory extends to forced systems in the sense that their periodic responses are born out of periodic orbits of the conservative limit that are in resonance with the forcing \cite{NAYFEH,GH}. Stated otherwise, the undamped NNMs can serve as very good approximations to the loci of the amplitudes of the resonance peaks \cite{VAKAKIS,HILL}. However, Cenedese and Haller highlighted through Melnikov analysis that NNMs do not necessarily perturb into forced-damped periodic responses \cite{CENEDESE}. Similarly, Hill et al. evidenced that only those NNMs which are able to transfer a large amount of energy through the fundamental components of the response strongly relate to the forced dynamics \cite{HILL2}.

The present study adopts Rosenberg’s framework and considers the periodic orbits of harmonically-forced, damped systems near resonance. Specifically, the focus is on the {\it phase resonance} of such systems. For linear systems, phase resonance takes place when the excitation force vector balances the damping forces with the result that all degrees of freedom vibrate synchronously and have a phase lag of $\pi/2$ with respect to the excitation \cite{GeradinBook}. Because the corresponding structural deformation is the normal mode of the conservative system, there is thus a 1:1 correspondence between undamped normal modes and phase resonances. This theoretical result was exploited in academia and industry for experimental modal analysis through force appropriation \cite{WRIGHT,GOGE}. In the presence of nonlinearity, Wright et al. \cite{ATKINS} used optimization to derive a multi-harmonic forcing vector that can isolate a target linear mode by counteracting nonlinear coupling to other modes. When it comes to the experimental identification of nonlinear modes, Peeters et al. investigated the theoretical link between undamped NNMs and phase resonances \cite{PEETERSJSV}. They demonstrated using harmonic balance that the mode shape at phase resonance is a NNM provided that multi-point, multi-harmonic forcing is utilized so as to impose quadrature between {\it all} harmonics of the displacement and of the forcing. Even though single-harmonic forcing with a single exciter was able to identify relatively accurate approximations of isolated NNMs through modal indicator functions \cite{PEETERS}, identified models \cite{NOEL,SZALAI}, control-based continuation \cite{RENSON} and phase-locked loops \cite{THOMAS,SCHEEL}, Renson et al. \cite{RENSON2} and Ehrhardt and Allen \cite{EHRHARDT} confirmed numerically and experimentally that significant errors between NNMs and quadrature curves can exist under such a forcing.

Because the need for distributed forcing with multiple harmonic components complicates the development of practical nonlinear experimental modal analysis methodologies, the phase resonance of nonlinear systems subjected to widely-used single-point, single-harmonic forcing is considered herein. In this context, the main contribution of this paper is to develop a computational framework which can provide accurate numerical approximations of the mode shapes and oscillation frequencies at phase resonance. The resulting nonlinear modes are different from NNMs; they are termed phase resonance nonlinear modes (PRNMs). A key property of PRNMs is that, besides primary resonances, they can accurately characterize superharmonic, subharmonic and ultra-subharmonic resonances for which, as shall be shown in this paper, phase lags at resonance may be different from $\pi/2$.

The paper is organized as follows. Section 2 considers a two-degree-of-freedom system featuring cubic nonlinearity as a motivating example. Section 3 analyzes in great detail the different families of $k:\nu$ resonances of a Duffing oscillator and pays particular attention to the phase lag between the $k$-th harmonic of the displacement and the forcing. Section 4 defines the concept of a PRNM and introduces a harmonic balance-based algorithm for computing PRNMs numerically. Section 5 reconsiders the two-degree-of-freedom system of Section 2 under the PRNM framework. Finally, conclusions are drawn in Section 6.

    \section{A motivating example}\label{motiv}

The primary resonances of forced dynamical systems occur when either the amplitude of the frequency response undergoes a local maximum (amplitude resonance) or phase lag quadrature takes place (phase resonance). This section aims to show that, although they exist in the neighborhood of primary resonances, the NNMs correspond neither to amplitude resonance nor to phase resonance under classical harmonic forcing. To support this assertion, a two-degree-of-freedom system with a cubic spring attached to the first mass is considered:
\begin{eqnarray}\label{eq:EOM2D_0}
    \ddot{x}_1+0.02\dot{x}_1-0.01\dot{x}_2+2x_1-x_2+x^3_1=f\sin{\omega t}\\
    \ddot{x}_2+0.11\dot{x}_2-0.01\dot{x}_1 + 2 x_2-x_1=0
\end{eqnarray}
The natural frequencies of the underlying linear system are 1 and $\sqrt{3}$ rad/s.

The NNMs of the undamped, unforced system and the damped, forced  responses, termed nonlinear frequency response curves (NFRCs), were calculated for different forcing amplitudes using numerical continuation \cite{PEETERSMSSP}. Figure 1a depicts the NFRCs in the vicinity of the first primary resonance together with the NNMs and the loci of amplitude/phase resonance points. At relatively low forcing, the agreement between NNMs, amplitude and phase resonances is excellent. However, for higher forcing, NNMs deviate significantly from amplitude and phase resonances. In the presence of modal interactions, the discrepancy can be even more important \cite{RENSON2}. Considering now the highest forcing level in Figure 1a, the energy balance principle proposed in \cite{HILL} can calculate which NNM motion and frequency correspond to this specific forcing amplitude. The result of energy balance, represented by the green point in Figure 2, evidences that the NNM does not lie on the considered NFRC.

The NFRCs in Figure 1a also feature resonance peaks around $\sqrt{3}/3$ rad/s, which signals the presence of a 3:1 superharmonic resonance of the second primary resonance. This secondary resonance cannot be captured using Rosenberg's NNMs, but it is interesting to observe that phase quadrature between the third harmonic of the displacement and the forcing exists in the vicinity of this resonance peak.

In summary, this simple example served to illustrate that (i) NNMs may present deficiencies when the resonant response under single-point, single-harmonic forcing is to be captured accurately and (ii) phase quadrature curves could provide a generic way of characterizing the different families of resonance of a nonlinear system.

\begin{figure}[htbp] 
  \begin{subfigure}[b]{1\linewidth}
    \centering
    \includegraphics[scale=0.74]{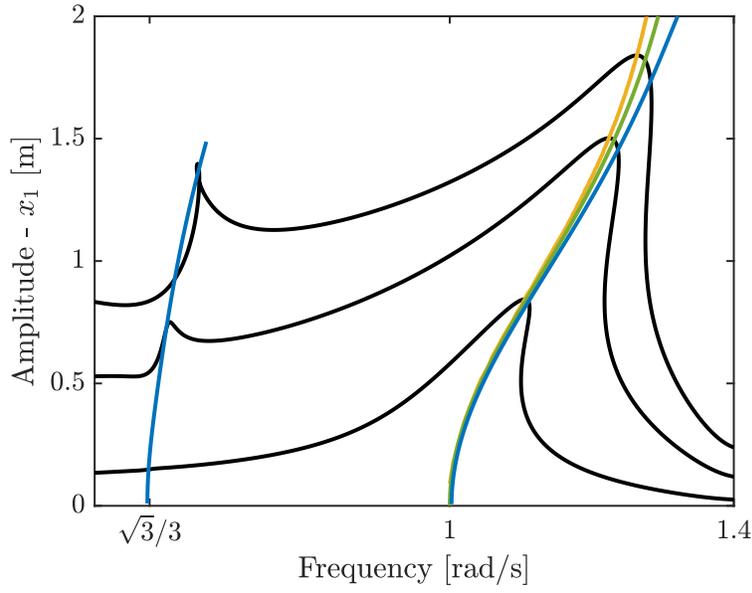}
    \caption{\label{fig:NNM_SUPERH_HIGH_FORCING}}
  \end{subfigure} 
  \begin{subfigure}[b]{1\linewidth}
    \centering
    \includegraphics[scale=0.74]{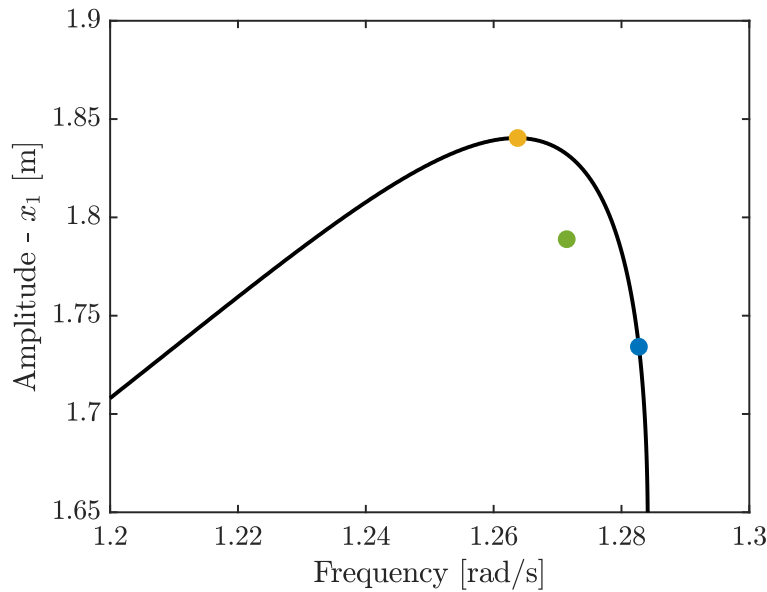}
    \caption{\label{fig:NNM_SUPERH_HIGH_FORCING_ZOOM}}
  \end{subfigure}
  \caption{Two degree-of-freedom system with a cubic spring: NFRCs (black), NNMs (green), amplitude resonance (yellow) and phase quadrature points (blue); (a) $f=0.161, 0.8, 1.5$ N; (b) close-up around the first primary resonance, $f=1.5$ N.}
  \label{fig:NNM_PRNM_2DOF}
\end{figure}

\section{Phase lags of the resonances of the Duffing oscillator}
\label{SECTION:DYNDUFFING}

The objective of this section is to analyze in detail the different resonance families of a harmonically-forced Duffing oscillator. The governing equation of motion is
\begin{equation}
    m\ddot{x}(t)+c\dot{x}(t)+kx(t)+k_{nl}x^3(t)=f\sin{\omega t}
    \label{eq:EOM1D}
\end{equation}
where $m$, $c$, $k$ and $k_{nl}$ represent the mass, damping, linear and nonlinear stiffness coefficients, respectively. $f$ is the forcing amplitude whereas $\omega$ is the excitation frequency of period $T$. The natural frequency of the undamped, linearized system is $\omega_0=\sqrt{\frac{k}{m}}$. The coefficients are set to $m=1$kg, $c=0.01$kg/s, $k=1$N/m and $k_{nl}=1$N/m$^3$ throughout the present study.

We consider the Fourier decomposition of the displacement
\begin{equation}
    x(t) = \frac{c_0}{\sqrt{2}}+\sum_{k=1}^{\infty} \left(s_k \sin \omega_k t + c_k \cos \omega_k t\right), 
    \label{eq:FourierDisp}
\end{equation}
The frequency of the $k$-th harmonic of the displacement is $\omega_k=\frac{k\omega}{\nu}$ where $\nu$ is a positive integer which allows to account for subharmonic components. We note that $k$ and $\nu$ should be relatively prime integers. The periods $T_k$ and $T$ are related through $T_k=\frac{\nu}{k}T$. Equation (\ref{eq:FourierDisp}) shows that each harmonic $k$ may trigger a resonance if $\omega_k$ corresponds to the (amplitude-dependent) frequency of the primary resonance of the system. According to Stoker \cite{Stoker}, the resonances can be divided into four categories, namely
$1:1$ primary resonance ($k=\nu=1$), $k:1$ superharmonic or ultraharmonic resonances, $1:\nu$ subharmonic resonances and $k:\nu$ ultra-subharmonic resonances. 

\begin{figure}[htbp] 
  \begin{subfigure}[b]{0.5\linewidth}
    \centering
    \includegraphics[width=1\linewidth]{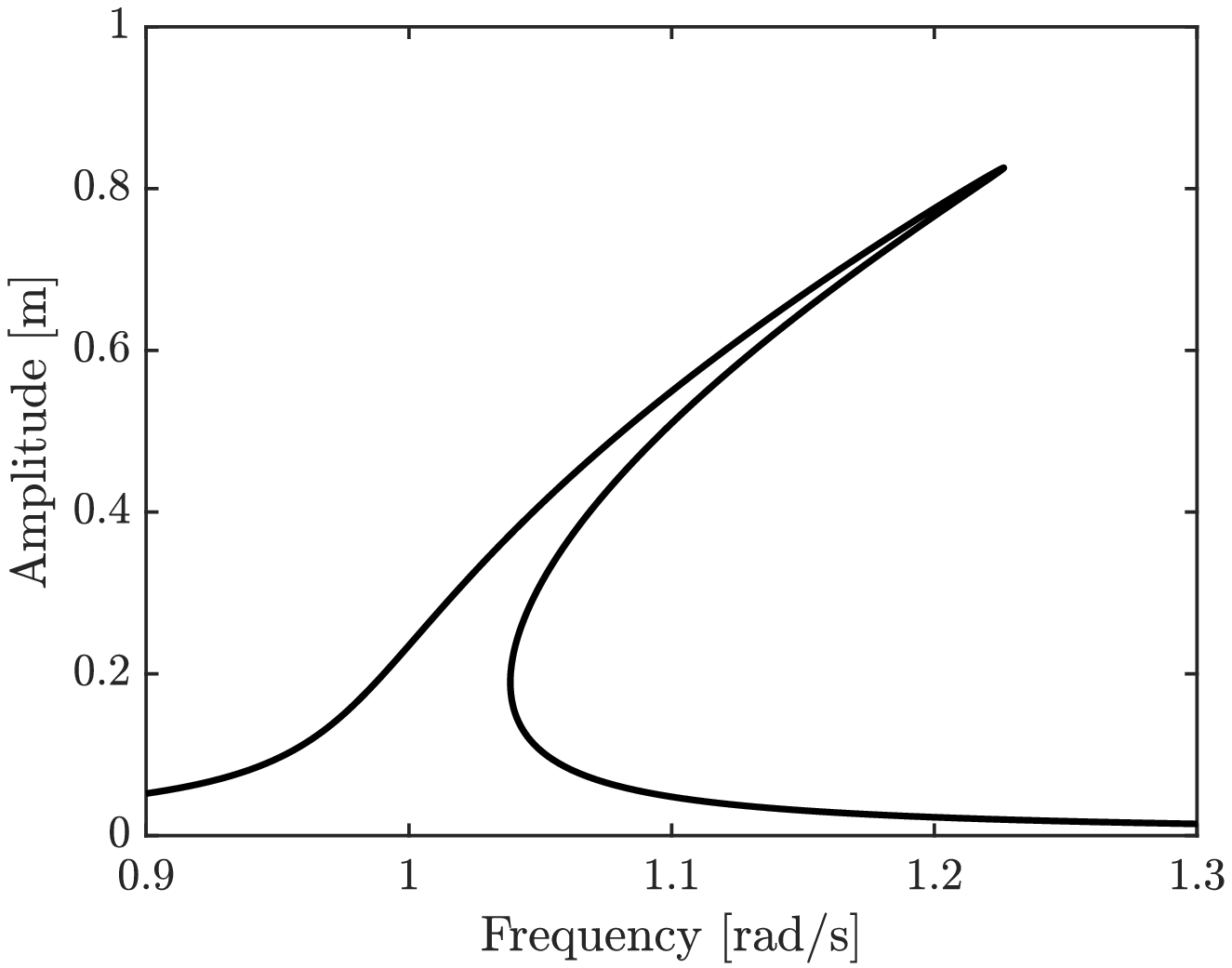} 
    \caption{\label{fig:F_0_01}} 
  \end{subfigure}
  \begin{subfigure}[b]{0.5\linewidth}
    \centering
    \includegraphics[width=1\linewidth]{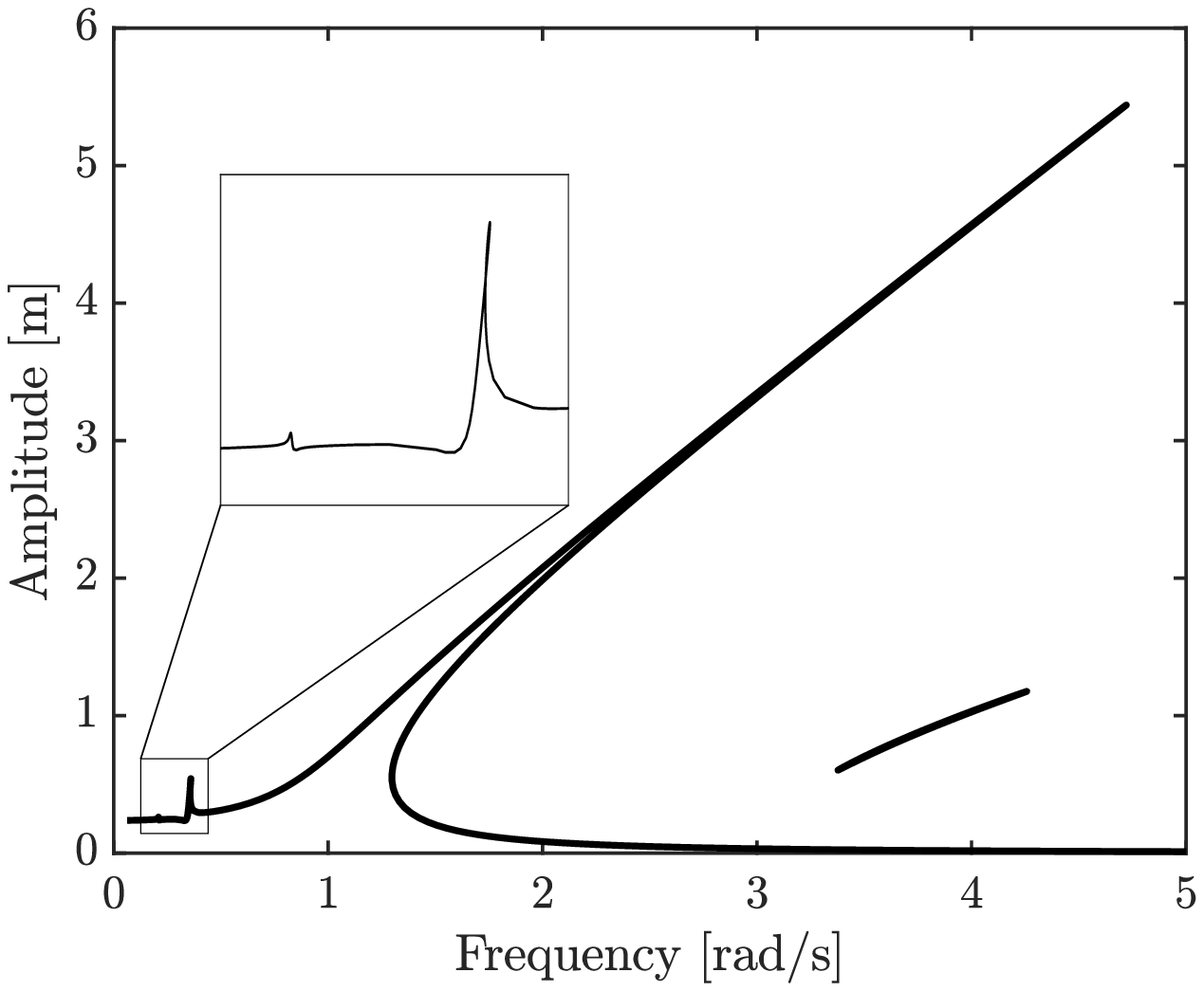} 
    \caption{\label{fig:F_0_25}} 
  \end{subfigure} 
  \begin{subfigure}[b]{0.5\linewidth}
    \centering
    \includegraphics[width=1\linewidth]{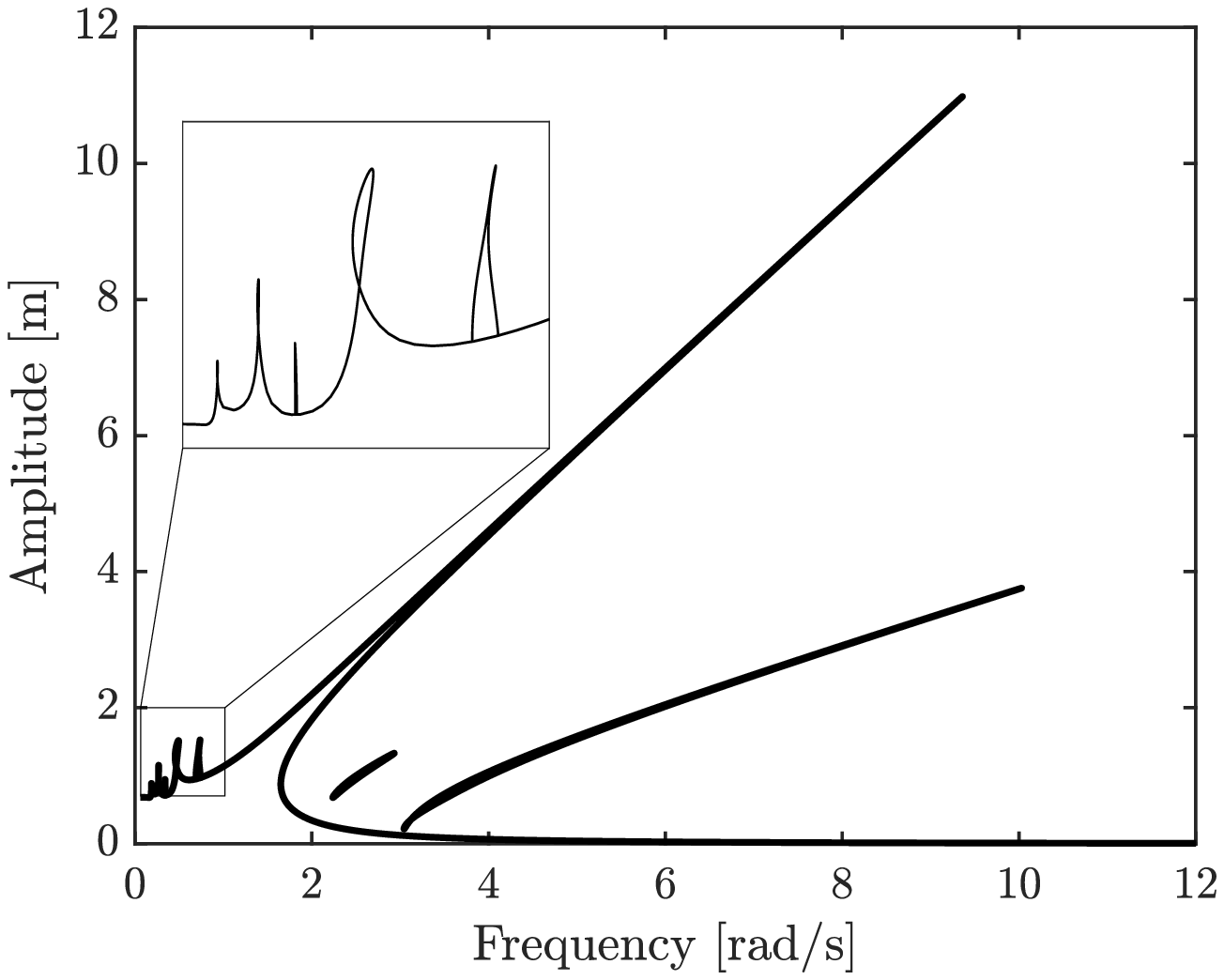} 
    \caption{\label{fig:F_1} } 
  \end{subfigure}
  \begin{subfigure}[b]{0.5\linewidth}
    \centering
    \includegraphics[width=1\linewidth]{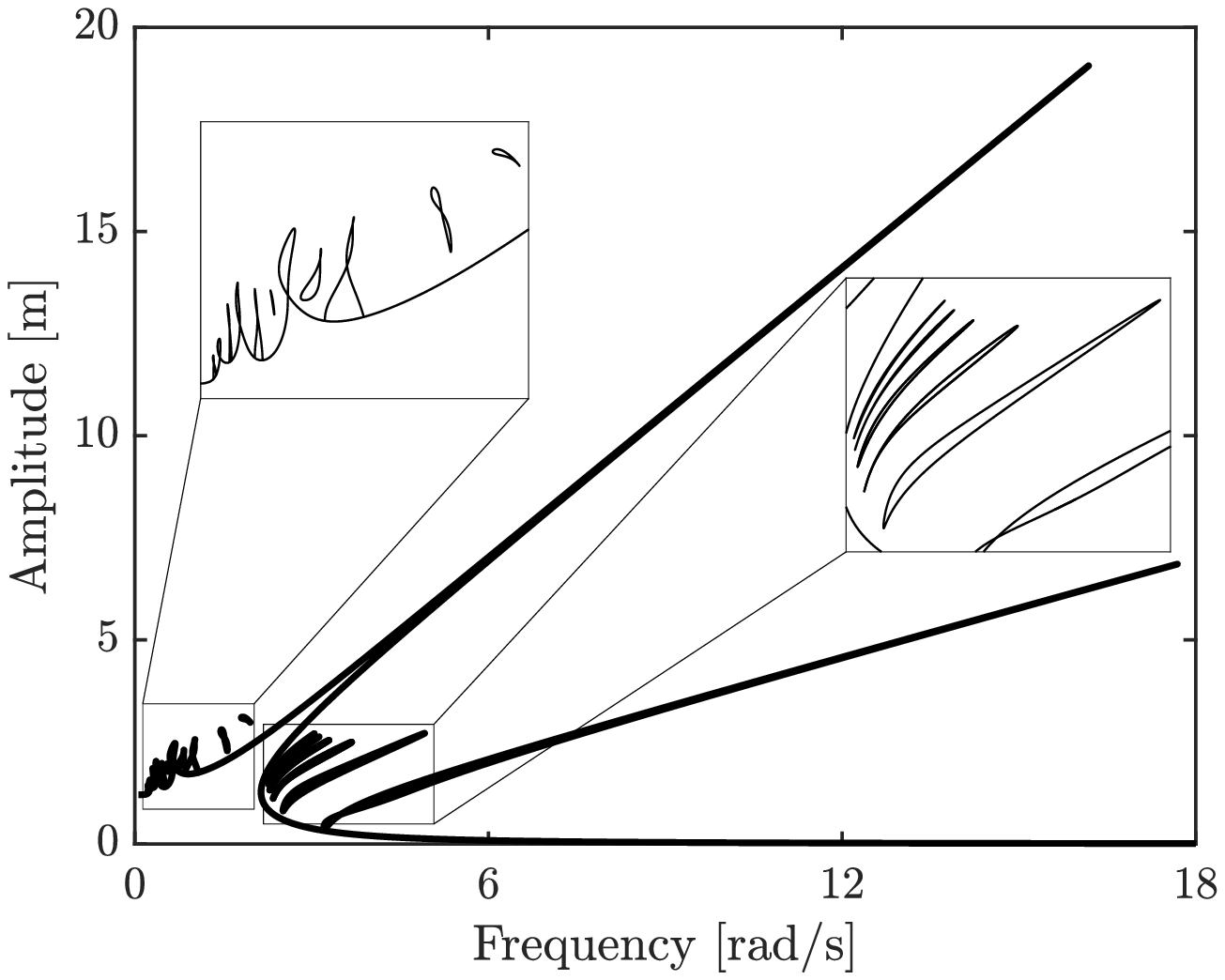} 
    \caption{\label{fig:F_3}} 
  \end{subfigure}
  \begin{subfigure}[b]{0.5\linewidth}
	\centering
	\includegraphics[width=1\linewidth]{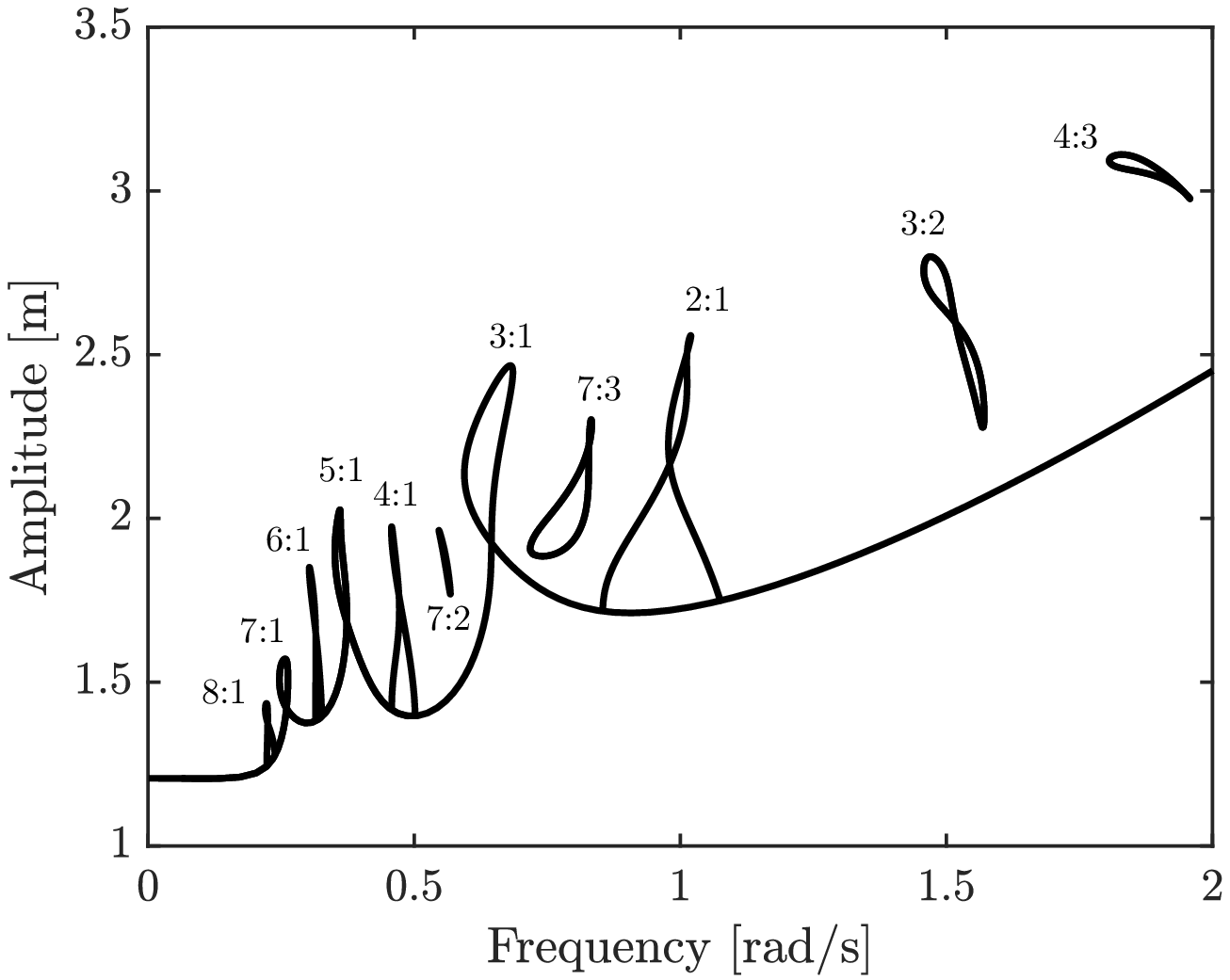} 
	\caption{\label{fig:F_3_ZOOM_SUP}}
  \end{subfigure}
  \begin{subfigure}[b]{0.5\linewidth}
	\centering
	\includegraphics[width=1\linewidth]{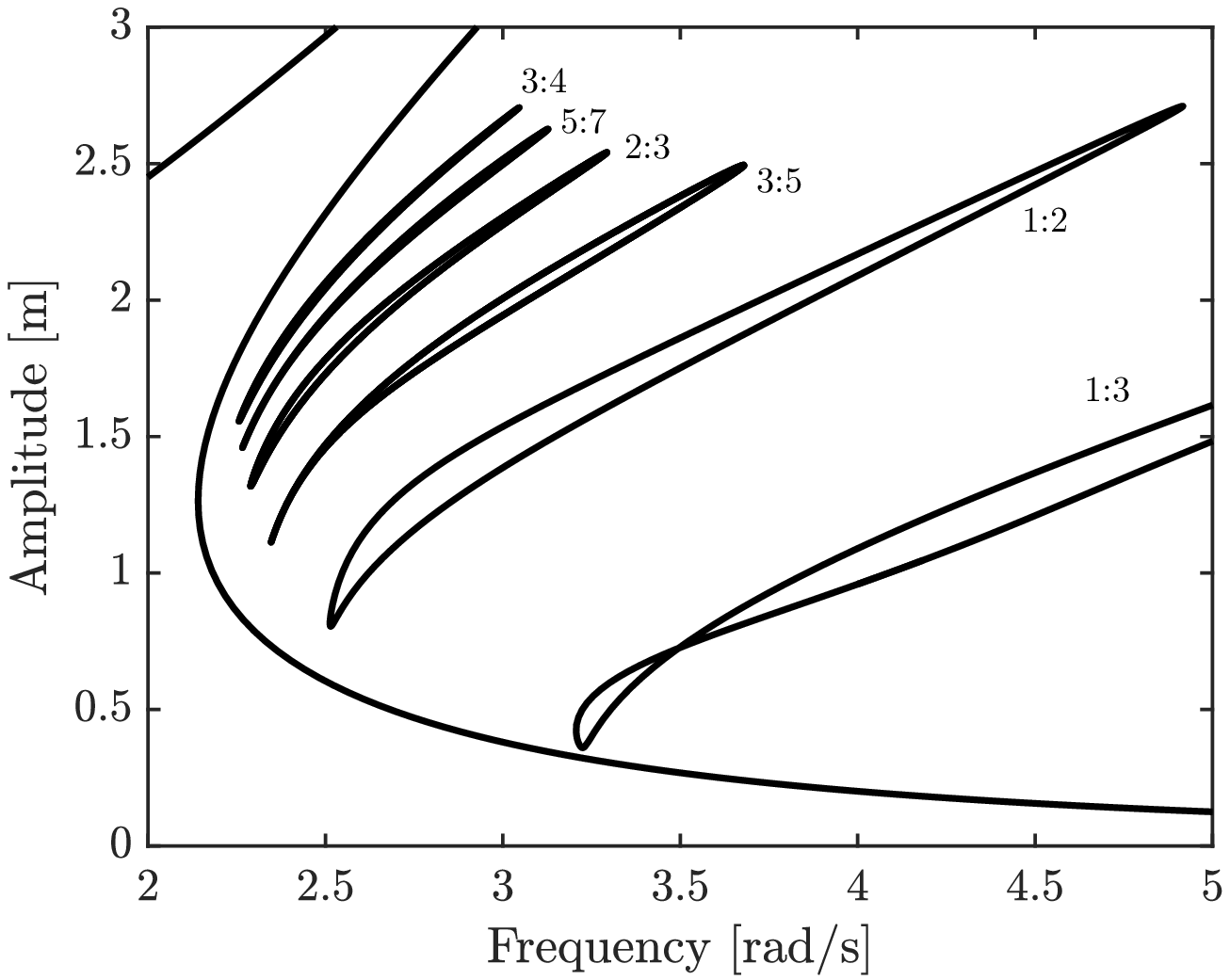}
	\caption{\label{fig:F_3_ZOOM_SUB}}
  \end{subfigure}
  \caption{NFRCs of the Duffing oscillator: (\subref{fig:F_0_01}) $f=0.01$N, (\subref{fig:F_0_25}) $f=0.25$N, (\subref{fig:F_1}) $f=1$N and (\subref{fig:F_3}) $f=3$N, (\subref{fig:F_3_ZOOM_SUP}) close-up before the primary resonance peak, (\subref{fig:F_3_ZOOM_SUP}) close-up after the primary resonance peak.}
  \label{fig:NFRCs} 
\end{figure}

The system (\ref{eq:EOM1D}) is analyzed at different forcing amplitudes $f$, i.e., 0.01N, 0.25N, 1N and 3N. The NFRCs computed using the harmonic balance method (HBM) with $8\times\nu$ subharmonics combined with pseudo arc-length continuation are depicted in Figure \ref{fig:NFRCs}. The computation of basins of attraction was necessary to initiate the continuation of the branches which are disconnected from the main NFRC. 

For a forcing amplitude of 0.01N in Figure \ref{fig:F_0_01}, the only nonlinear effect appearing in the NFRC is the hardening of the primary resonance. At 0.25N in Figure \ref{fig:F_0_25}, 3:1 superharmonic and 1:3 subharmonic resonance branches appear before and after the primary resonance, respectively. It should be noted that the subharmonic resonance is isolated from the main curve. Additional branches corresponding to 2:1, 4:1, 5:1 and 7:1 superharmonic and 1:2, 1:3 subharmonic resonances arise in Figure \ref{fig:F_1} at 1N. Finally, when the forcing amplitude is 3N,
7:2, 7:3, 3:2, 4:3, 3:4, 5:7, 2:3 and 3:5 ultra-subharmonic resonances which are all isolated from the main branch can be observed in Figure \ref{fig:F_3}. Interestingly, each resonance branch $k:\nu$ is an isolated solution when $\nu>1$ \cite{Marchionne}. 


All these resonances are examined in greater detail hereafter. Similar analyses were achieved in a number of studies, e.g., analytically using multiple scales \cite{NAYFEH}, higher-order averaging \cite{YAGA1,YAGA2} or harmonic balance \cite{LEUNG}, and numerically \cite{Marchionne,PARLITZ,GUILLOT}. However, with the exception of \cite{LEUNG} for odd superharmonic and subharmonic resonances, very few attention was devoted in the existing body of literature to the careful analysis of the phase lag in the vicinity of the $k:\nu$ resonance. Specifically, our interest is in the phase lag between the $k$-th harmonic of the displacement and the forcing.

\begin{figure}[htpb]
    \centering
    \includegraphics[width=1\linewidth]{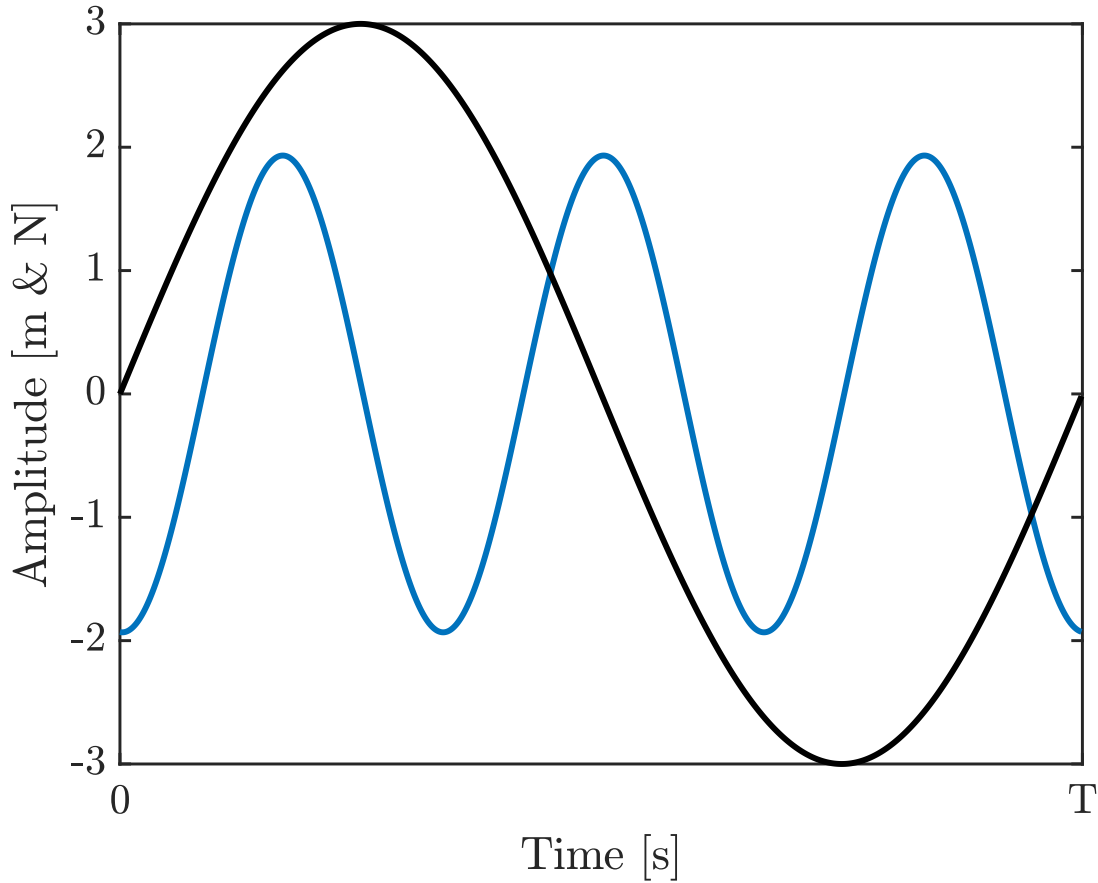}
    \caption{Time series corresponding to the 3:1 superharmonic resonance: forcing ($f=3$N, black) and third harmonic of the displacement ($k=3$, blue).}
    \label{fig:HAR3_FORCING}
\end{figure}

\begin{figure}[ht] 
  \begin{subfigure}[b]{0.5\linewidth}
    \centering
    \includegraphics[width=1\linewidth]{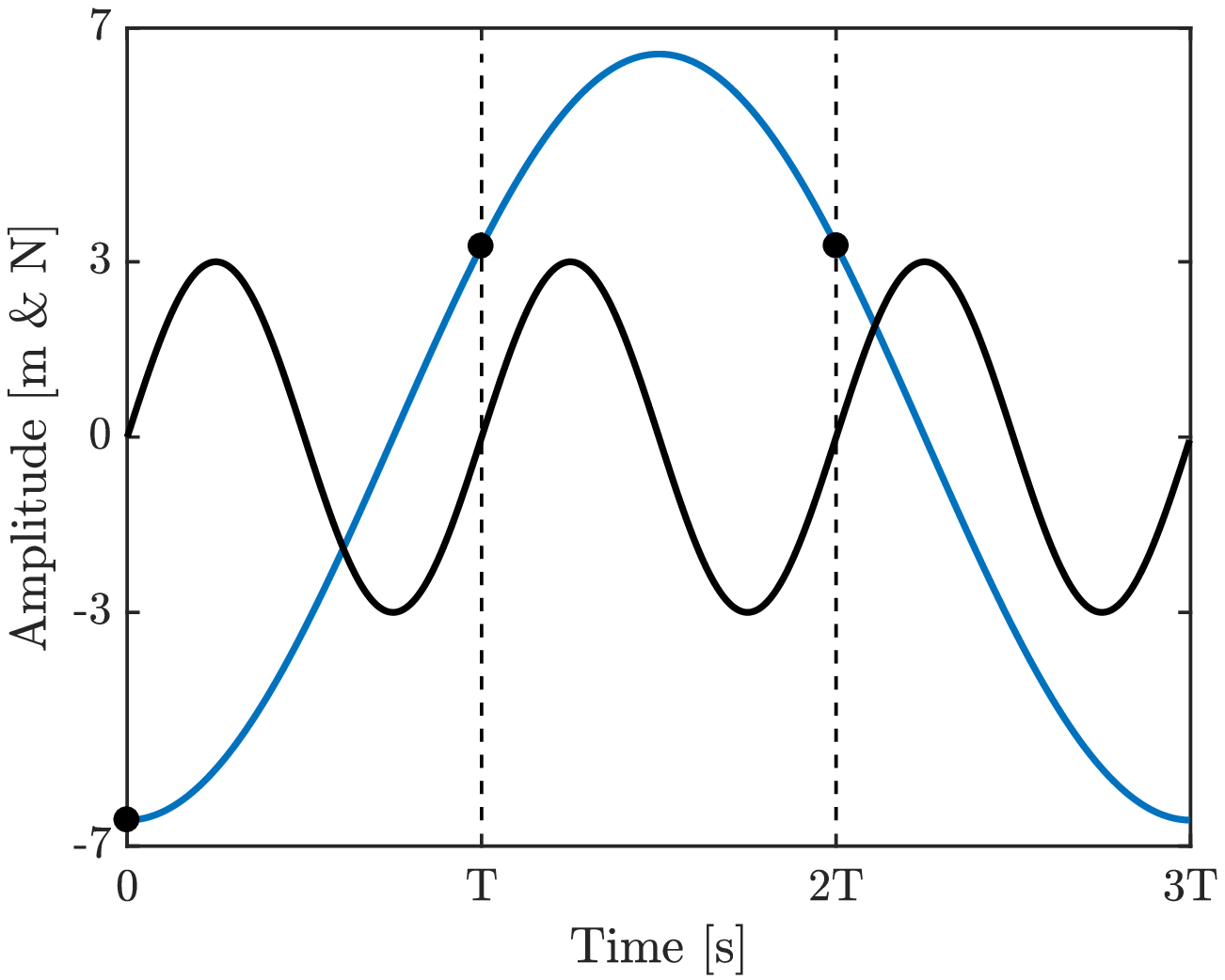} 
    \caption{\label{fig:HAR1_3_FORCING_90}}
  \end{subfigure}
  \begin{subfigure}[b]{0.5\linewidth}
    \centering
    \includegraphics[width=1\linewidth]{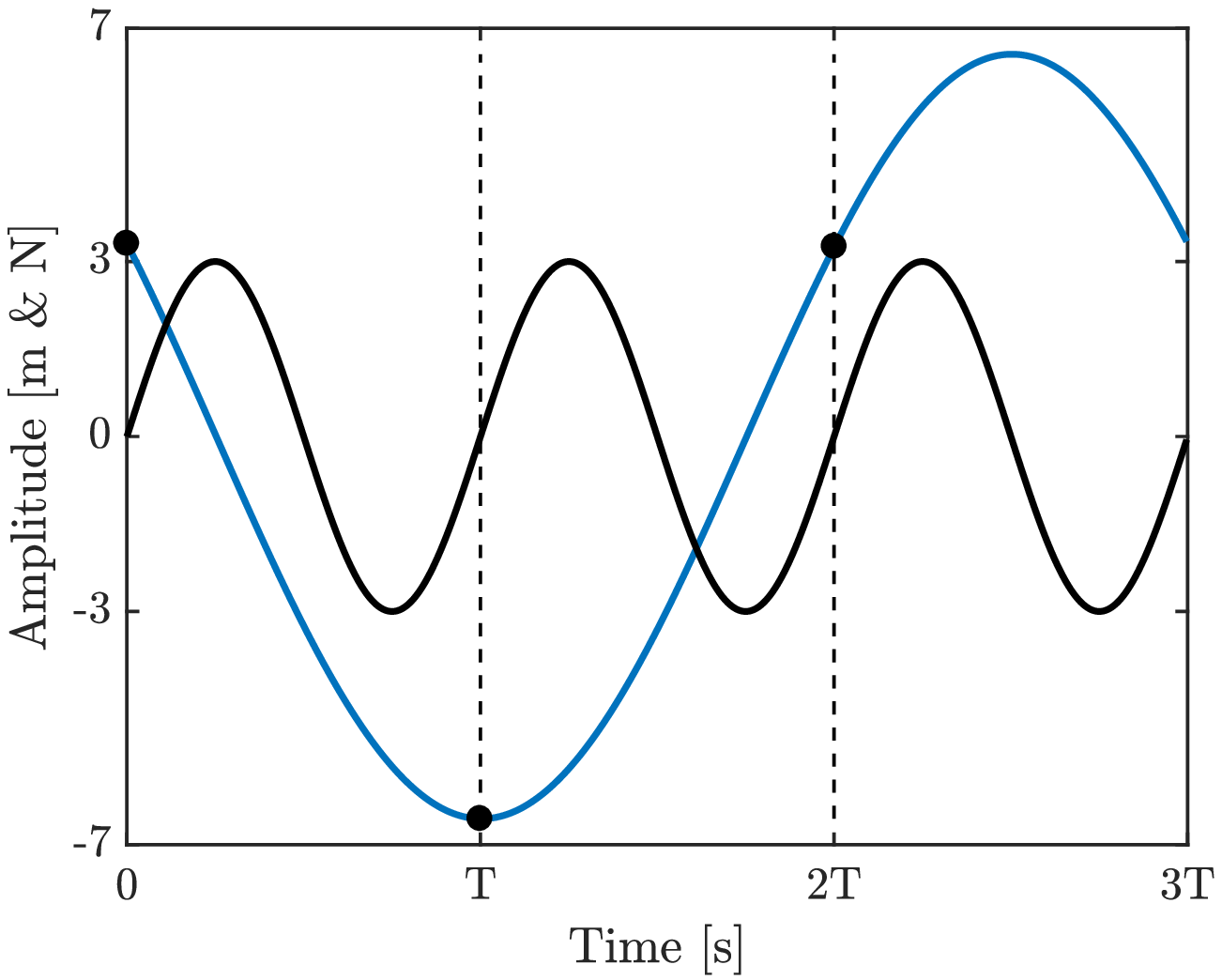}
    \caption{\label{fig:HAR1_3_FORCING_210}}
  \end{subfigure} 
  \begin{subfigure}[b]{0.5\linewidth}
    \centering
    \includegraphics[width=1\linewidth]{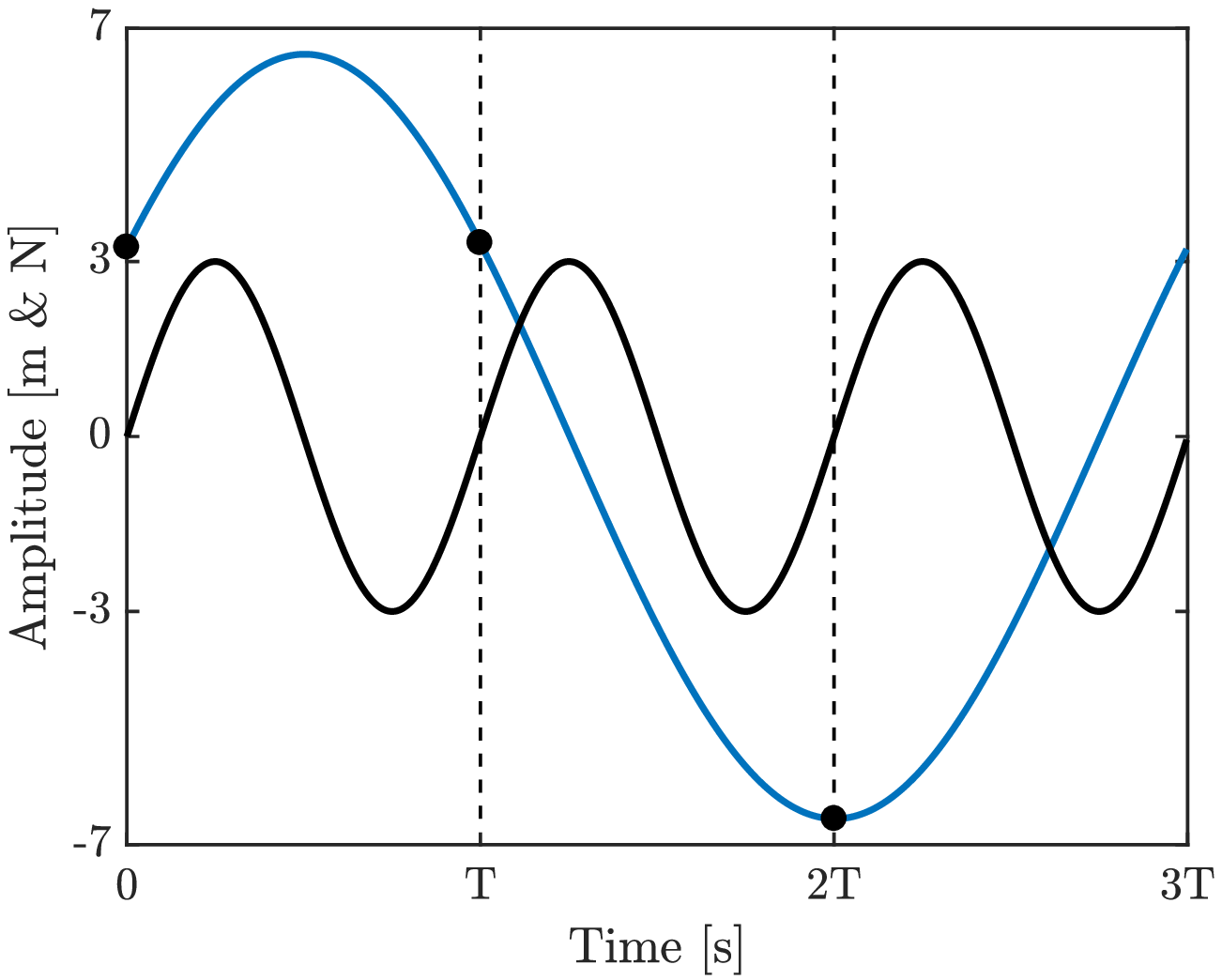}
    \caption{\label{fig:HAR1_3_FORCING_330}}
  \end{subfigure}
  \caption{Time series corresponding to the 1:3 subharmonic resonance: forcing ($f=3$N, black) and first harmonic of the displacement ($k=1$, blue): (\subref{fig:HAR1_3_FORCING_90}) $\phi_1=\frac{\pi}{2}$, (\subref{fig:HAR1_3_FORCING_210}) $\phi_1=\frac{5\pi}{6}$ and (\subref{fig:HAR1_3_FORCING_330}) $\phi_1=\frac{3\pi}{2}$.}
  \label{fig:phase_examples} 
\end{figure}

With no loss of generality, if 
\begin{equation}
    f(t)=f\sin{\omega t}, \quad x_k(t) = A_k \sin \left(\omega_k t-\phi_k\right)
    \label{eq:FourierDispPhase}
\end{equation}
where $A_k=\sqrt{s_k^2+c_k^2}$, $\phi_k=\atantwo(-c_k,s_k)$ and $f>0$, then the phase lag between the $k$-th harmonic of the displacement and the forcing corresponds to $\phi_k$. Though the $\atantwo$ function is defined in the $]-\pi;\pi]$ interval, we shift it herein in the $[0;2\pi[$ interval by adding $2\pi$ to negative values of $\phi_k$. For illustration, Figure \ref{fig:HAR3_FORCING} shows time series corresponding to the 3:1 superharmonic resonance where $f(t)=3\sin{\omega t}$ and $x_3(t) = -A_3\cos {3\omega t}$. In this case, $\phi_k$ is obviously equal to $\pi/2$, and there is therefore quadrature between the forcing and the third harmonic of the displacement. When $\nu>1$, the phase lag is not uniquely defined, as illustrated for the 1:3 subharmonic resonance in Figure \ref{fig:phase_examples}. In this case, the phase lag can take three values, namely $\phi_1=\frac{\pi}{2}$, $\phi_1=\frac{5\pi}{6}$ or $\phi_1=\frac{3\pi}{2}$. 


Although sine functions are considered throughout this study, we remark that the phase lags for cosine functions 
$f(t)=f\cos{\omega t}, \quad x_k(t) = A_k \cos \left(\omega_kt-\phi_k\right)$ should be calculated with $\phi_k=\atantwo(s_k,c_k)$, as in \cite{NAYFEH,LEUNG}.

\subsection{Primary resonance (1:1)}

\begin{figure}[ht] 
  \begin{subfigure}[b]{0.5\linewidth}
    \centering
    \includegraphics[width=1\linewidth]{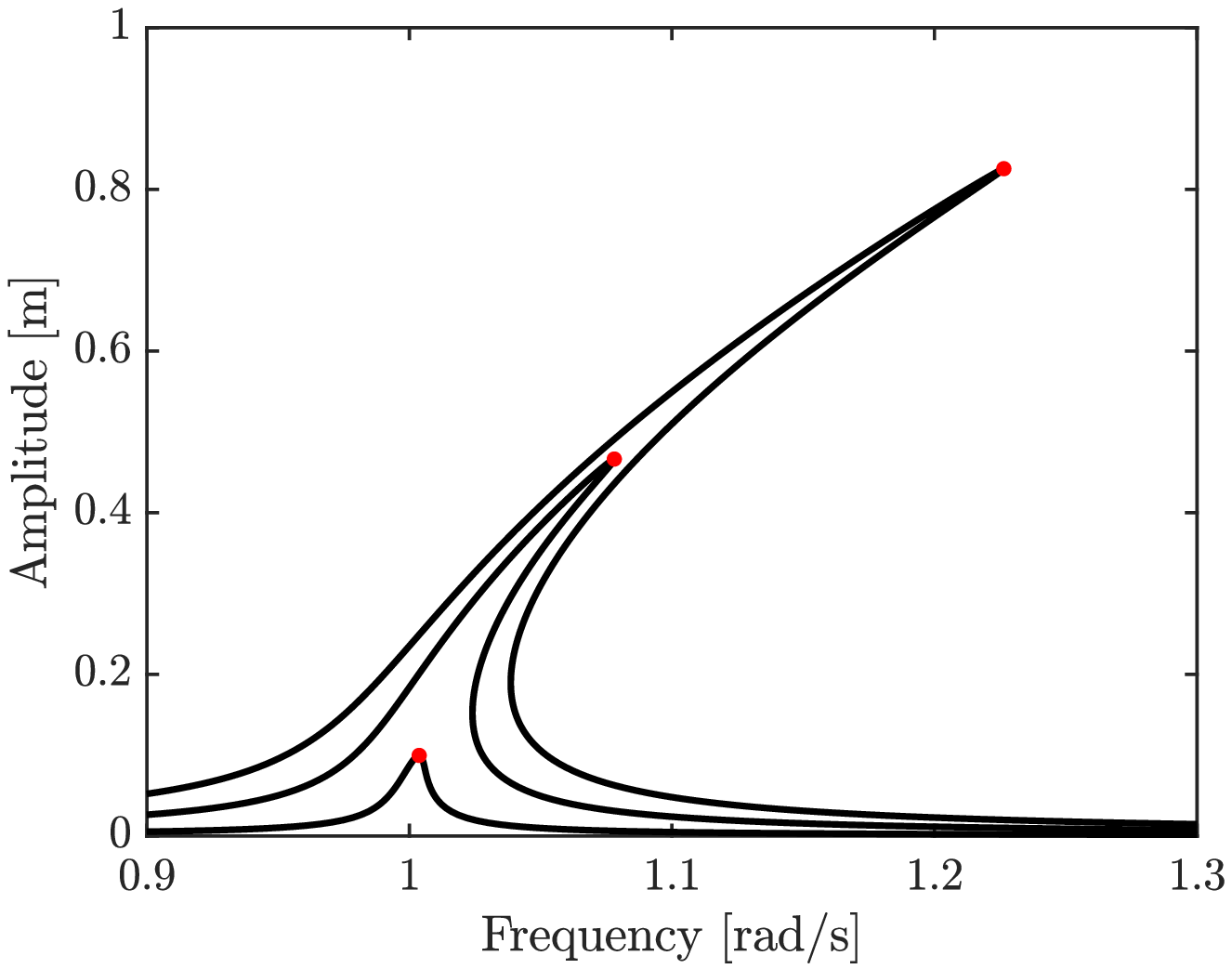} 
    \caption{\label{fig:FUNDRED_NFRC_NO_PRNM}}
  \end{subfigure}
  \begin{subfigure}[b]{0.5\linewidth}
    \centering
    \includegraphics[width=1\linewidth]{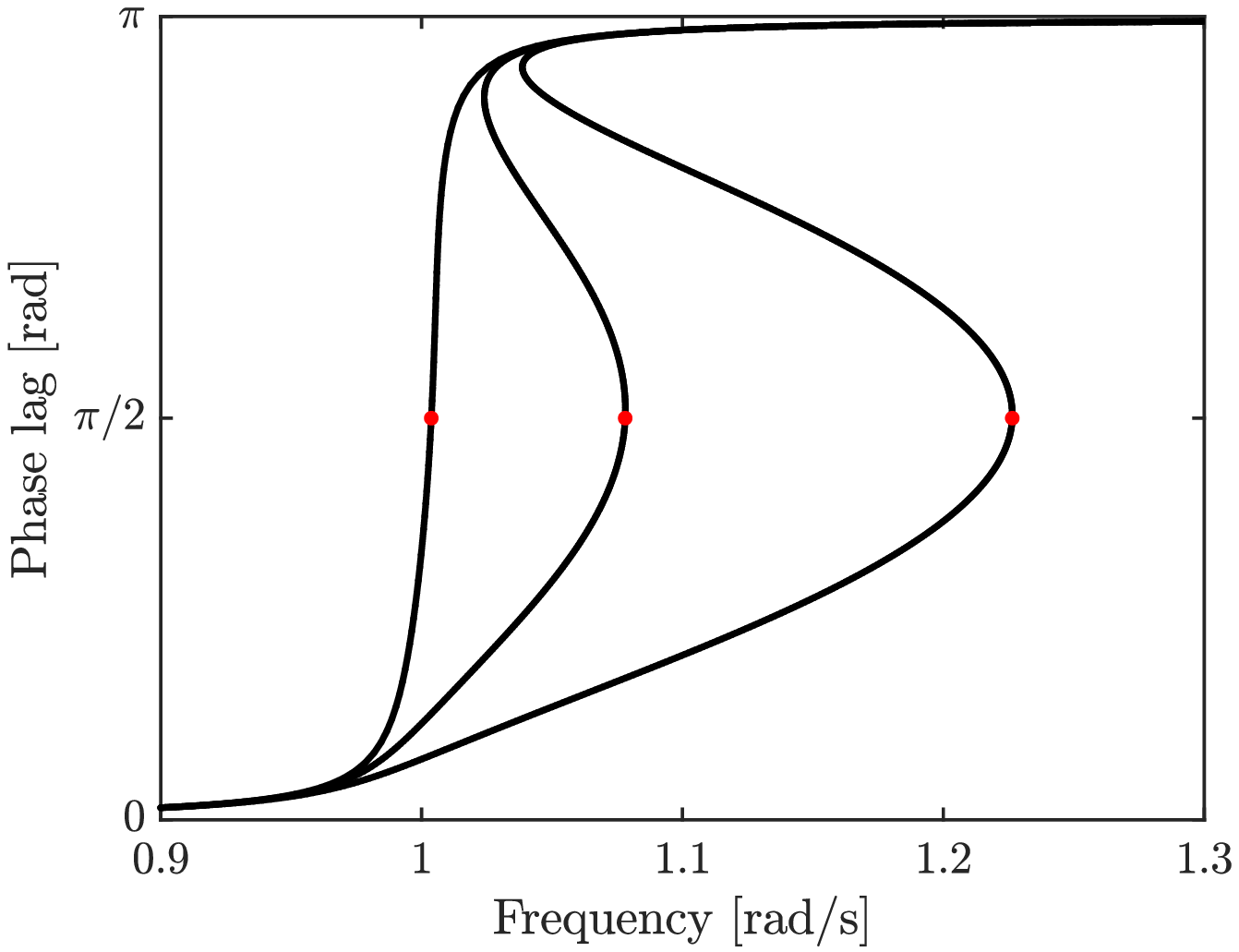}
    \caption{\label{fig:FUNDRES_PHASELAG_NOPRNM}}
  \end{subfigure} 
  \caption{NFRCs of the primary resonance of the Duffing oscillator for 3 forcing amplitudes (0.001N, 0.005N and 0.01N): (\subref{fig:FUNDRED_NFRC_NO_PRNM}) amplitude and (\subref{fig:FUNDRES_PHASELAG_NOPRNM}) phase lag of the first harmonic component. The red dots correspond to phase quadrature.}
  \label{fig:FUNDRES} 
\end{figure}

The amplitude and phase lag of the first harmonic of the displacement in the neighborhood of the primary resonance are displayed in Figures \ref{fig:FUNDRED_NFRC_NO_PRNM} and \ref{fig:FUNDRES_PHASELAG_NOPRNM}, respectively. These figures show that the phase lag varies between $0$ and $\pi$ and that nonlinear phase resonance corresponding to a phase lag of $\pi/2$ occurs very near amplitude resonance, at least for the amount of damping considered herein. 

\subsection{Superharmonic resonances (\texorpdfstring{$k:1$}{})}
In this case, the resonances are located for values of $\omega$ lower than $\omega_0$. 

\subsubsection{Odd superharmonic resonances}
\begin{figure}[ht] 
  \begin{subfigure}[b]{0.5\linewidth}
    \centering
    \includegraphics[width=1\linewidth]{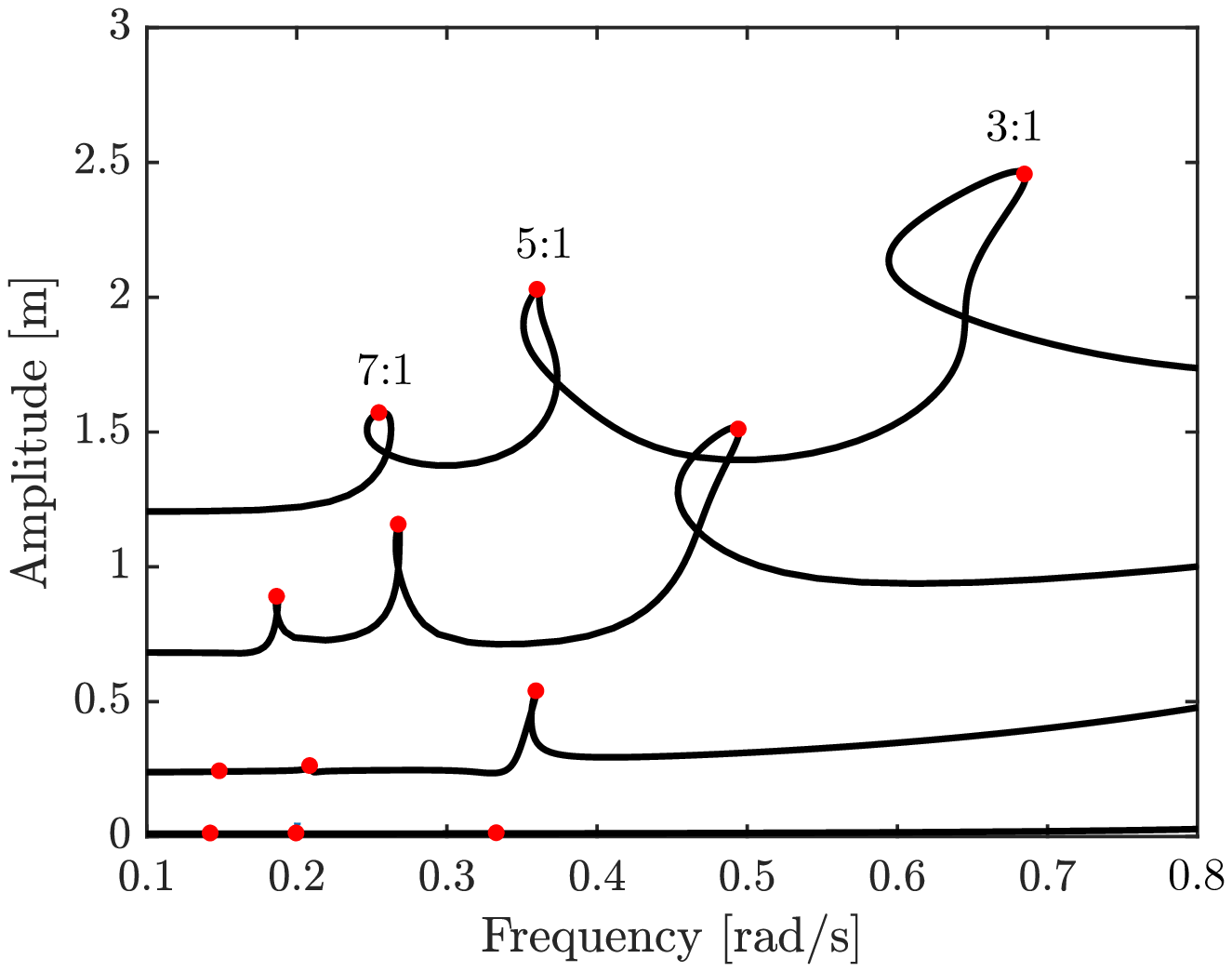} 
    \caption{\label{fig:SUPERHODD_NFRC_NO_PRNM}}
  \end{subfigure}
  \begin{subfigure}[b]{0.5\linewidth}
    \centering
    \includegraphics[width=1\linewidth]{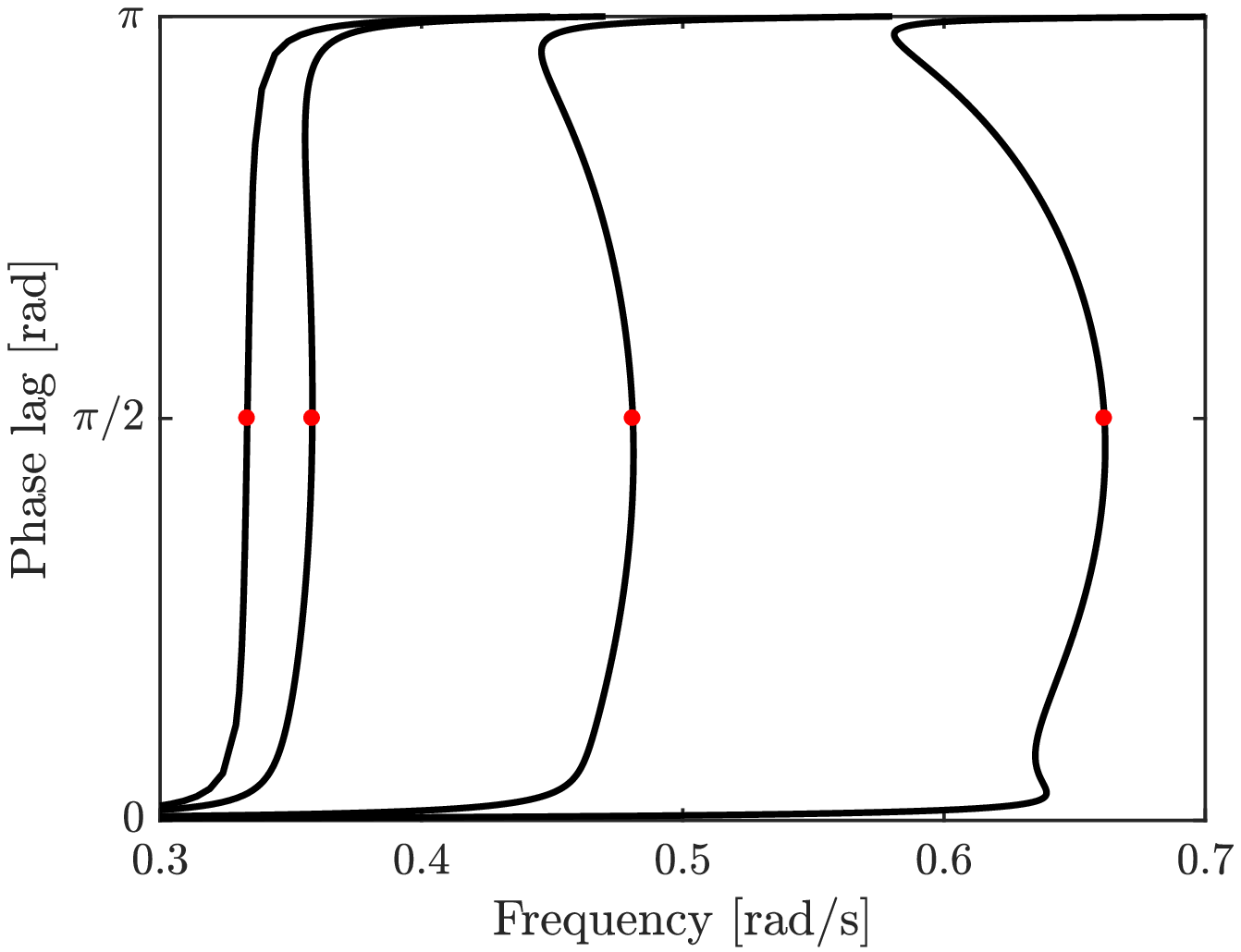}
    \caption{\label{fig:SUPERHODD_PHASELAG_NOPRNM}}
  \end{subfigure} 
  \caption{NFRCs of the 3:1, 5:1 and 7:1 superharmonic resonances for 4 forcing amplitudes (0.01N, 0.25N, 1N and 3N): (\subref{fig:SUPERHODD_NFRC_NO_PRNM}) amplitude and (\subref{fig:SUPERHODD_PHASELAG_NOPRNM}) phase lag of the $3^{rd}$ harmonic component of the 3:1 resonance (computed with only 3 harmonics for better readability). The points where the phase lag is equal to $\pi/2$ are marked by red dots. }
  \label{fig:SUPHERHODD} 
\end{figure}

Odd superharmonic resonances are easily calculated, because they appear in the direct continuation of the main branch of the NFRC. Three superharmonic resonances, namely 3:1, 5:1 and 7:1, are represented in Figure \ref{fig:SUPERHODD_NFRC_NO_PRNM}. More resonances could have been found if more harmonics were considered in the HBM. 

For the 3:1 resonance, the phase lag is represented in Figure \ref{fig:SUPERHODD_PHASELAG_NOPRNM}. A similar behavior as for the primary resonance is observed, namely the phase lag evolves from 0 to $\pi$, and the phase quadrature points are in excellent agreement with the maxima of the frequency response. The same finding holds for the other superharmonic resonances. These results suggest that phase quadrature between the $k$-th harmonic of the displacement and the forcing exists at resonance for odd superharmonic branches, which is fully consistent with the results in \cite{LEUNG}.


\subsubsection{Even superharmonic resonances}

\begin{figure}[ht] 
  \begin{subfigure}[b]{0.5\linewidth}
    \centering
    \includegraphics[width=1\linewidth]{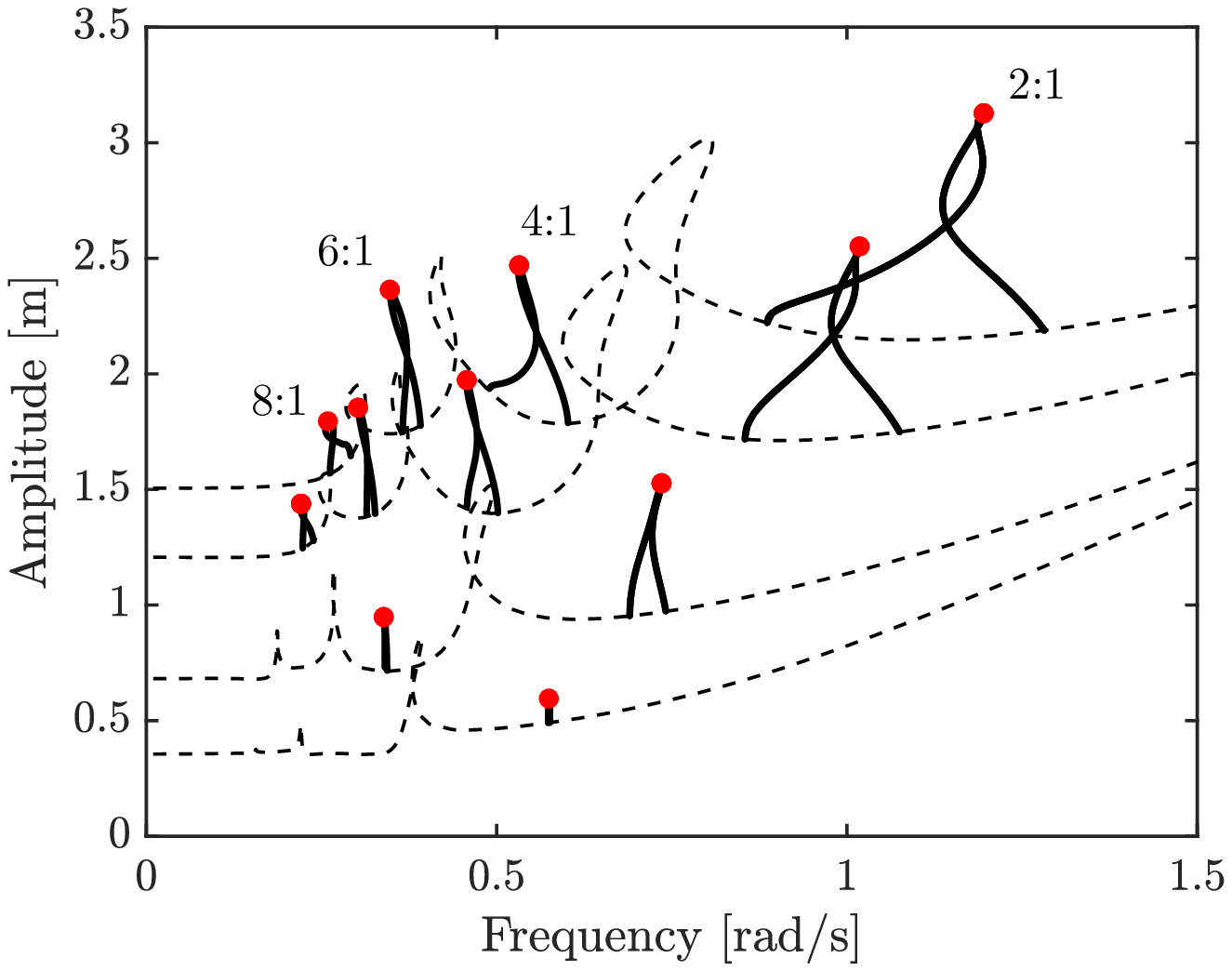}
    \caption{\label{fig:SUPERHEVEN_NFRC_NO_PRNM}}
  \end{subfigure}
  \begin{subfigure}[b]{0.5\linewidth}
    \centering
    \includegraphics[width=1\linewidth]{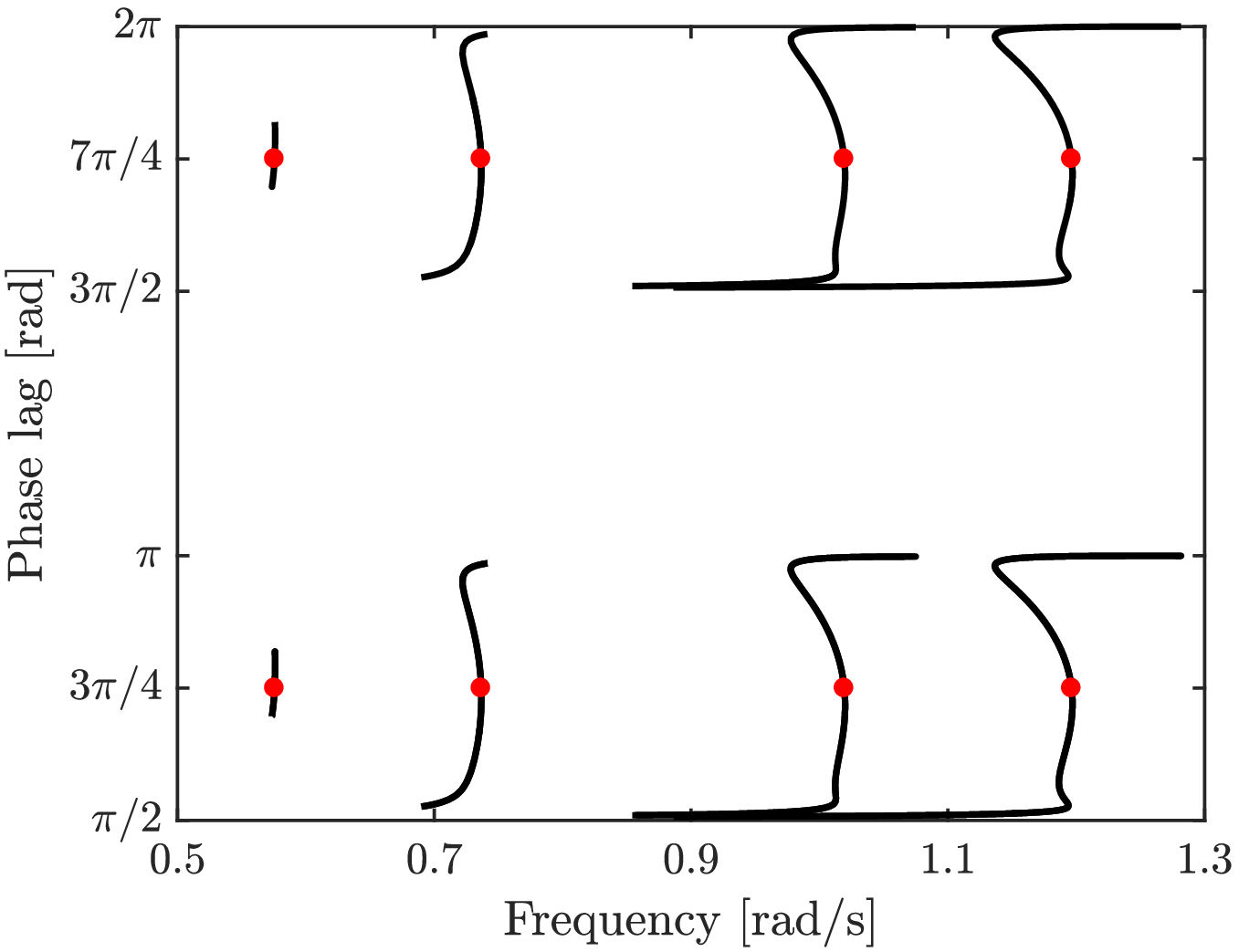}
    \caption{\label{fig:SUPERHEVEN_PHASELAG_NO_PRNM}}
  \end{subfigure} 
  \caption{NFRCs of the 2:1, 4:1, 6:1 and 8:1 superharmonic resonances for 4 forcing amplitudes (0.4N, 1N, 3N and 5N): (\subref{fig:SUPERHEVEN_NFRC_NO_PRNM}) amplitude and (\subref{fig:SUPERHEVEN_PHASELAG_NO_PRNM}) phase lag of the $2^{nd}$ harmonic component of the 2:1 resonance. The points where the phase lag is equal to $\frac{3\pi}{4}$ and $\frac{7\pi}{4}$ are marked by red dots.}
  \label{fig:SUPHERHEVEN} 
\end{figure}

Contrary to odd superharmonic resonances, the 2:1, 4:1, 6:1 and 8:1 superharmonic resonances in Figure \ref{fig:SUPERHEVEN_NFRC_NO_PRNM} bifurcate out from the main NFRC. The phase lag of the 2:1 resonance in Figure \ref{fig:SUPERHEVEN_PHASELAG_NO_PRNM} shows a fundamental difference with the phase lag of odd resonances, i.e., phase quadrature is no longer observed in the vicinity of the resonance. Specifically, the phase lag, comprised between $\pi/2$ and $\pi$, passes through $3\pi/4$ at resonance. The same observation holds for the 4:1, 6:1 and 8:1 resonances. There is thus the need to generalize the phase resonance condition for such resonances. 

\begin{figure}[ht] 
  \begin{subfigure}[b]{0.5\linewidth}
    \centering
    \includegraphics[width=1\linewidth]{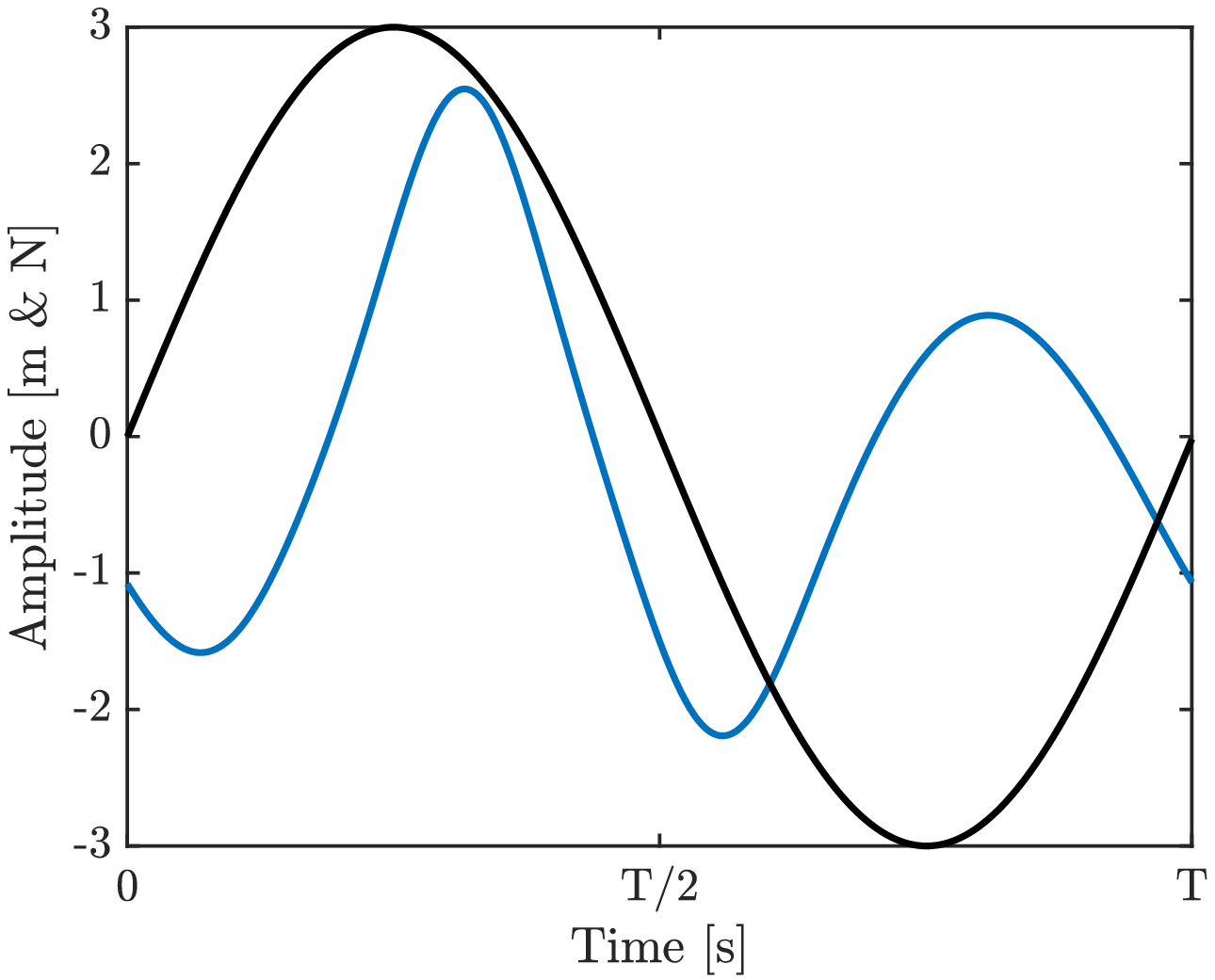}
    \caption{\label{fig:SUPERHEVEN_TIME_SERIES_1}}
  \end{subfigure}
  \begin{subfigure}[b]{0.5\linewidth}
    \centering
    \includegraphics[width=1\linewidth]{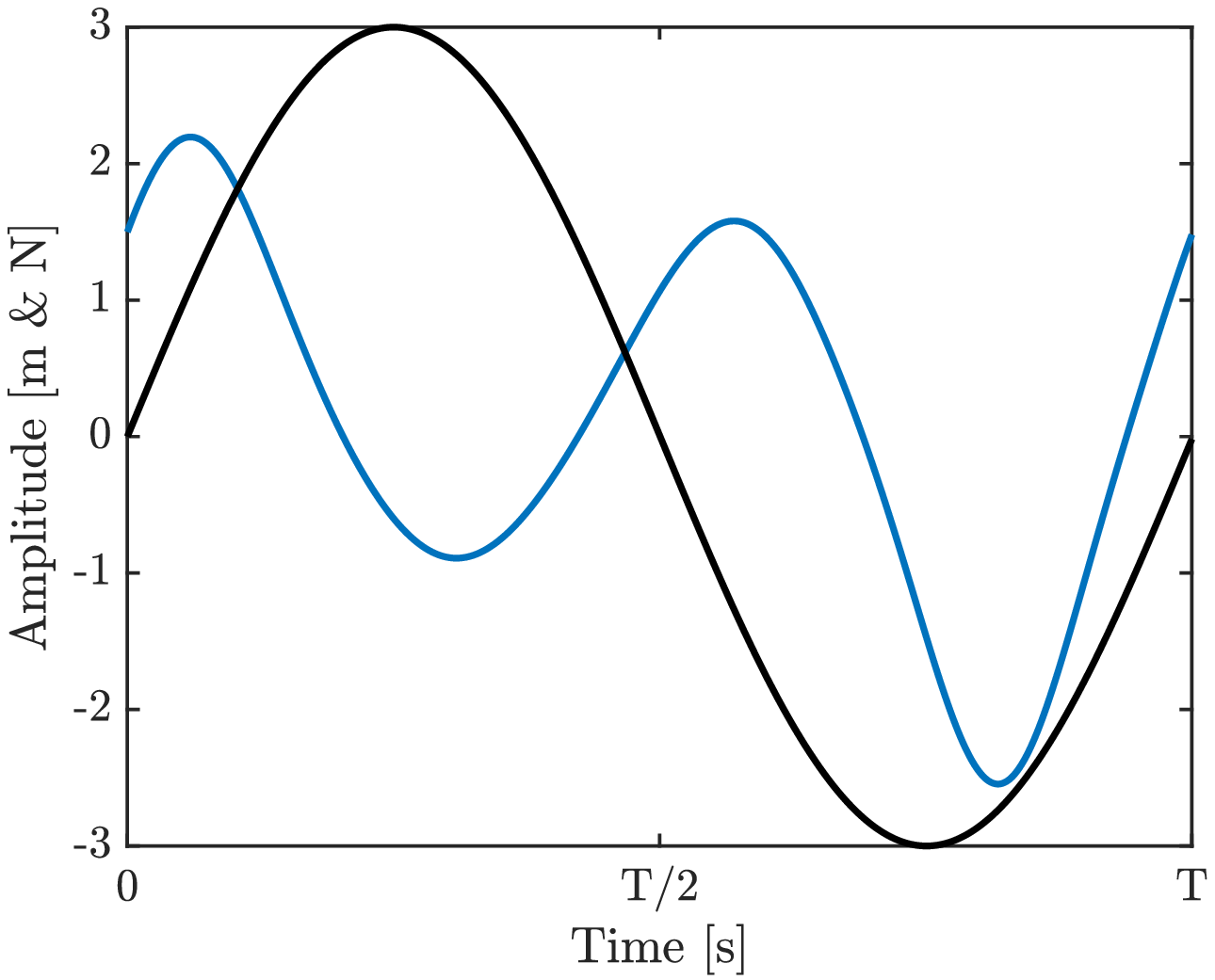}
    \caption{\label{fig:SUPERHEVEN_TIME_SERIES_2}}
  \end{subfigure} 
  \caption{Time series (blue) and forcing (black) corresponding to two different branches of the 2:1 superharmonic resonance: (\subref{fig:SUPERHEVEN_TIME_SERIES_1}) phase lag = $3\pi/4$ and (\subref{fig:SUPERHEVEN_TIME_SERIES_2}) phase lag = $7\pi/4$.}
  \label{fig:SUPERHEVEN_TIME_SERIES} 
\end{figure}

We also underline that, for each resonance, the continuation procedure goes two times through the branch, with a phase shift of $\pi$ between the two solutions. This phase shift is clearly seen in Figure \ref{fig:SUPERHEVEN_PHASELAG_NO_PRNM}. The corresponding time series in Figure \ref{fig:SUPERHEVEN_TIME_SERIES} show that the solutions are not symmetric with respect to the mid period and are opposite of each other with a phase shift of $\pi$. They share, however, the same maximum amplitude in absolute value. This behavior is typical of branches generated through symmetry-breaking bifurcations.

\subsection{Subharmonic resonances (\texorpdfstring{$1:\nu$}{})}
In this case, the resonances are located for values of $\omega$ greater than $\omega_0$. 

\subsubsection{Odd subharmonic resonances}
\begin{figure}[ht] 
  \begin{subfigure}[b]{0.5\linewidth}
    \centering
    \includegraphics[width=1\linewidth]{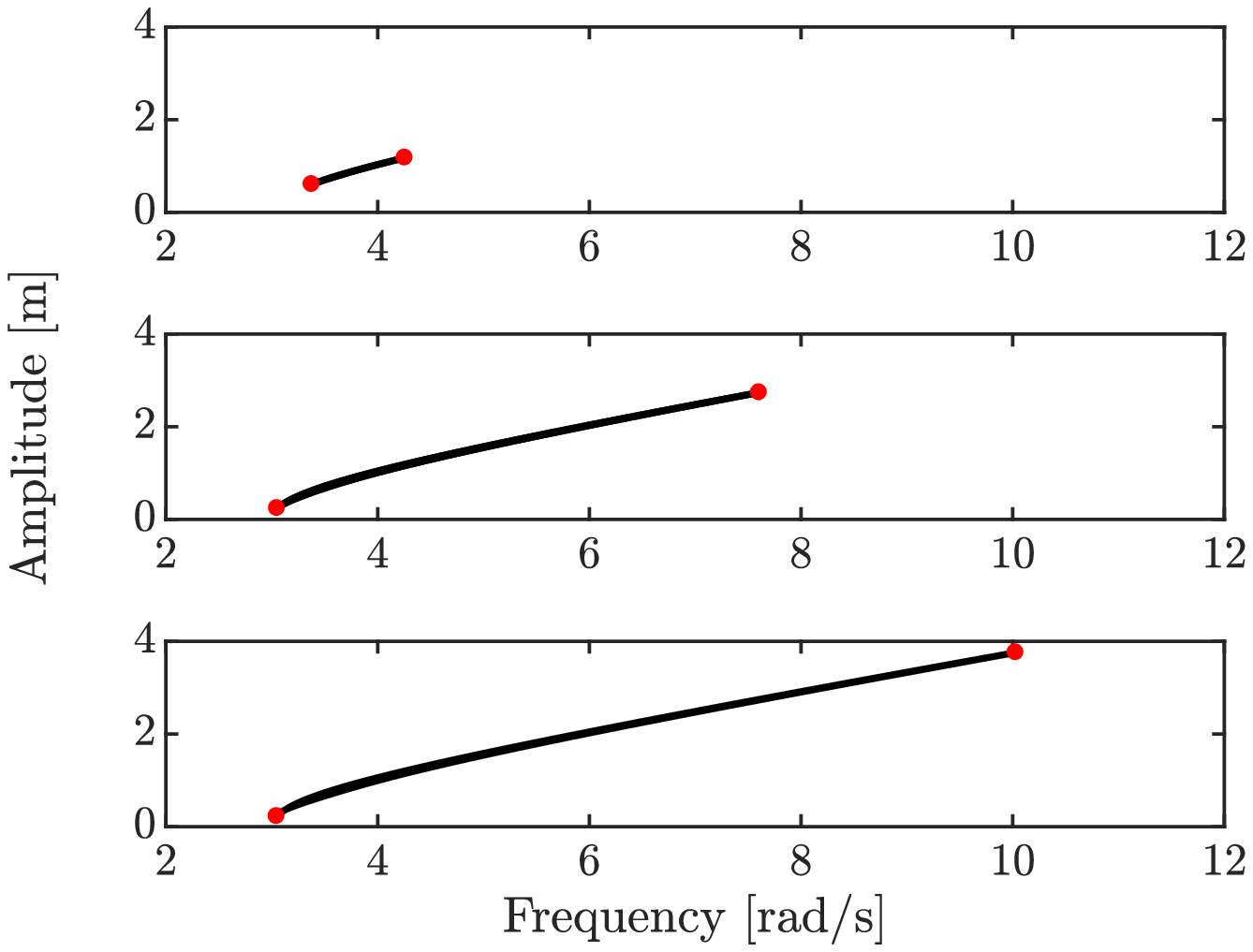}
    \caption{\label{fig:SUBHODD_NFRC_NO_PRNM}}
  \end{subfigure}
  \begin{subfigure}[b]{0.5\linewidth}
    \centering
    \includegraphics[width=1\linewidth]{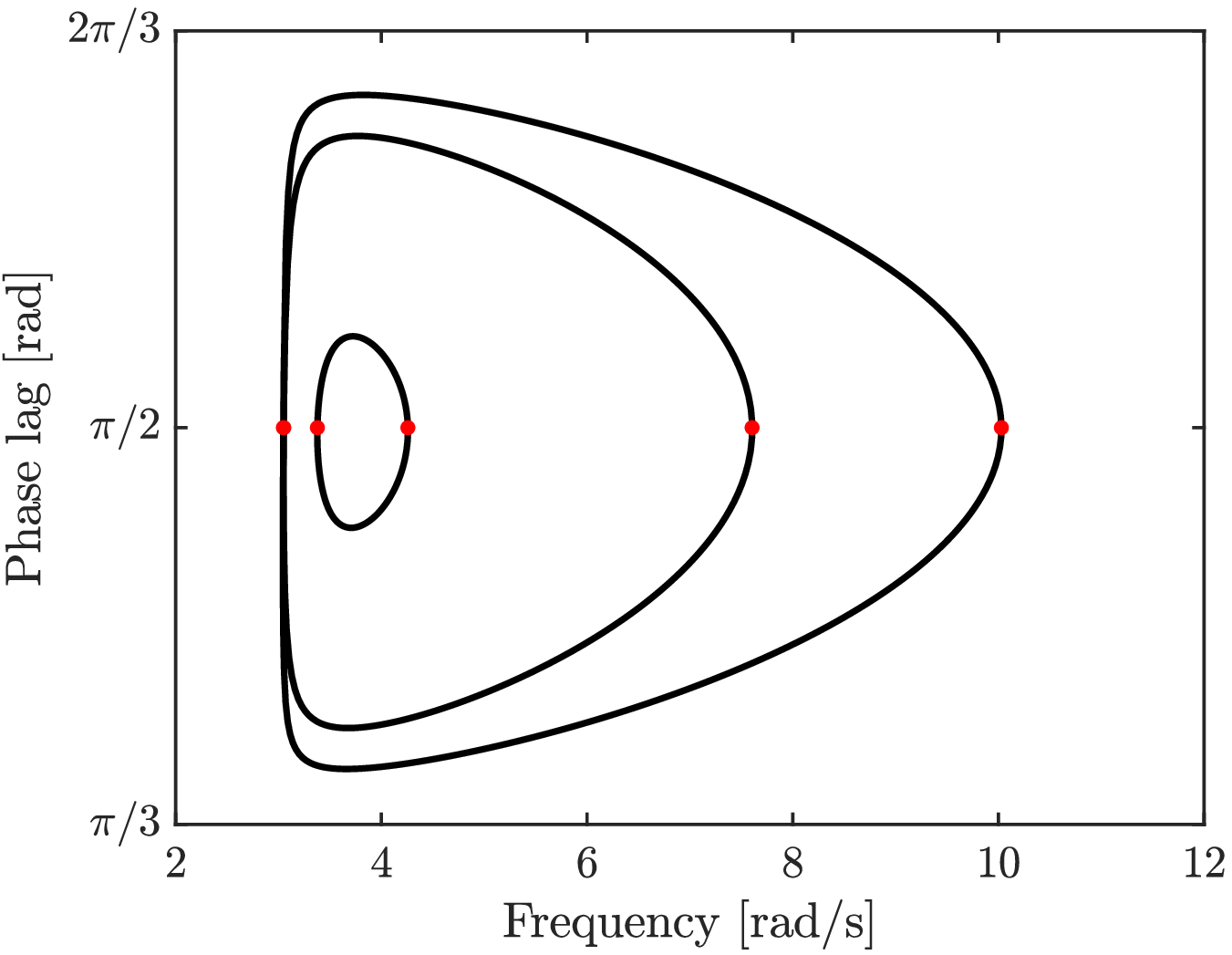}
    \caption{\label{fig:SUBHODD_PHASE_NO_PRNM}}
  \end{subfigure} 
  \caption{NFRCs of the 1:3 subharmonic resonance for 3 forcing amplitudes (0.25N, 0.6N and 1N): (\subref{fig:SUBHODD_NFRC_NO_PRNM}) amplitude and (\subref{fig:SUBHODD_PHASE_NO_PRNM}) phase lag of the first harmonic.}
  \label{fig:SUBHODD} 
\end{figure}

The 1:3 subharmonic resonance observed in Figure \ref{fig:F_0_25} is displayed in Figure \ref{fig:SUBHODD}. The phase lag, bounded by $\pi/3$ and $2\pi/3$, goes through $\pi/2$ at the extremities of the isola. Phase quadrature is thus observed twice, corresponding to the minimum and maximum response amplitudes on the branch, in accordance with \cite{LEUNG}. This is also the case for higher-order subharmonic resonances, but the phase lag is located in the interval $[\pi/2\pm\pi/2\nu]$.

\subsubsection{Even subharmonic resonances}
\begin{figure}[ht] 
  \begin{subfigure}[b]{0.5\linewidth}
    \centering
    \includegraphics[width=1\linewidth]{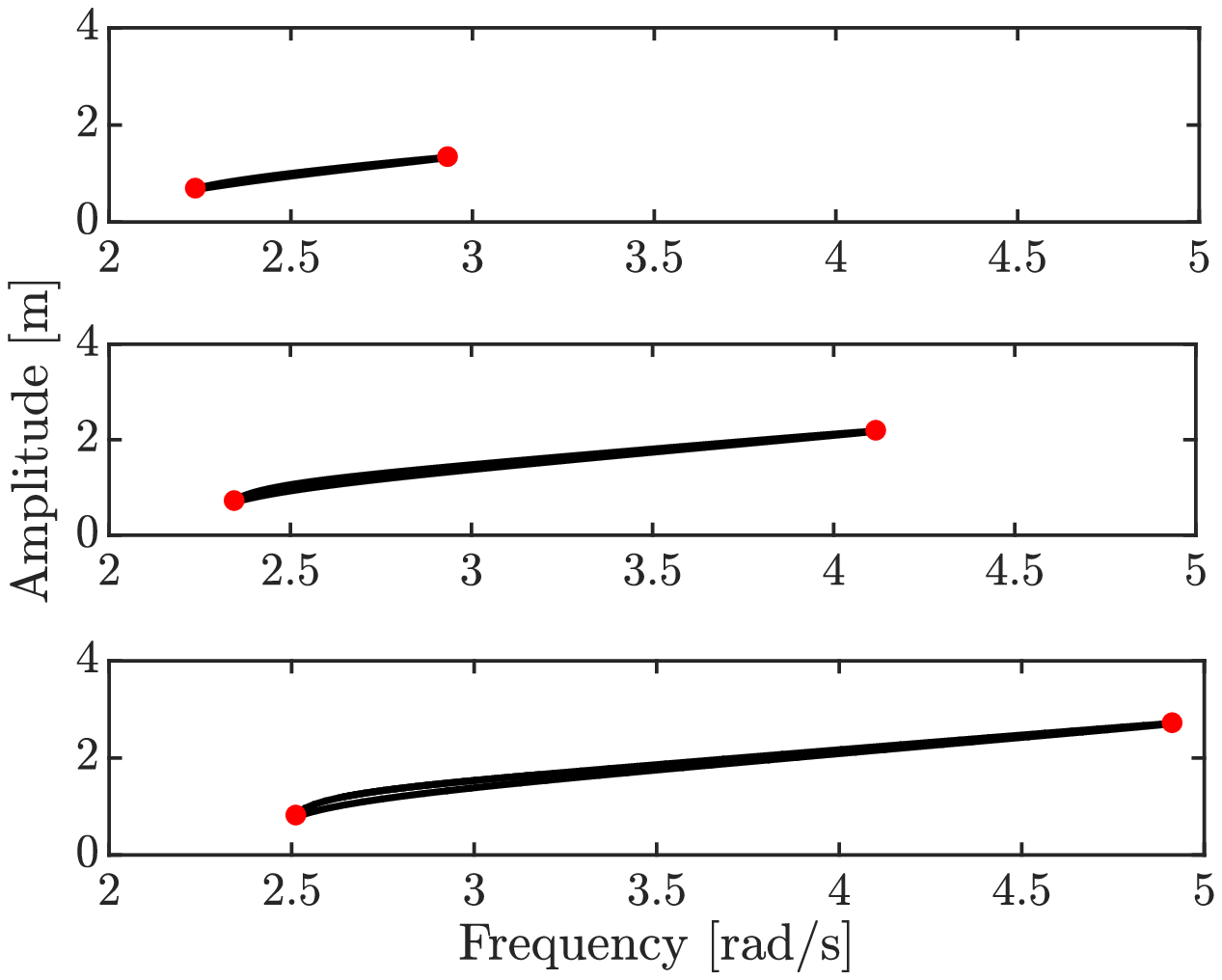}
    \caption{\label{fig:SUBHEVEN_NFRC_NO_PRNM}}
  \end{subfigure}
  \begin{subfigure}[b]{0.5\linewidth}
    \centering
    \includegraphics[width=1\linewidth]{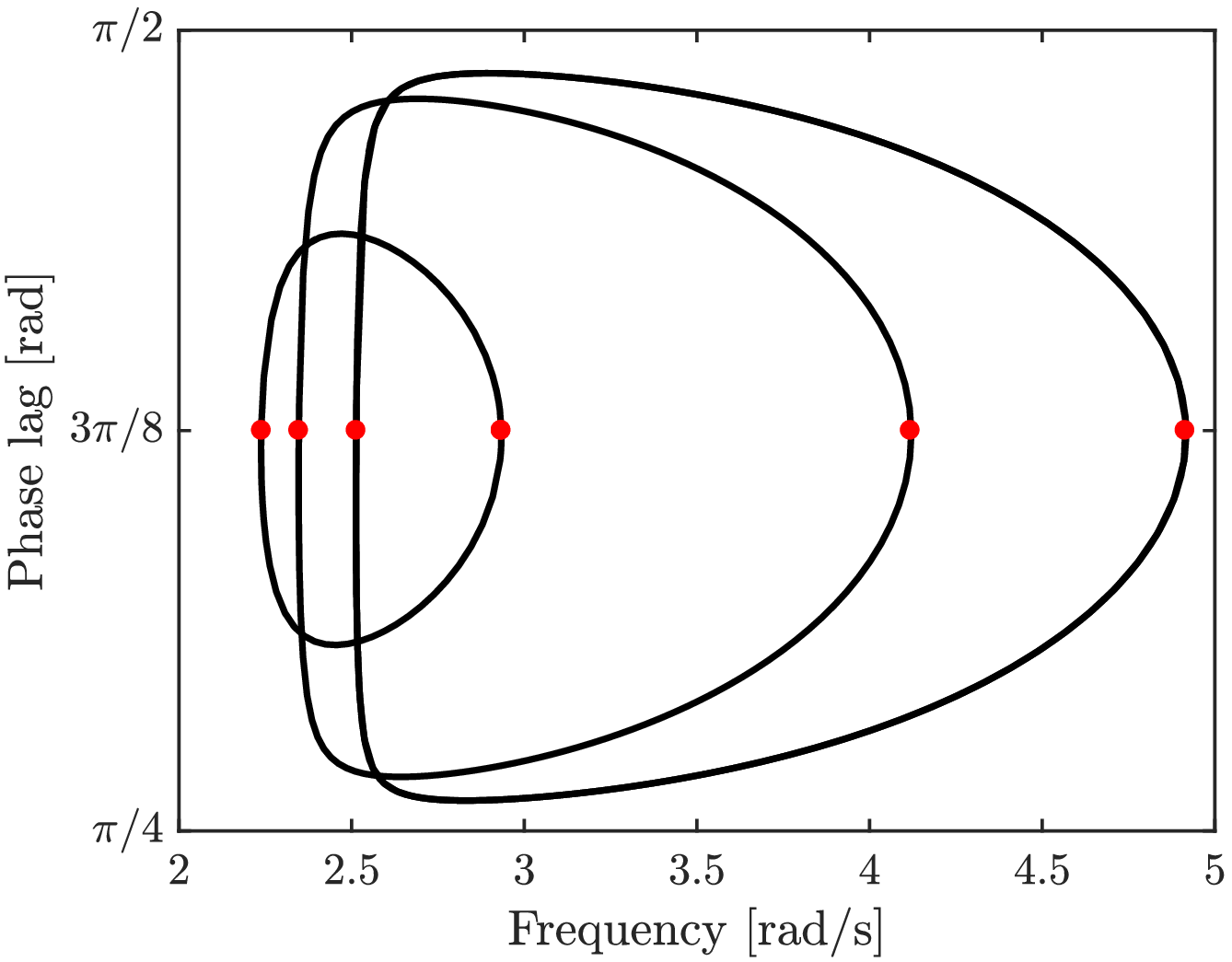}
    \caption{\label{fig:SUBHEVEN_PHASE_NO_PRNM}}
  \end{subfigure} 
  \caption{NFRCs of the 1:2 subharmonic resonance for 3 forcing amplitudes (1N, 2N and 3N):  (\subref{fig:SUBHEVEN_NFRC_NO_PRNM}) amplitude and (\subref{fig:SUBHEVEN_PHASE_NO_PRNM}) phase lag of the first harmonic.}
  \label{fig:SUBHEVEN} 
\end{figure}

As for even superharmonic resonances, the phase lag is not centered around $\pi/2$. Specifically, for the 1:2 resonance in Figure \ref{fig:SUBHEVEN}, it is centered around $3\pi/8$ and comprised between $\pi/4$ and $\pi/2$. We note that a second solution, whose maximum amplitude in absolute value is the same as in Figure \ref{fig:SUBHEVEN_NFRC_NO_PRNM}  and whose phase lag is shifted by $\pi/\nu$ can be found. For higher-order resonances, the phase lag is located in the interval $[3\pi/4\nu\pm\pi/4\nu]$. 

\subsection{Ultra-subharmonic resonances (\texorpdfstring{$k:\nu$}{})}

\subsubsection{\texorpdfstring{$k$}{} and \texorpdfstring{$\nu$}{} are odd}
\begin{figure}[ht] 
  \begin{subfigure}[b]{0.5\linewidth}
    \centering
    \includegraphics[width=1\linewidth]{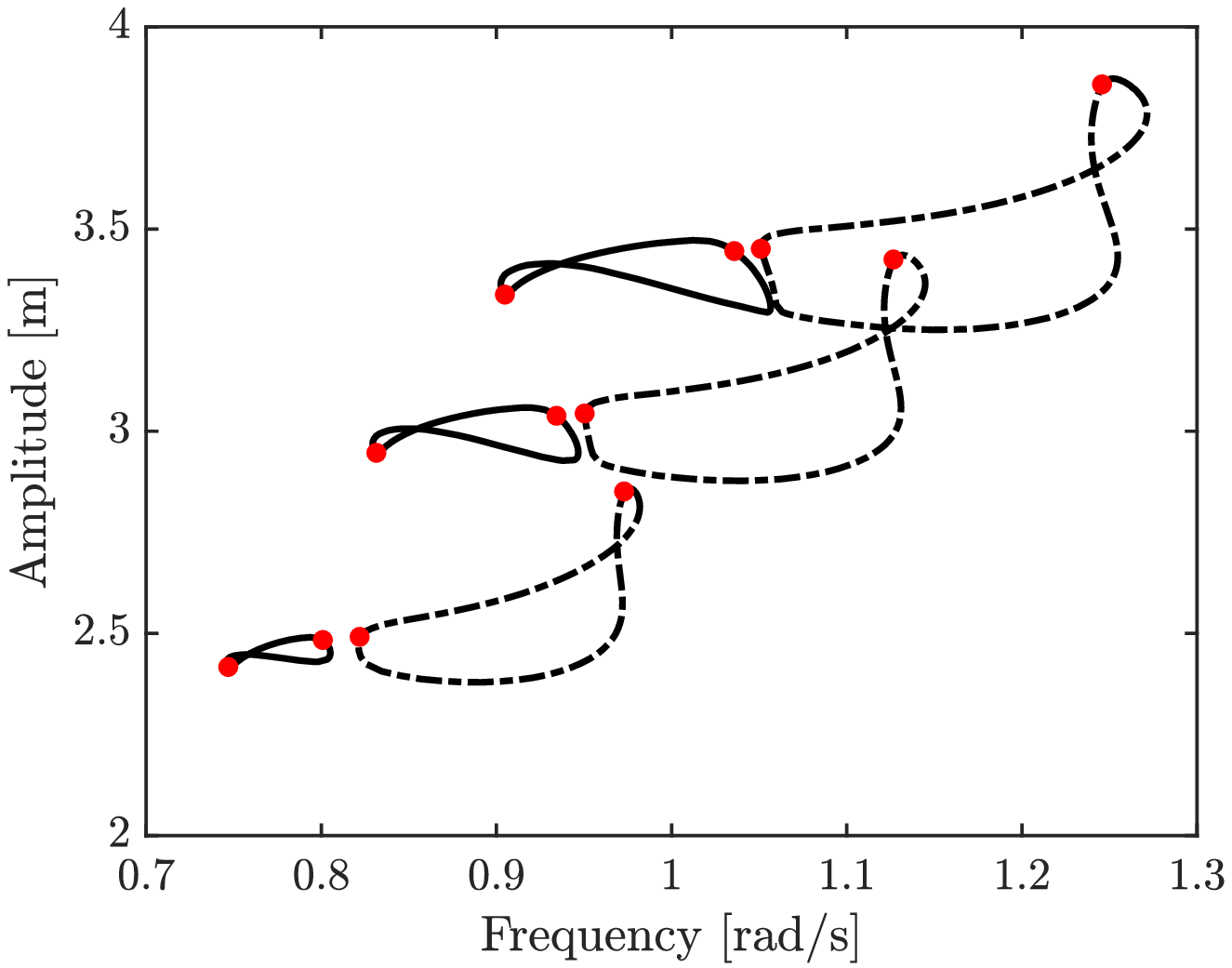} 
    \caption{\label{fig:RATIONAL_ODD_7_3_NFRC_NO_PRNM}} 
  \end{subfigure}
  \begin{subfigure}[b]{0.5\linewidth}
    \centering
    \includegraphics[width=1\linewidth]{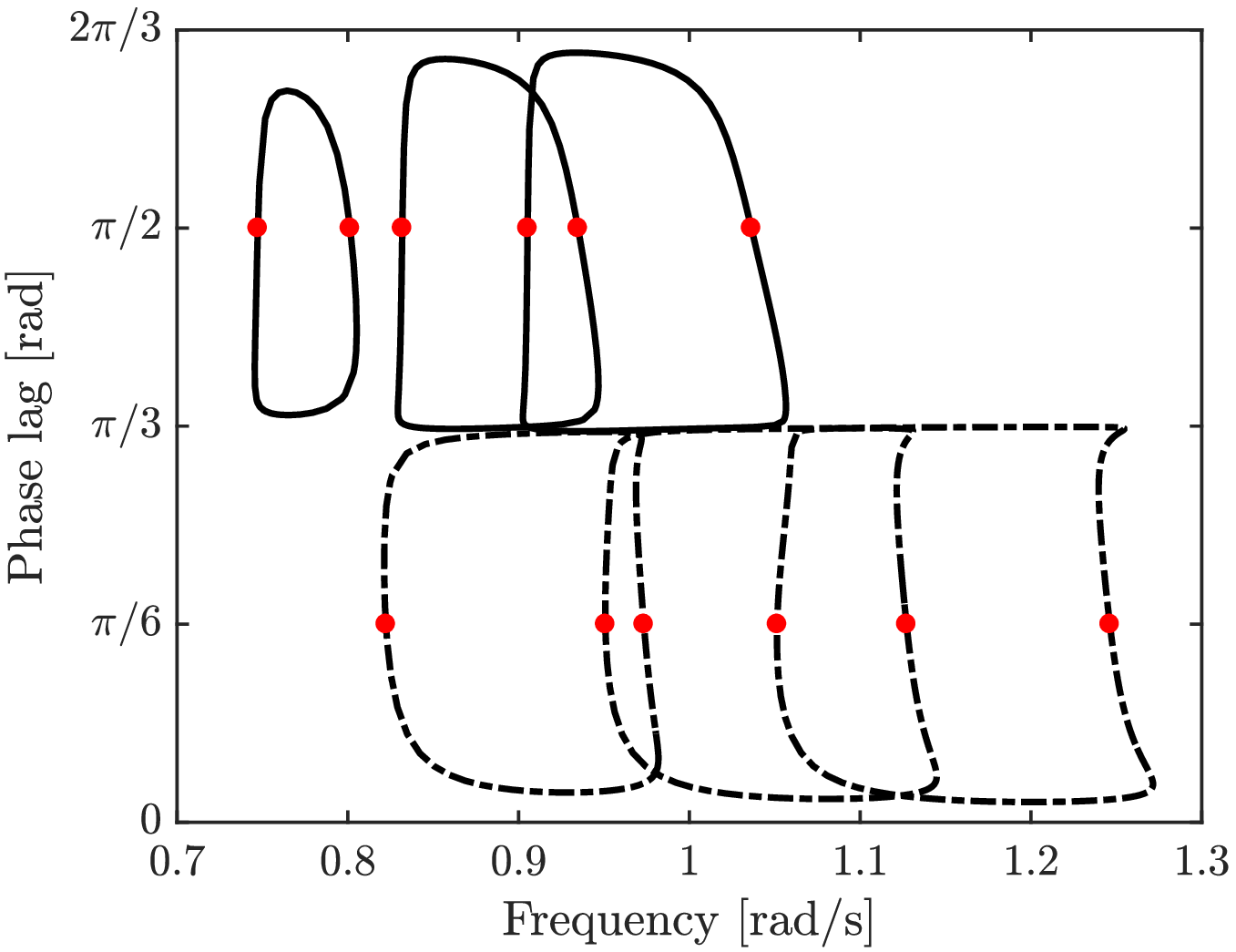} 
    \caption{\label{fig:RATIONAL_ODD_7_3_PHASE_NO_PRNM}} 
  \end{subfigure} 
  
  \begin{subfigure}[b]{0.5\linewidth}
    \centering
    \includegraphics[width=1\linewidth]{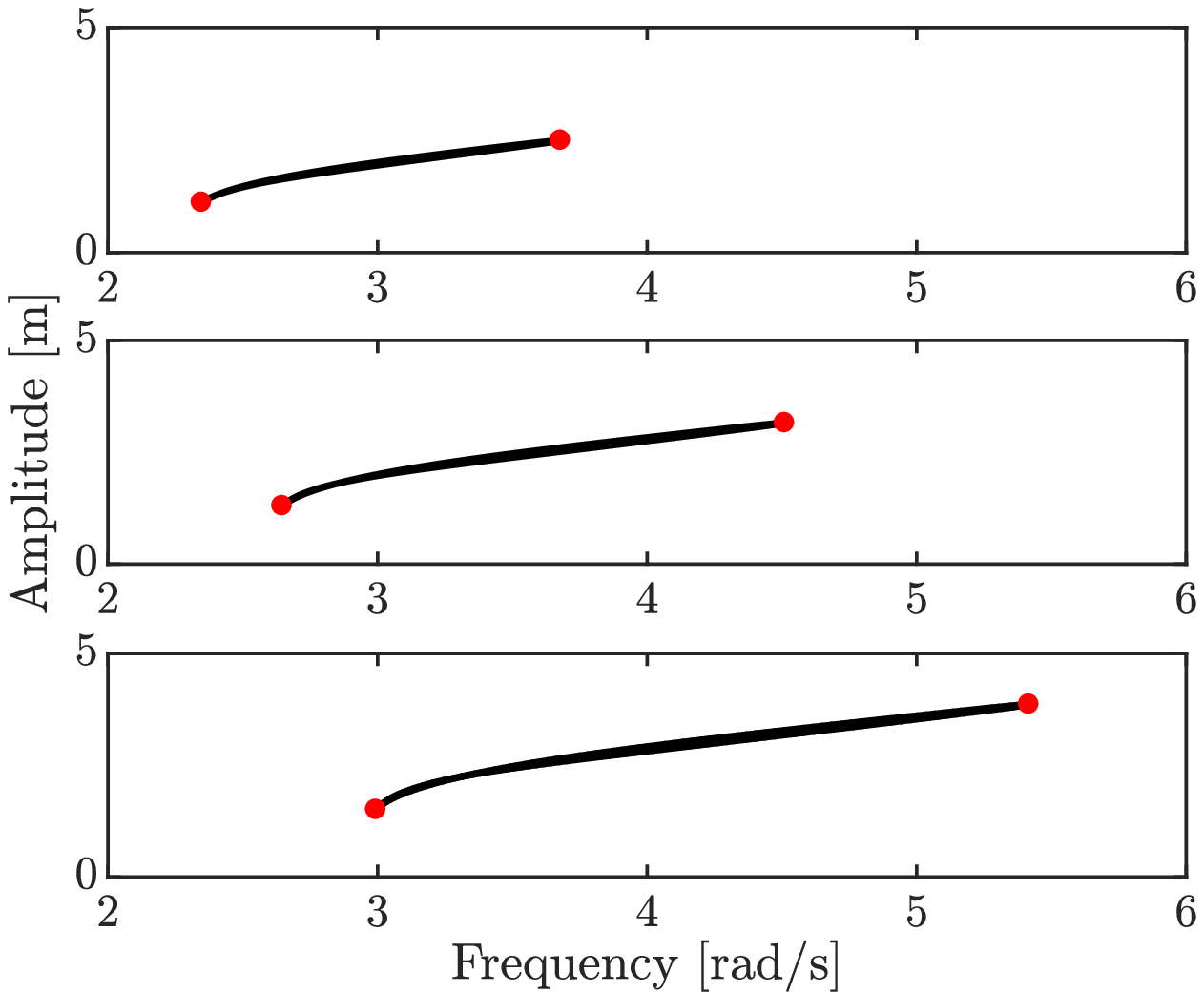}
    \caption{\label{fig:RATIONAL_ODD_SUBPLOT_3_5_NFRC_NO_PRNM}} 
  \end{subfigure}
  \begin{subfigure}[b]{0.5\linewidth}
    \centering
    \includegraphics[width=1\linewidth]{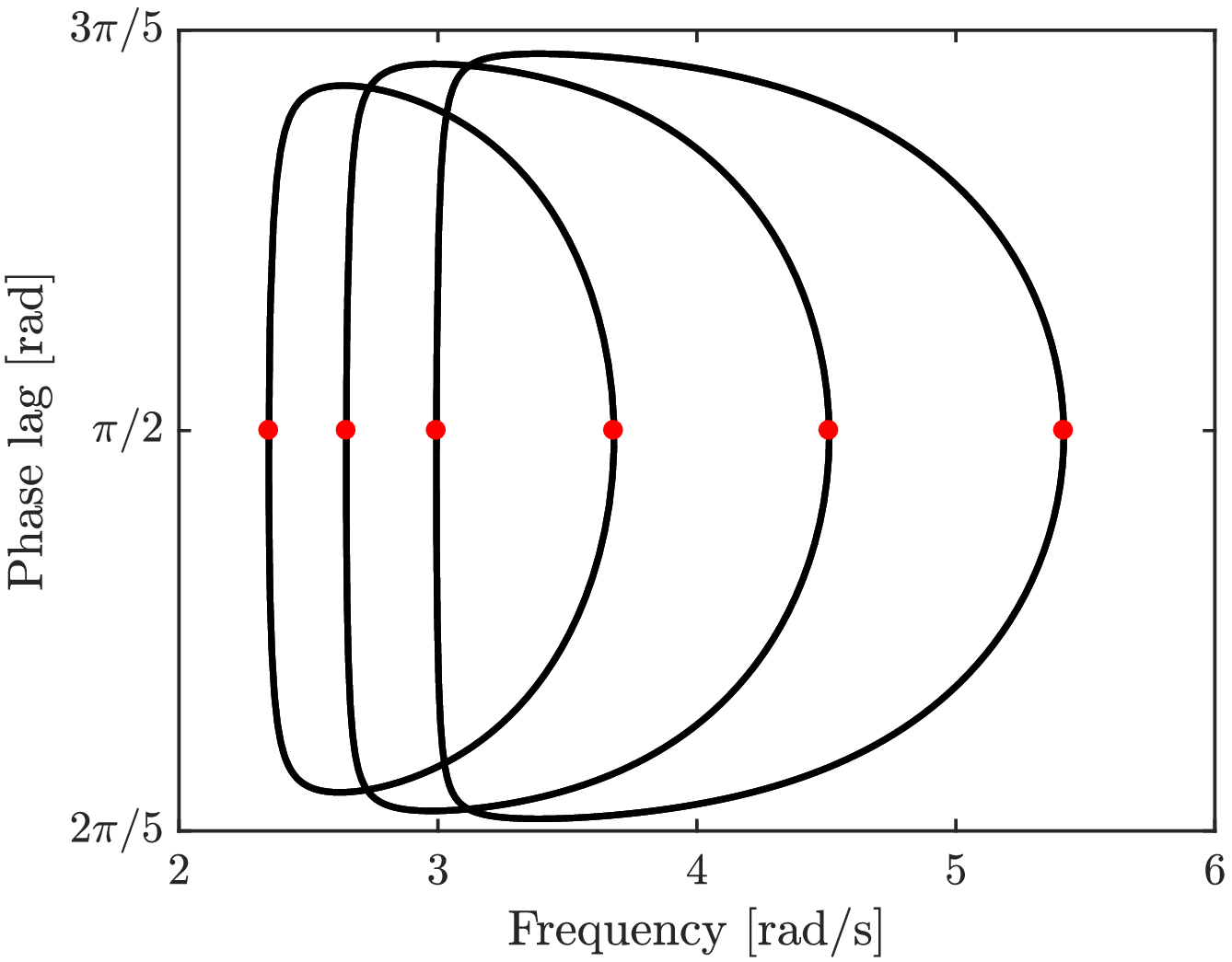} 
    \caption{\label{fig:RATIONAL_ODD_3_5_PHASE_NO_PRNM}} 
  \end{subfigure} 
  \caption{NFRCs of the 7:3 and 3:5 resonances for 3 forcing amplitudes (5N, 8N and 11N \& 3N, 5N and 8N, respectively): (\subref{fig:RATIONAL_ODD_7_3_NFRC_NO_PRNM}) amplitude - 7:3, (\subref{fig:RATIONAL_ODD_7_3_PHASE_NO_PRNM}) phase lag - 7:3, (\subref{fig:RATIONAL_ODD_SUBPLOT_3_5_NFRC_NO_PRNM}) amplitude - 3:5 and (\subref{fig:RATIONAL_ODD_3_5_PHASE_NO_PRNM}) phase lag - 3:5.}
  \label{fig:RATIONAL_KODD} 
\end{figure}

\begin{figure}[ht] 
  \begin{subfigure}[b]{0.5\linewidth}
    \centering
    \includegraphics[width=1\linewidth]{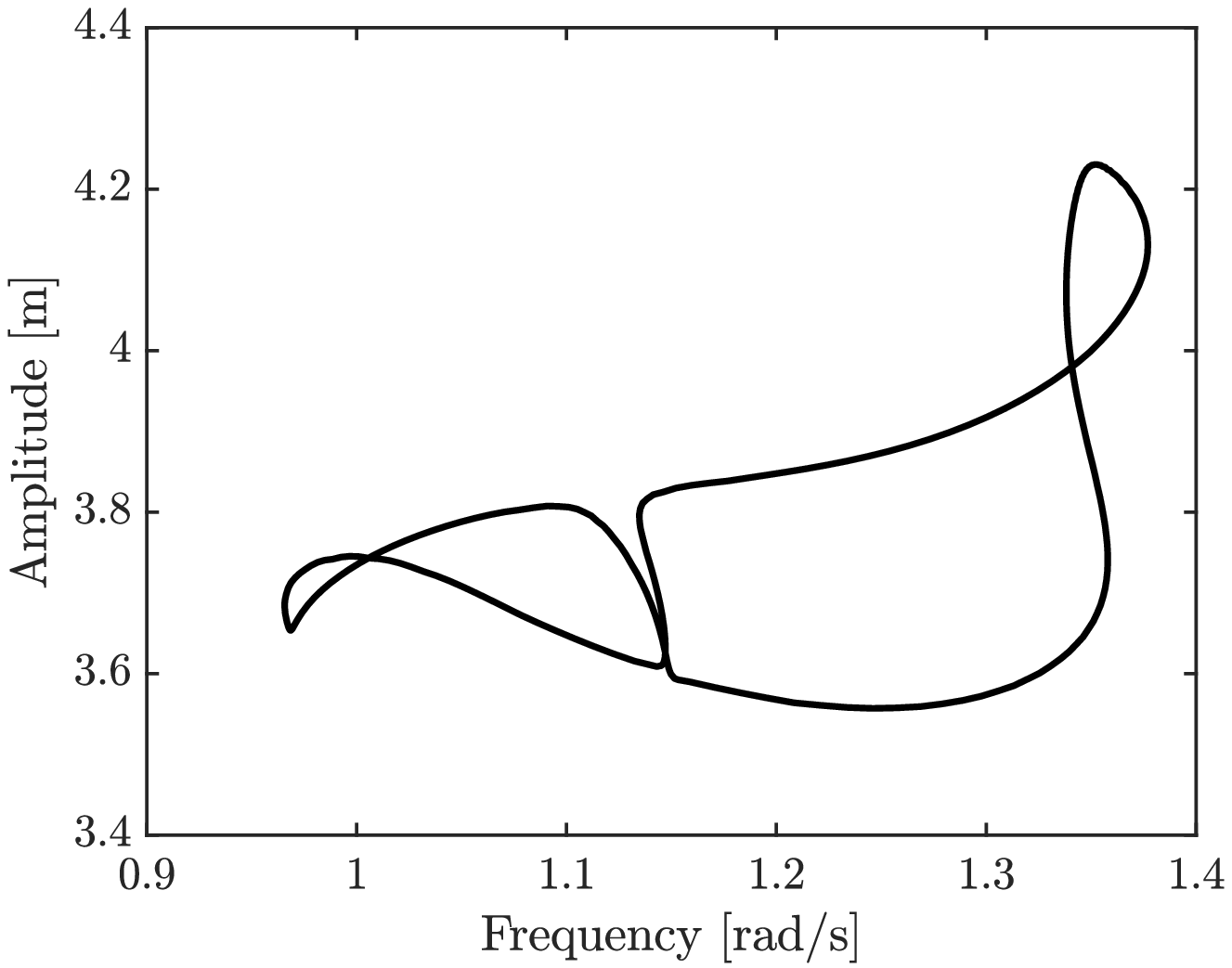}
    \caption{\label{fig:RATIONAL_ODD_7_3_NFRC_MERGING}}
  \end{subfigure}
  \begin{subfigure}[b]{0.5\linewidth}
    \centering
    \includegraphics[width=1\linewidth]{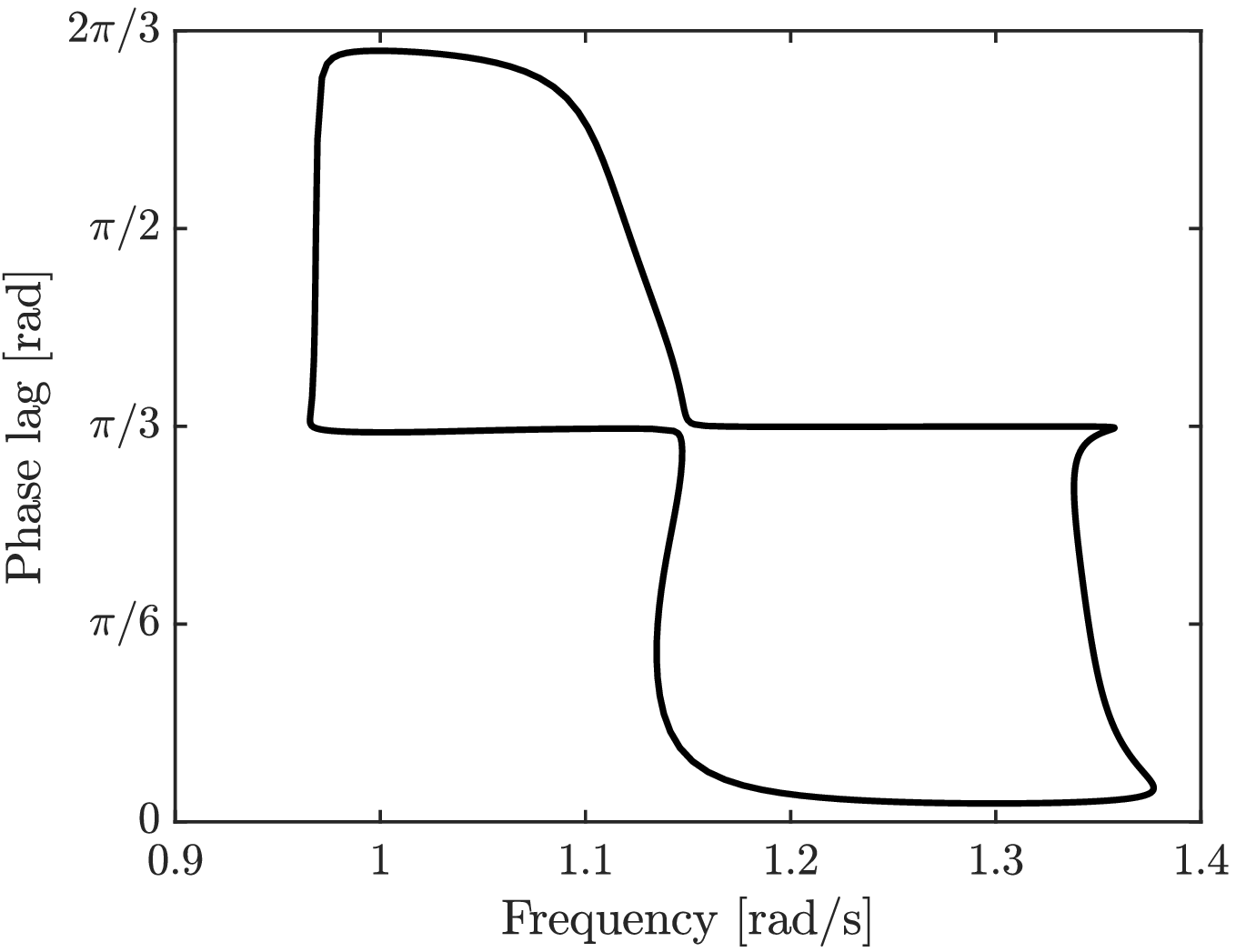}
    \caption{\label{fig:RATIONAL_ODD_7_3_PHASE_MERGING}}
  \end{subfigure} 
  \caption{Merging of the 7:3 resonances: (\subref{fig:RATIONAL_ODD_7_3_NFRC_MERGING}) amplitude and (\subref{fig:RATIONAL_ODD_7_3_PHASE_MERGING}) phase lag of the 7th harmonic component.}
  \label{fig:RATIONAL_ODD_MERGING} 
\end{figure}

Figure \ref{fig:RATIONAL_ODD_7_3_NFRC_NO_PRNM} and \ref{fig:RATIONAL_ODD_7_3_PHASE_NO_PRNM} present the amplitude and phase lag of the 7:3 resonance, respectively. Two different branches were found oscillating around phase lags of $\pi/2$ and $\pi/6$. Their geometry is more complex than the previously observed isolas with the consequence that the $\pi/2$ and $\pi/6$ points are not necessarily located at the extremities of the isolas. As the forcing amplitude increases, the two solutions get closer to each other; they eventually merge, as shown in Figure \ref{fig:RATIONAL_ODD_MERGING}. For the 3:5 resonance in Figures \ref{fig:RATIONAL_ODD_SUBPLOT_3_5_NFRC_NO_PRNM} and \ref{fig:RATIONAL_ODD_3_5_PHASE_NO_PRNM}, the phase lag centers around $\pi/2$.

\subsubsection{Either \texorpdfstring{$k$}{} or \texorpdfstring{$\nu$}{} is even}

\begin{figure}[ht] 
  \begin{subfigure}[b]{0.5\linewidth}
    \centering
    \includegraphics[width=1\linewidth]{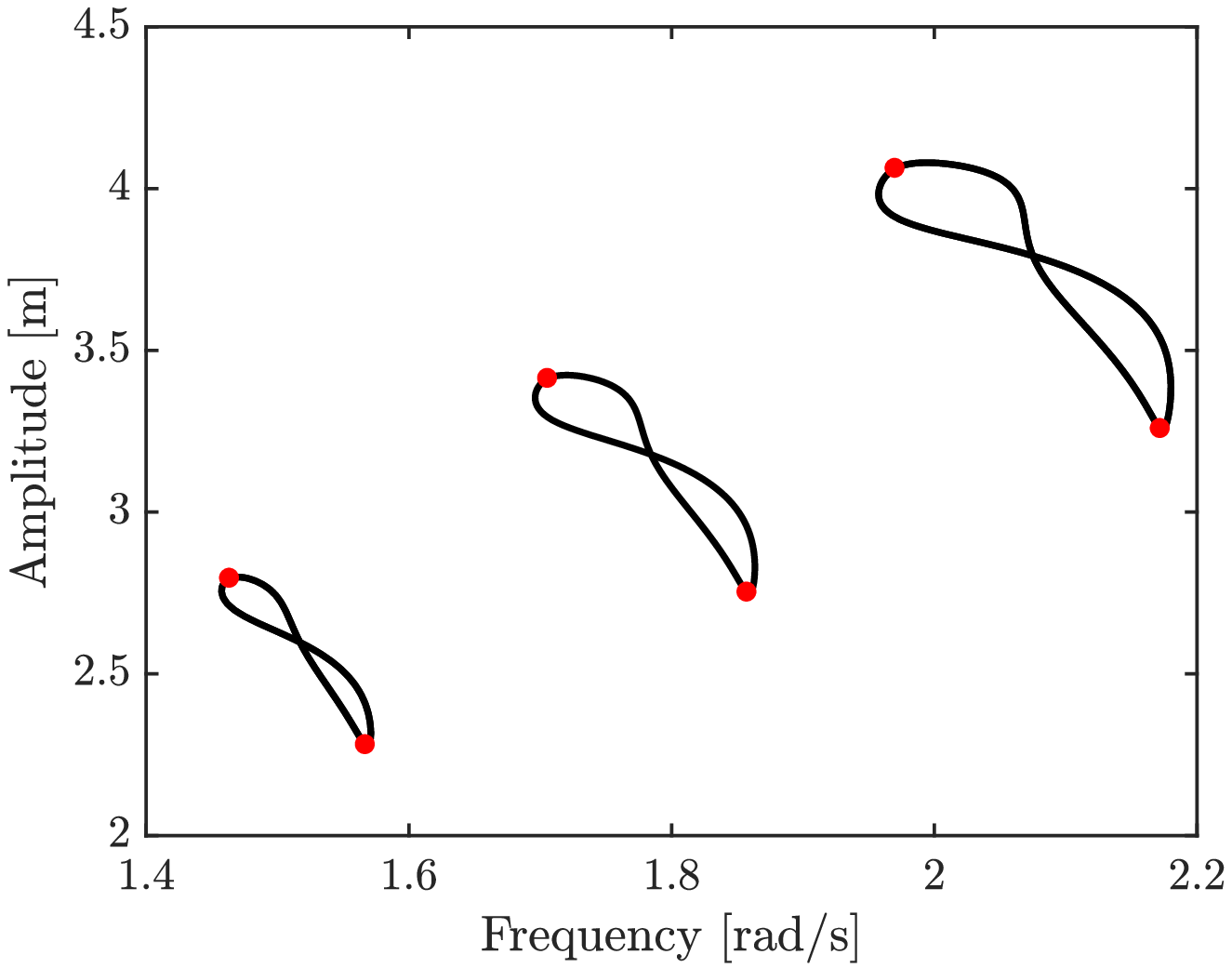} 
    \caption{\label{fig:RATIONAL_KODD_NUEVEN_3_2_NFRC_NO_PRNM}} 
  \end{subfigure}
  \begin{subfigure}[b]{0.5\linewidth}
    \centering
    \includegraphics[width=1\linewidth]{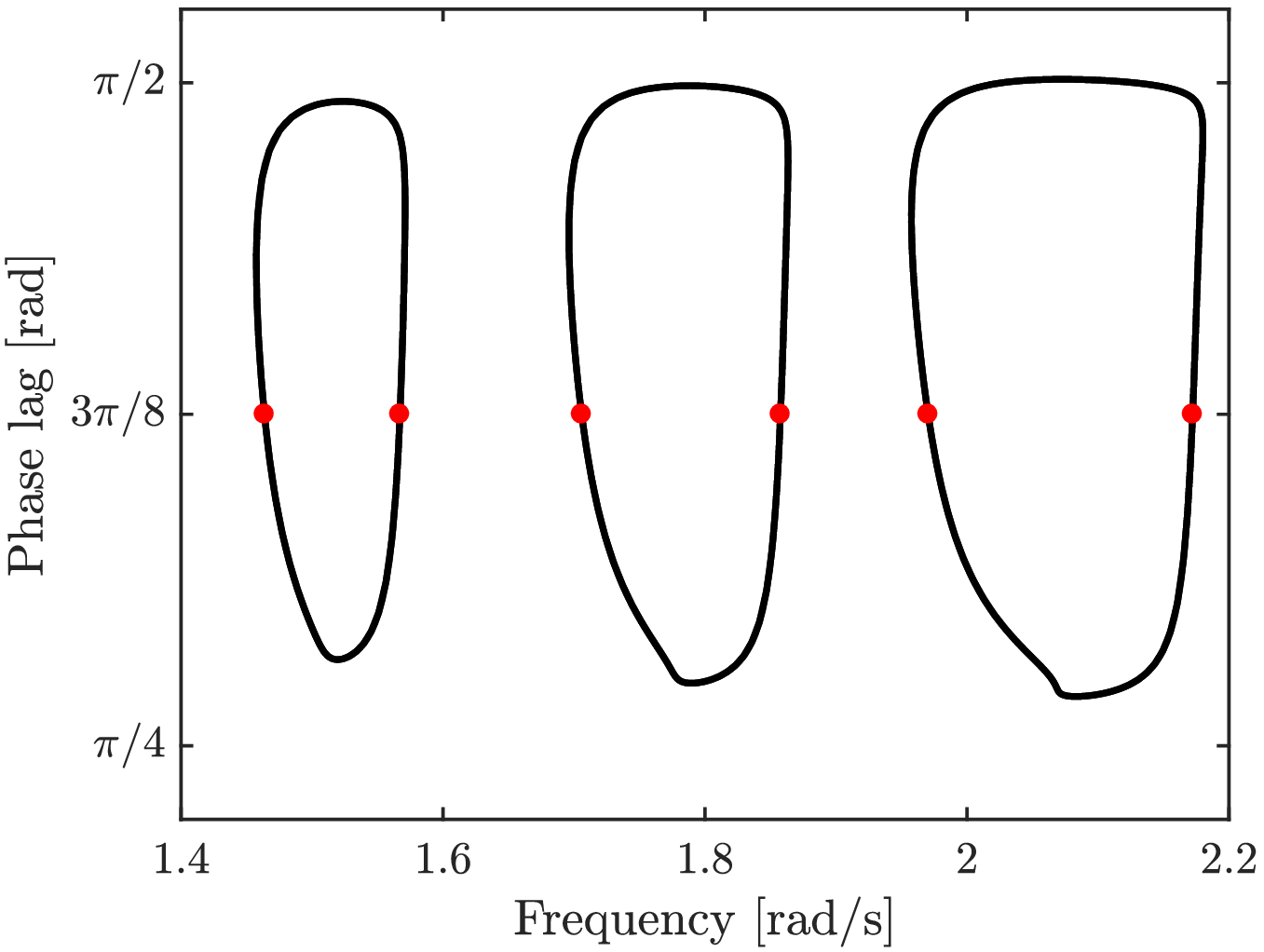} 
    \caption{\label{fig:RATIONAL_KODD_NUEVEN_3_2_PHASE_NO_PRNM}} 
  \end{subfigure} 
  
  \begin{subfigure}[b]{0.5\linewidth}
    \centering
    \includegraphics[width=1\linewidth]{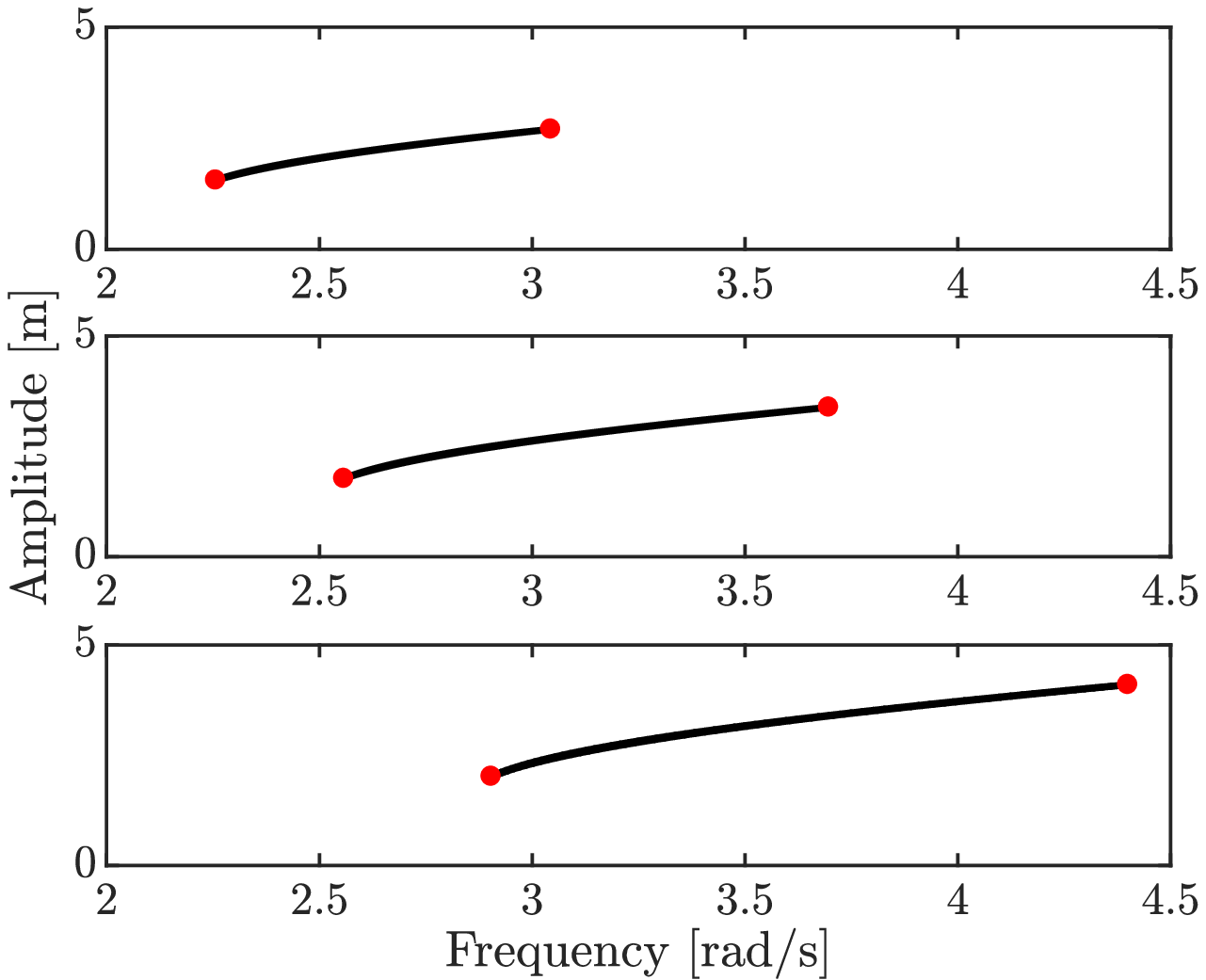}
    \caption{\label{fig:RATIONAL_KODD_NUEVEN_SUBPLOT_3_4_NFRC_NO_PRNM}} 
  \end{subfigure}
  \begin{subfigure}[b]{0.5\linewidth}
    \centering
    \includegraphics[width=1\linewidth]{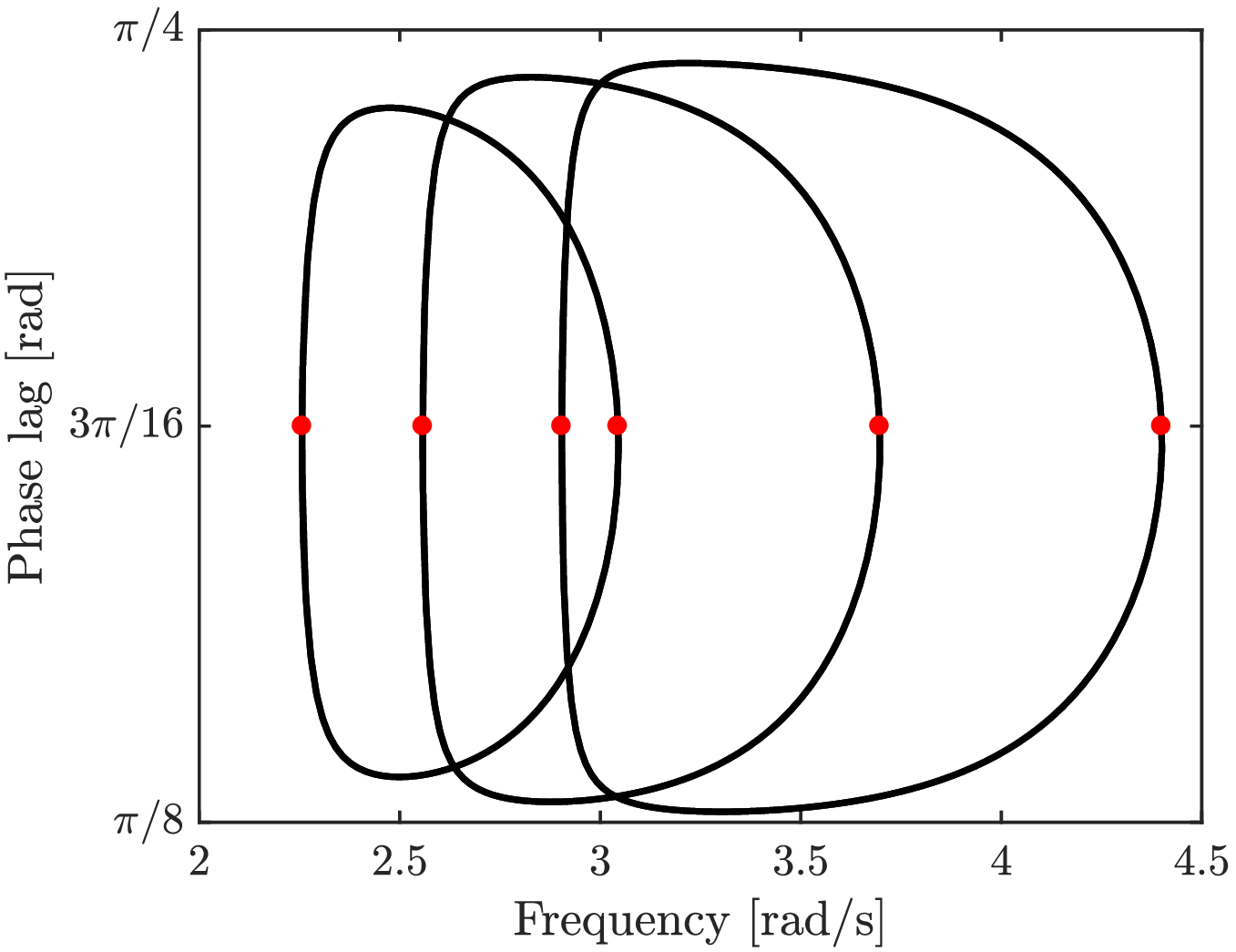} 
    \caption{\label{fig:RATIONAL_KODD_NUEVEN_3_4_PHASE_NO_PRNM}} 
  \end{subfigure} 
  \caption{NFRCs of the 3:2 and 3:4 resonances for 3 forcing amplitudes (3N, 5N and 8N): (\subref{fig:RATIONAL_KODD_NUEVEN_3_2_NFRC_NO_PRNM}) amplitude - 3:2, (\subref{fig:RATIONAL_KODD_NUEVEN_3_2_PHASE_NO_PRNM}) phase lag - 3:2, (\subref{fig:RATIONAL_KODD_NUEVEN_SUBPLOT_3_4_NFRC_NO_PRNM}) amplitude - 3:4 and (\subref{fig:RATIONAL_KODD_NUEVEN_3_4_PHASE_NO_PRNM}) phase lag - 3:4.}
  \label{fig:RATIONAL_KODD_NUEVEN} 
\end{figure}

\begin{figure}[ht] 
  \begin{subfigure}[b]{0.5\linewidth}
    \centering
    \includegraphics[width=1\linewidth]{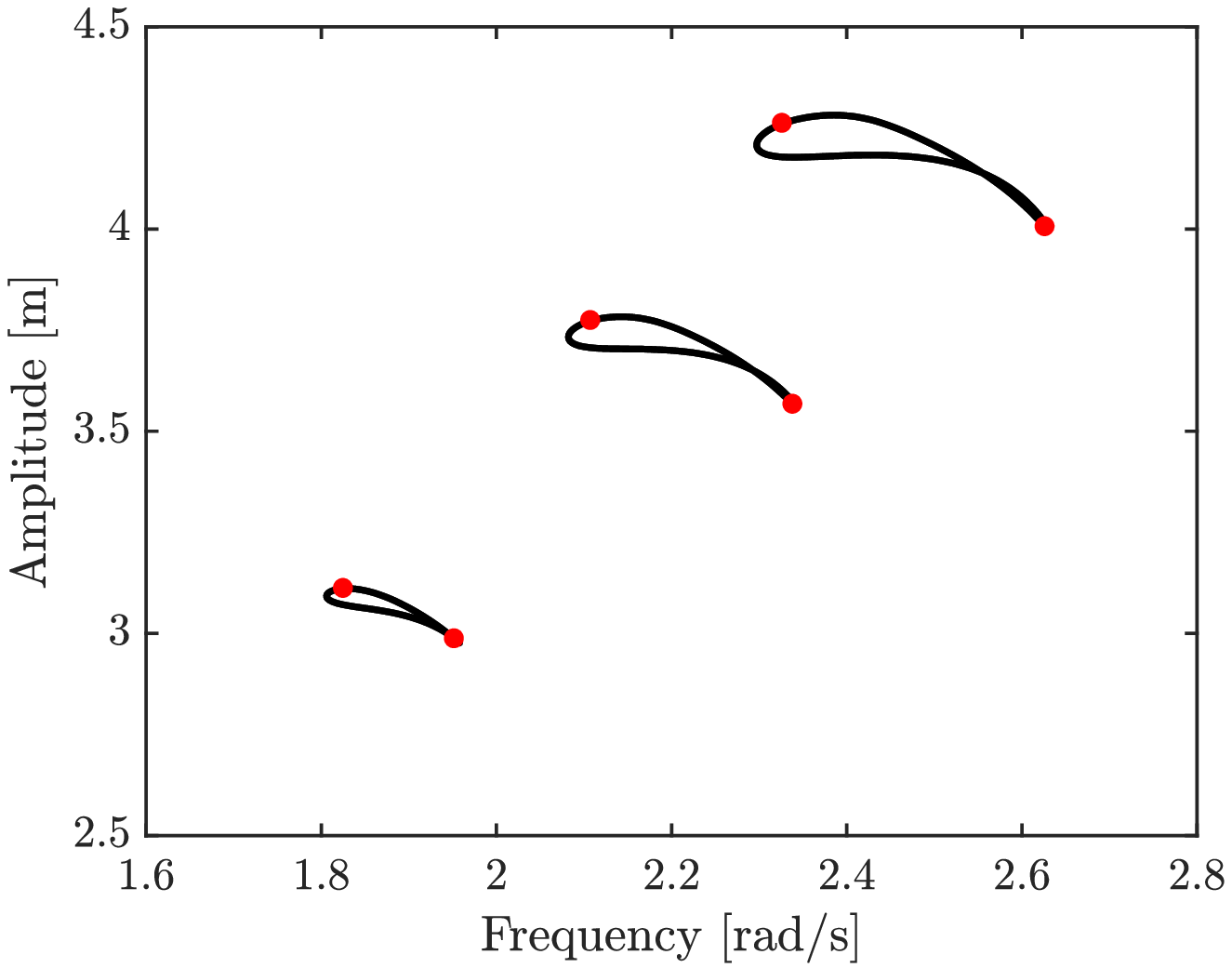} 
    \caption{\label{fig:RATIONAL_KEVEN_NUODD_4_3_NFRC_NO_PRNM}} 
  \end{subfigure}
  \begin{subfigure}[b]{0.5\linewidth}
    \centering
    \includegraphics[width=1\linewidth]{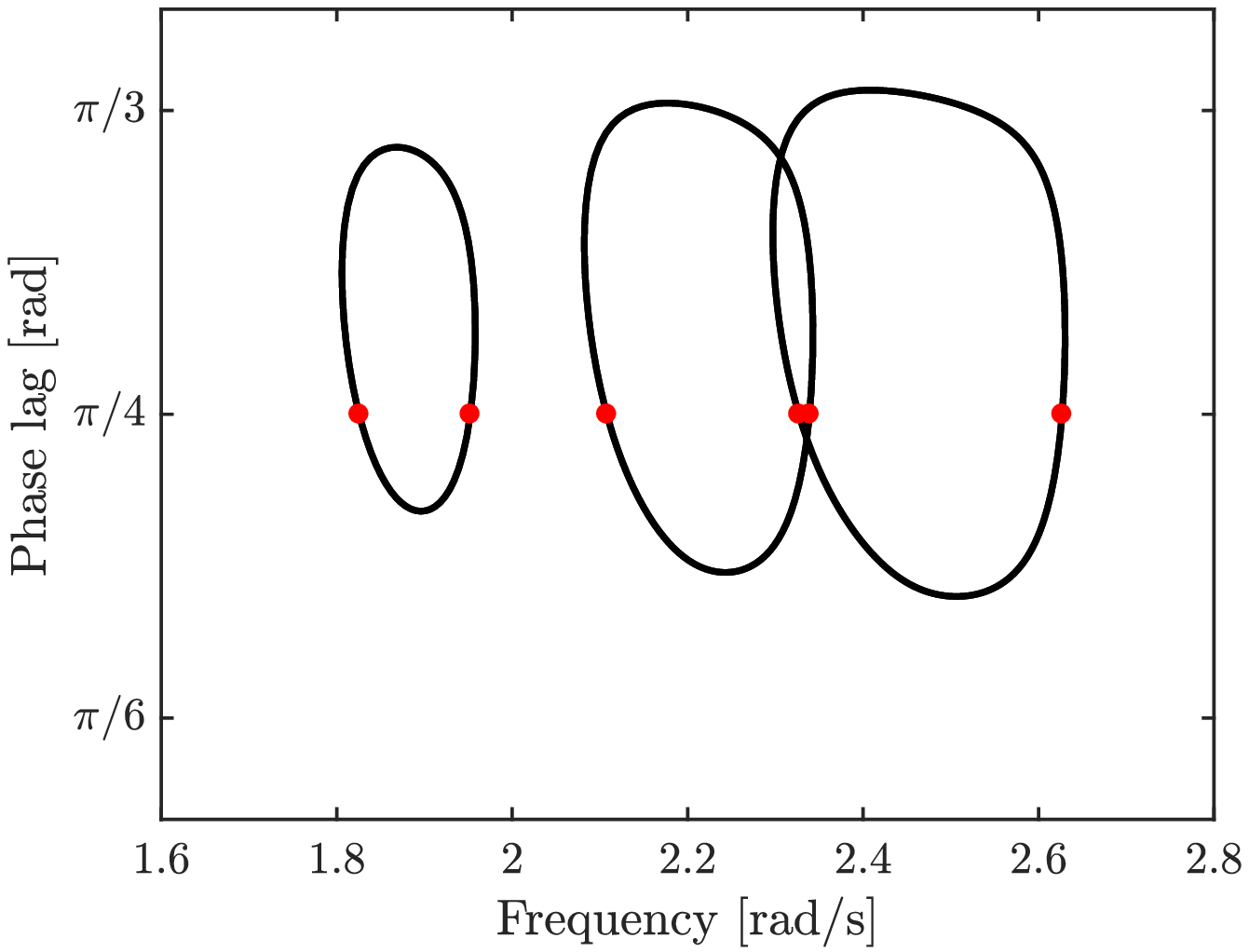} 
    \caption{\label{fig:RATIONAL_KEVEN_NUODD_4_3_PHASE_NO_PRNM}}
  \end{subfigure} 
  
  \begin{subfigure}[b]{0.5\linewidth}
    \centering
    \includegraphics[width=1\linewidth]{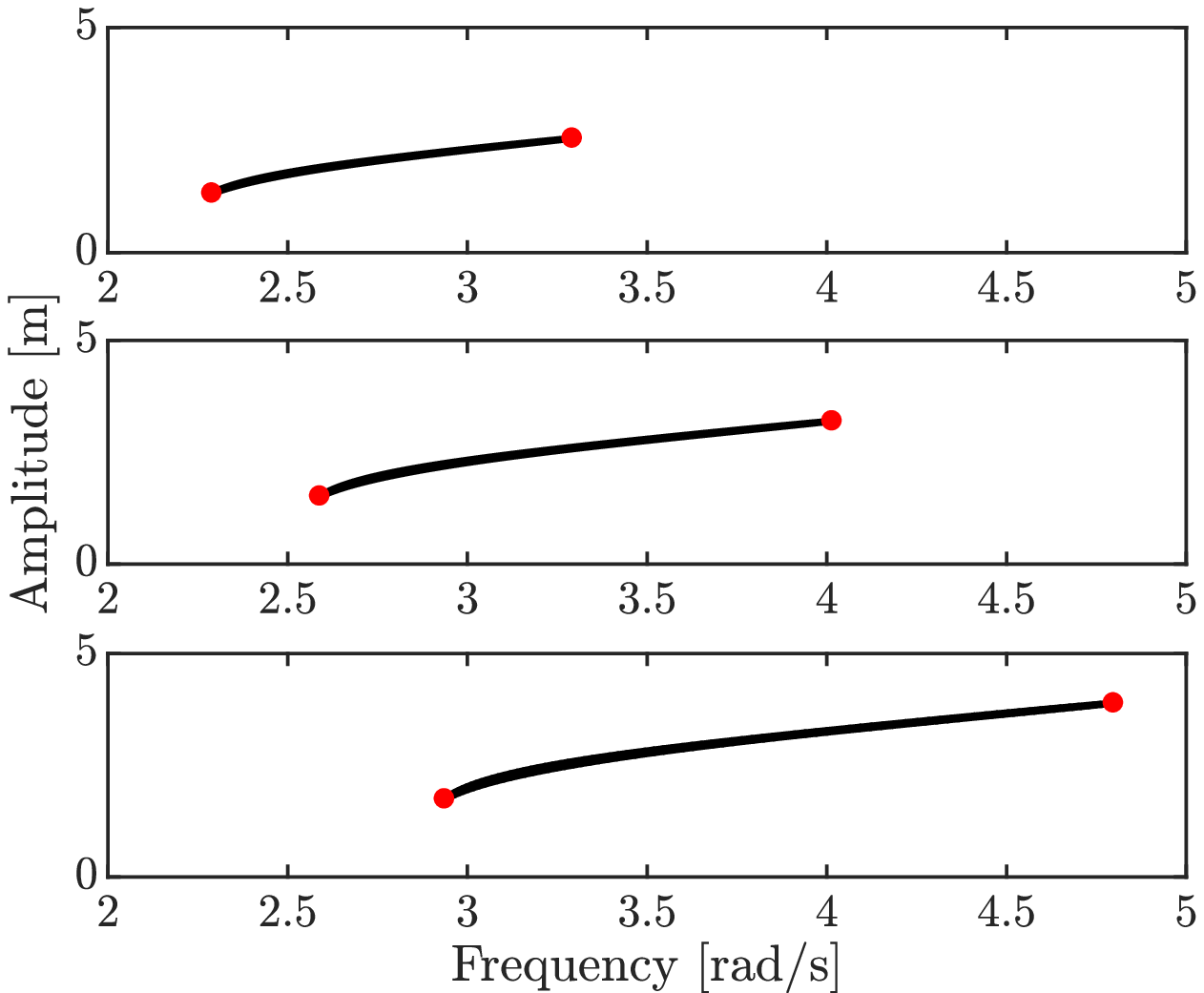}
    \caption{\label{fig:RATIONAL_KEVEN_NUODD_SUBPLOT_2_3_NFRC_NO_PRNM} } 
  \end{subfigure}
  \begin{subfigure}[b]{0.5\linewidth}
    \centering
    \includegraphics[width=1\linewidth]{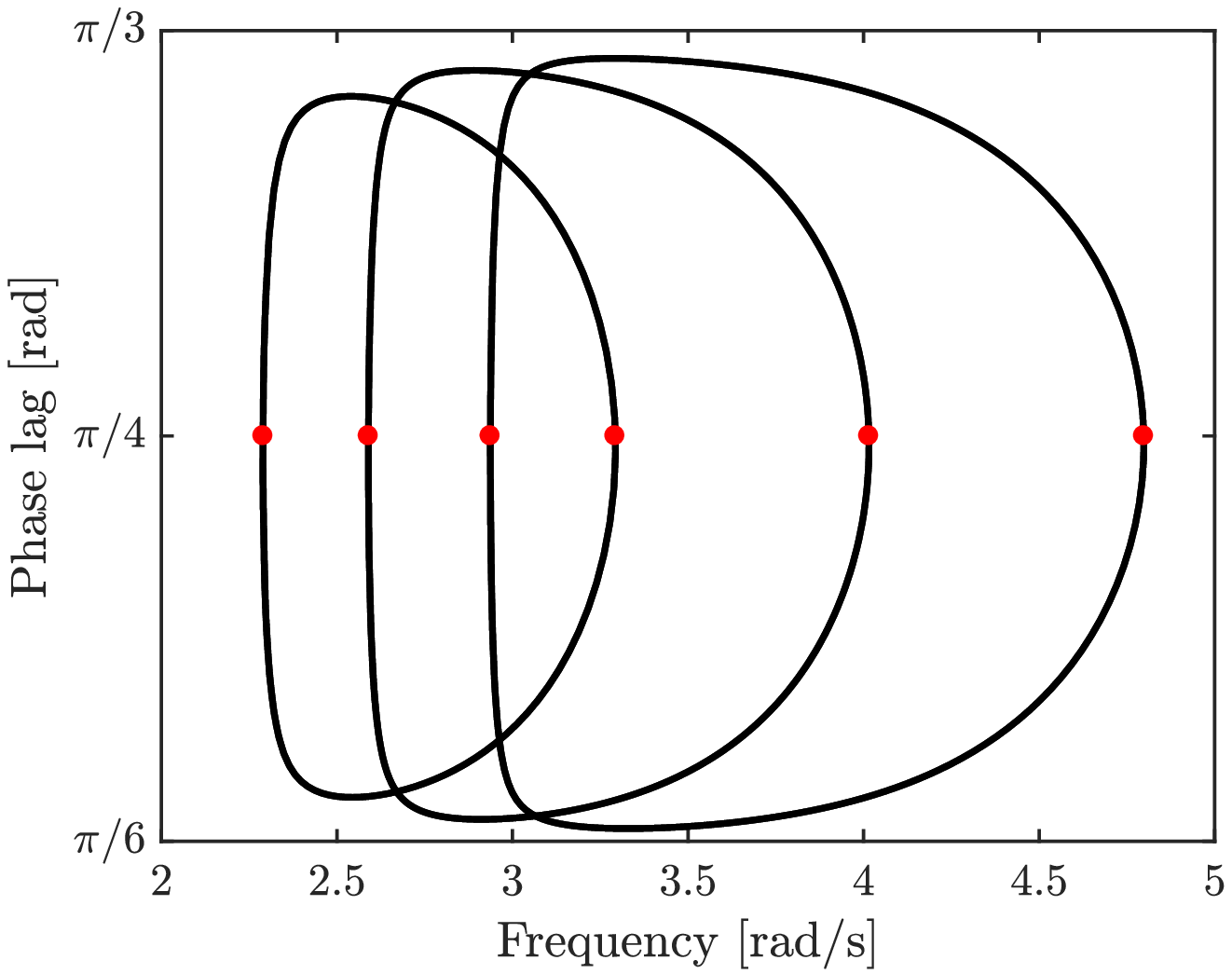} 
    \caption{\label{fig:RATIONAL_KEVEN_NUODD_2_3_PHASE_NO_PRNM}} 
  \end{subfigure} 
  \caption{NFRCs of the 4:3 and 2:3 resonances for 3 forcing amplitudes (3N, 5N and 7N \& 3N, 5N and 8N, respectively): (\subref{fig:RATIONAL_KEVEN_NUODD_4_3_NFRC_NO_PRNM}) amplitude - 4:3, (\subref{fig:RATIONAL_KEVEN_NUODD_4_3_PHASE_NO_PRNM}) phase lag - 4:3, (\subref{fig:RATIONAL_KEVEN_NUODD_4_3_PHASE_NO_PRNM}) amplitude - 2:3 and (\subref{fig:RATIONAL_KEVEN_NUODD_2_3_PHASE_NO_PRNM}) phase lag - 2:3}
  \label{fig:RATIONAL_KEVEN_NUODD} 
\end{figure}

Figures \ref{fig:RATIONAL_KODD_NUEVEN} and \ref{fig:RATIONAL_KEVEN_NUODD} represent the 3:2/3:4 and 4:3/2:3 resonances, respectively. For all of them, the phase lag is comprised in the interval [$3\pi/4\nu\pm\pi/4\nu$]; a second solution shifted by $\pi/\nu$ also exists and is the opposite of the first one.

\subsection{Summary}

The goal of this section was to carry out a detailed study of the phase lag between the forcing and the $k$-th harmonic of the displacement for the $k:\nu$ resonance of the Duffing oscillator.

Two important cases must be distinguished. First, when $k$ and $\nu$ are odd, the phase lag takes values in the interval $[\pi/2 \pm \pi/2\nu]$. There is thus always at least one point on the resonance branch for which phase quadrature is verified. Second, when either $k$ or $\nu$ is even, the phase lag oscillates around $[3\pi/4\nu \pm \pi/4\nu]$. Phase quadrature is no longer achieved for such branches. We note that several branches may exist for the $k:\nu$ resonance, but at least one of them verifies these findings.

In the different scenarios investigated, the resonance branch takes the form of either a resonance peak ($\nu=1$) or an isola ($\nu>1$). In the case of a resonance peak, there is a single point where the phase is $\pi/2$ or $3\pi/4$; it is located in the immediate vicinity of the point featuring the maximum amplitude on the branch, at least for the amount of damping considered herein. In the case of an isola, there are two points where the phase is $\pi/2$ or $3\pi/4\nu$. For an isola with simple geometry, e.g., for $1:\nu$ resonances, the two points are located at the extremities of the isola. For more complex geometries, i.e., for ultra-subharmonic resonances with $k>\nu$, the two points are not necessarily located at the extremities of the isola, but they remain relatively close from it. 

In summary, these results allow us to generalize the concept of a nonlinear phase resonance, defined originally for primary resonances under multi-point, multi-harmonic forcing \cite{PEETERSJSV}, to superharmonic, subharmonic and ultra-subharmonic resonances under single-point, single-harmonic forcing. To do so, a generalized phase resonance condition must be considered, i.e., the points where the phase lag is $\pi/2$ ($k$ and $\nu$ odd) or $3\pi/4\nu$ ($k$ or $\nu$ even) define the locus of nonlinear phase resonance. The numerical computation of these loci is the objective of the next section.

\section{Phase resonance nonlinear modes}\label{SECTION:PRNM}

\subsection{Velocity feedback}

We start from the harmonically-forced linear oscillator
\begin{equation}
    m\ddot{x}(t)+c\dot{x}(t)+kx(t)=f\sin{\omega t}.
    \label{eq:UnforcedDuffing}
\end{equation}
To drive this system into resonance requires careful tuning of the excitation frequency $\omega$. A more efficient strategy to operate the system into one of its normal modes is to apply direct velocity feedback \cite{BABI}
\begin{equation}
    m\ddot{x}(t)+c\dot{x}(t)+kx(t)-\mu\dot{x}(t)=0
    \label{eq:UnforcedDuffingFund}
\end{equation}
where the feedback term $\mu\dot{x}(t)$ plays the role of {\it virtual harmonic forcing}. Because this virtual forcing and the velocity are collinear, phase quadrature with the displacement $x(t)$, and, hence, phase resonance, is naturally enforced when $\mu=c$. This strategy was also used for driving nonlinear systems into primary resonances \cite{KRACK,BABI2}.

We leverage the results in Section 3 to generalize the concept of a velocity feedback for the different families of resonance of a nonlinear system. Taking the Duffing oscillator as an illustrative example, we obtain
\begin{itemize}
    \item For the primary resonance,
\begin{equation}
    m\ddot{x}(t)+c\dot{x}(t)+kx(t)+k_{nl}x^3(t)-\mu\dot{x}_{1,T}(t)=0
    \label{eq:PRNM_1}
\end{equation}
where, unlike \cite{KRACK,BABI2}, the velocity feedback $\dot{x}_{1,T}$ contains only the first harmonic of the velocity. The subscript $T$ indicates that the feedback is a $T$-periodic signal.
\item For $k:\nu$ secondary resonances with phase quadrature points,
\begin{equation}
    m\ddot{x}(t)+c\dot{x}(t)+kx(t)+k_{nl}x^3(t)-\mu\dot{x}_{k,T}\left(t\right)=0
    \label{eq:PRNM_k}
\end{equation}
where the velocity feedback is the $k$-th harmonic of the velocity transformed into a $T$-periodic signal.


\item For $k:\nu$ secondary resonances without quadrature points,
\begin{equation}
    m\ddot{x}(t)+c\dot{x}(t)+kx(t)+k_{nl}x^3(t)-\mu\dot{x}_{k,T}\left(t-\alpha\right)=0
    \label{eq:PRNM_k_nu}
\end{equation}
In this case, the feedback is delayed by $\alpha=\frac{1}{\omega_k}\left(\frac{\pi}{2}-\frac{3\pi}{4\nu}\right)$.
\end{itemize}

Considering the $k$-th harmonic of the displacement in Equation \eqref{eq:FourierDispPhase}, 
\begin{equation}
    {x}_k(t)= A_k\sin{(\omega_k t - \phi_k)}
\end{equation}
the velocity feedback $\mu\dot{x}_{k,T}(t)$ can be obtained after transforming $\dot{x}_k(t)$ into a $T$-periodic signal   
\begin{equation}
    \mu\dot{x}_{k,T}(t)  = \mu \omega_k A_k\cos{(\omega t -\phi_k)}
    \label{eq:velocity}
\end{equation}
Equation \eqref{eq:velocity} proves that the velocity feedback and the $k$-th harmonic of the displacement are in quadrature, as sought.

Considering now the general case, i.e., when the velocity feedback is
\begin{equation}
    \mu\dot{x}_{k,T}(t-\alpha)  = \mu \omega_k A_k\cos{(\omega \left(t-\alpha\right) -\phi_k)}
    \label{eq:velocity_forcing}
\end{equation}
evidences that the periodic solutions of Equations \eqref{eq:EOM1D} and \eqref{eq:PRNM_k_nu} are identical since the velocity feedback is equivalent to classical harmonic forcing of frequency $\omega$. 

\begin{figure}[ht] 
  \begin{subfigure}[b]{0.5\linewidth}
    \centering
    \includegraphics[width=1\linewidth]{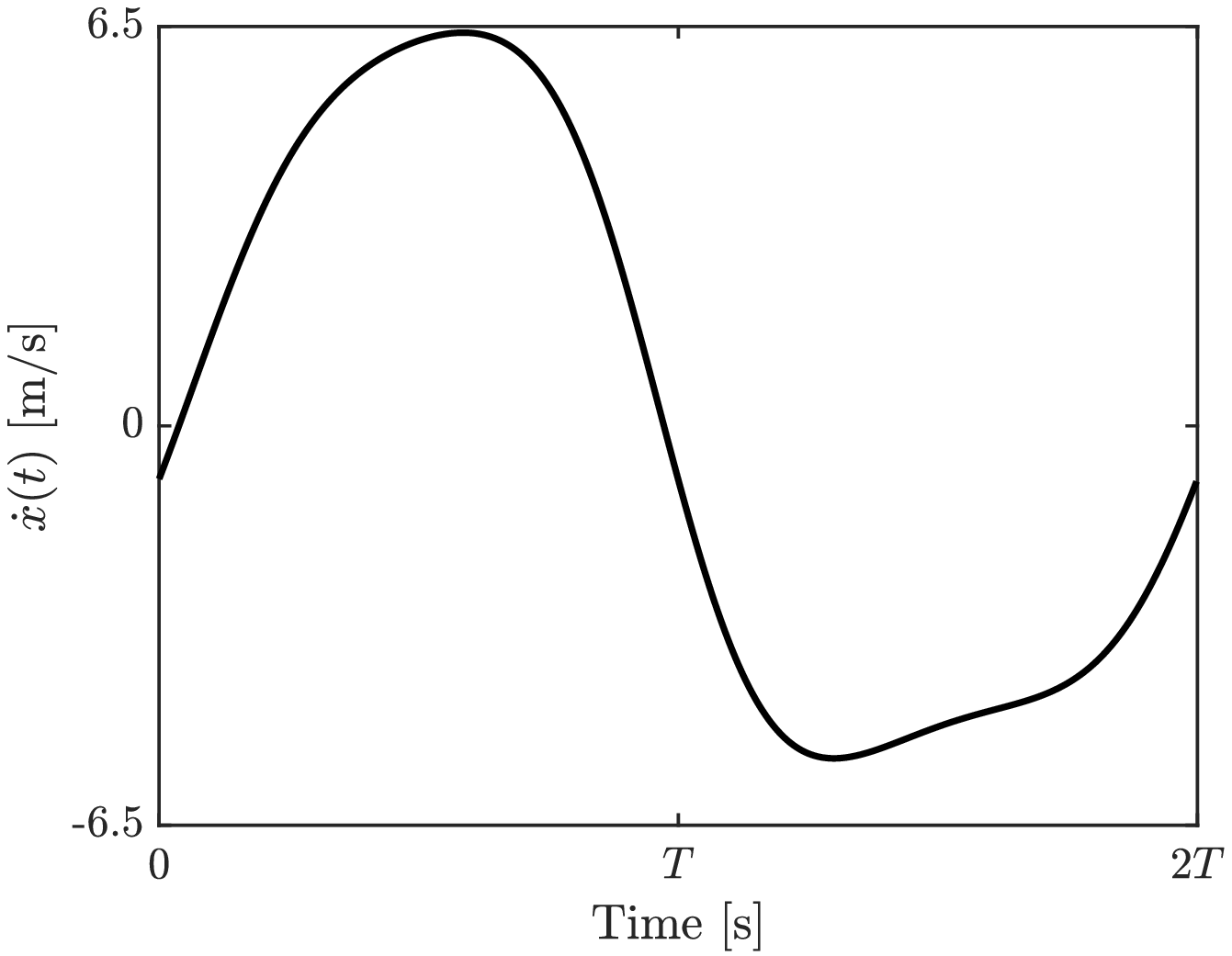} 
    \caption{\label{fig:FORCING_STEP1}} 
  \end{subfigure}
  \begin{subfigure}[b]{0.5\linewidth}
    \centering
    \includegraphics[width=1\linewidth]{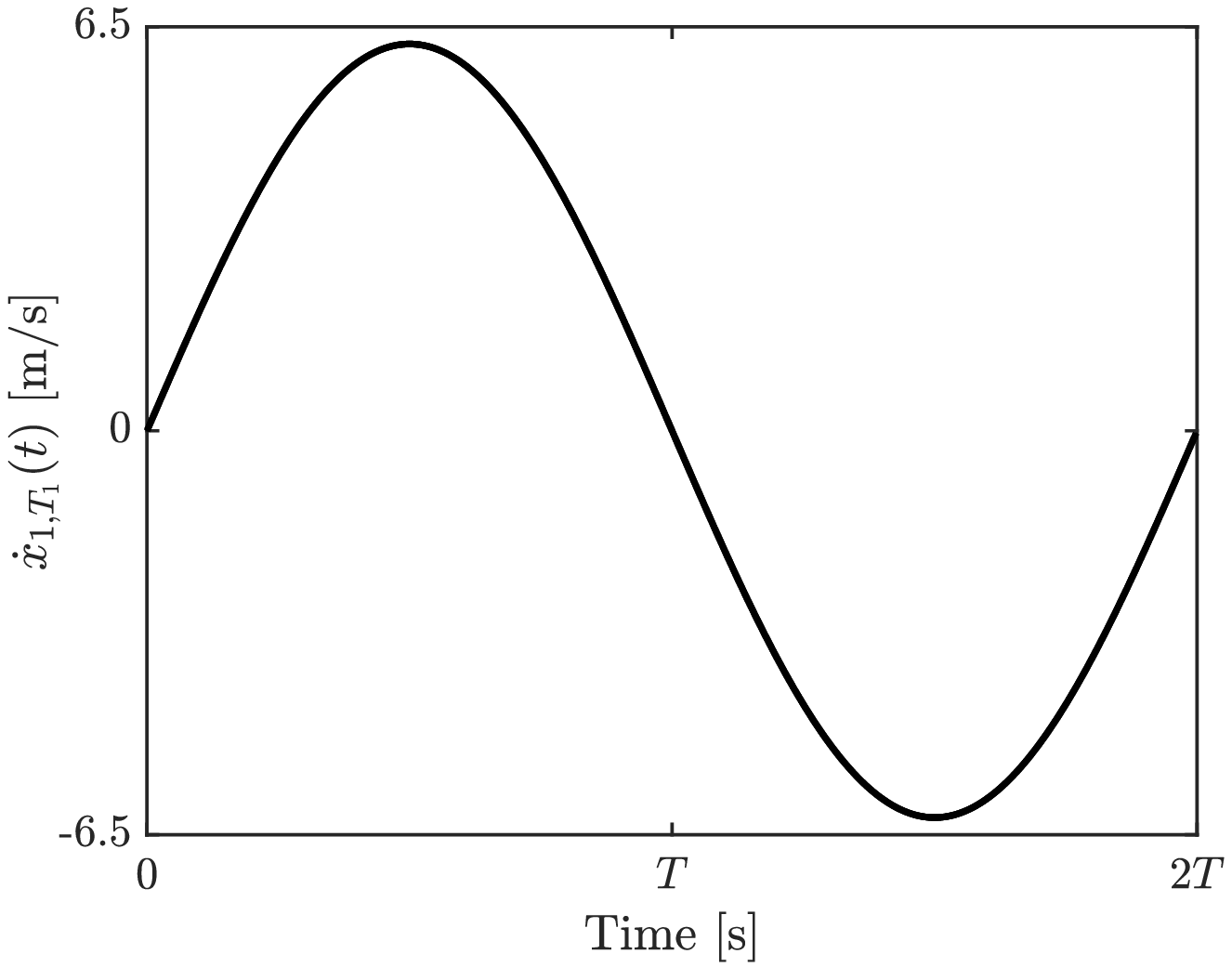} 
    \caption{\label{fig:FORCING_STEP2}} 
  \end{subfigure} 
  \begin{subfigure}[b]{0.5\linewidth}
    \centering
    \includegraphics[width=1\linewidth]{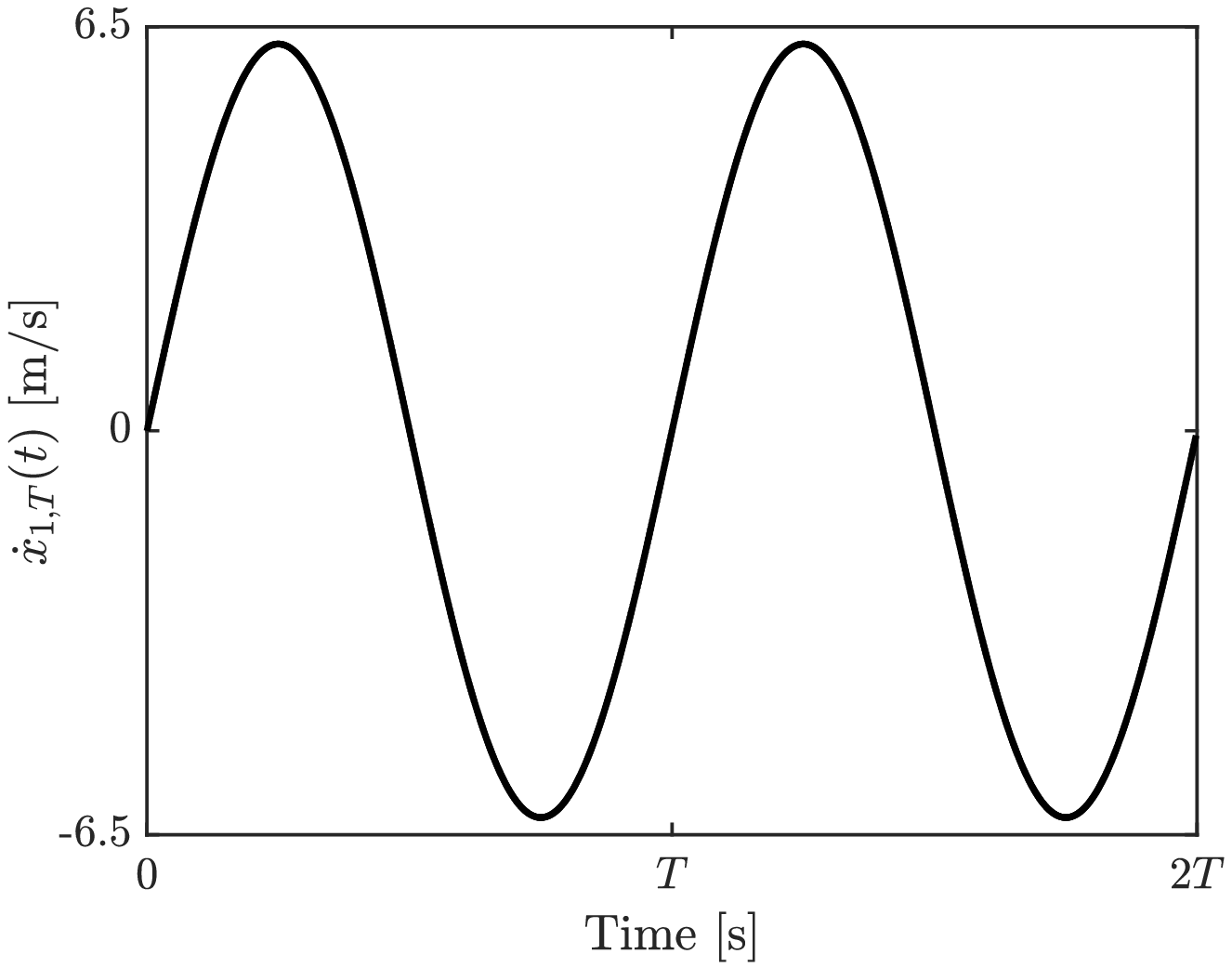} 
    \caption{\label{fig:FORCING_STEP3} } 
  \end{subfigure}
  \begin{subfigure}[b]{0.5\linewidth}
    \centering
    \includegraphics[width=1\linewidth]{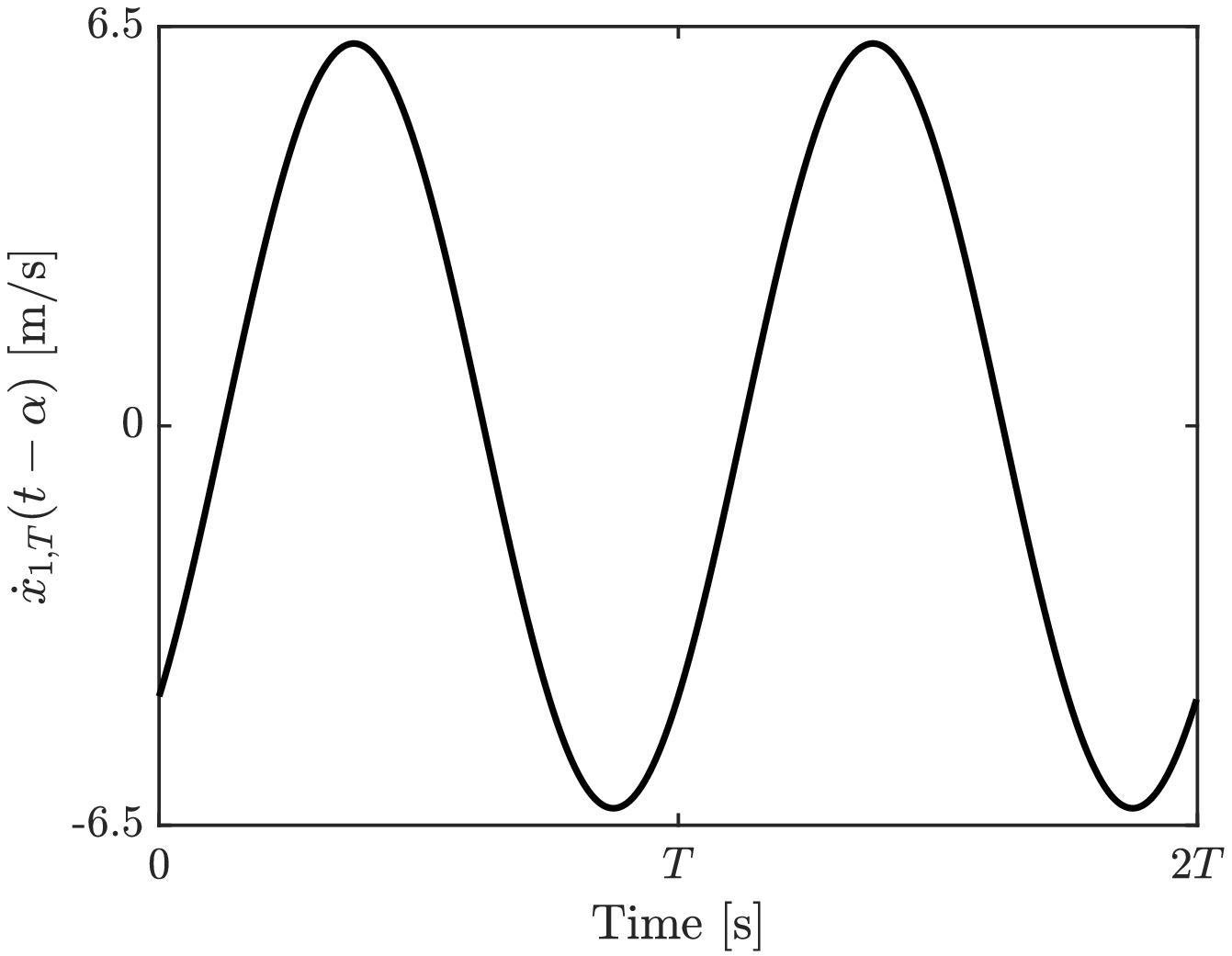} 
    \caption{\label{fig:FORCING_STEP4}} 
  \end{subfigure} 
  \caption{Calculation of the velocity feedback: (\subref{fig:F_0_01}) original velocity, (\subref{fig:F_0_25}) after step 1 (filtering), (\subref{fig:F_1}) after step 2 ($T$-periodic) and (\subref{fig:F_3}) after step 3 (delay).}
  \label{fig:FORCING} 
\end{figure}

For illustration, the 1:2 subharmonic resonance ($k=1$, $\nu=2$) is taken as an example. Figure \ref{fig:FORCING} shows the three steps to calculate the velocity feedback $\dot{x}_{1,T}(t-\alpha)$ in Figure \ref{fig:FORCING_STEP4} from the original velocity $\dot{x}(t)$ in Figure \ref{fig:FORCING_STEP1}. The first step filters out all but the first harmonic of the velocity to obtain $\dot{x}_{1,T_1}(t)$ in Figure \ref{fig:FORCING_STEP2}. This signal has a period $T_1=2T$ and is thus transformed during the second step into the $T$-periodic signal $\dot{x}_{1,T}(t)$ shown in Figure \ref{fig:FORCING_STEP3}. The third step shifts the resulting signal by the delay $\alpha= \pi/8\omega_k=\pi/4\omega$ to obtain the final feedback in Figure \ref{fig:FORCING_STEP4}.

\subsection{Computational framework}

Generalizing the results in the previous section to systems with $n$ degrees of freedom (DOFs) subjected to harmonic forcing at the $l$-th DOF
\begin{equation}
    \mathbf{M} \ddot{\mathbf{x}}(t) + \mathbf{C}\dot{\mathbf{x}}(t) + \mathbf{K}\mathbf{x}(t)+ \mathbf{f}_{nl}(\mathbf{x}(t) ,\mathbf{\dot{x}}(t)) = \mathbf{f}\sin{\omega t},
    \label{eq:EOMt}
\end{equation}
we define the PRNMs of the $k:\nu$ resonance as 
\begin{center}
    \textit{the periodic responses obtained by replacing the harmonic forcing by a $T$-periodic velocity feedback comprising the $k$-th harmonic. The feedback is to be delayed by $\alpha=\frac{1}{\omega_k}\left(\frac{\pi}{2}-\frac{3\pi}{4\nu}\right)$ when $k$ or $\nu$ is even.}
\end{center}
Mathematically, we solve:
\begin{equation}
    \mathbf{M} \ddot{\mathbf{x}}(t) + \mathbf{C}\dot{\mathbf{x}}(t) + \mathbf{K}\mathbf{x}(t)+ \mathbf{f}_{nl}(\mathbf{x}(t) ,\mathbf{\dot{x}}(t)) - \mu\dot{\mathbf{x}}_{k,T}(t-\alpha) = \mathbf{0},
    \label{eq:PRNMt}
\end{equation}
where the feedback vector $\dot{\mathbf{x}}_{k,T}(t-\alpha)$ has only one non-zero entry at the DOF $l$, and, without loss of generality, we assume that the feedback gain $\mu$ is strictly positive. Because the velocity feedback is a single-harmonic component of period $T=2\pi/\omega$, the periodic orbits of \eqref{eq:PRNMt} are actual periodic orbits of \eqref{eq:EOMt}, i.e., those with the assigned phase lag.  Equation \eqref{eq:PRNMt} should thus be interpreted as an effective reformulation of \eqref{eq:EOMt}  which targets the calculation of the locus of nonlinear phase resonance of the NFRCs.

We now aim to develop an efficient computational framework to obtain accurate numerical approximations of the PRNMs and their oscillations frequencies. In this context, the HBM is a particularly appropriate method because it naturally separates the responses into different harmonics and requires no interpolation to render the signal $T$-periodic. The displacement and nonlinear force vectors are thus approximated by truncated Fourier series:
\begin{equation}
    \mathbf{x}(t) = \frac{\mathbf{c}_0^x}{\sqrt{2}}+\sum_{k=1}^{N_h}\left( \mathbf{s}_{k}^x \sin \omega_k t + \mathbf{c}_{k}^x \cos \omega_k t\right)
    \label{eq:FourierDispNh}
\end{equation}
\begin{equation}
    \mathbf{f}_{nl}(t) = \frac{\mathbf{c}_0^f}{\sqrt{2}}+\sum_{k=1}^{N_h}\left( \mathbf{s}_{k}^f \sin \omega_k t + \mathbf{c}_{k}^f \cos \omega_k t\right)
    \label{eq:FourierFnlNh}
\end{equation}
or, in a more compact form,
  \begin{equation}
     \mathbf{x}(t) = (\mathbf{Q}(t)\otimes \mathbf{I}_n)\mathbf{X}
     \label{xt}
 \end{equation}
  \begin{equation}
     \mathbf{f}_{nl}(t) = (\mathbf{Q}(t)\otimes \mathbf{I}_n)\mathbf{F}_{nl}
     \label{fnlt}
 \end{equation}
 where $\mathbf{Q}(t)=\left[1/\sqrt(2) \,\,\sin \omega_1 t \,\,\cos \omega_1 t \,\,\dots\,\, \sin \omega_{N_h} t \,\,\cos \omega_{N_h} t\right]$, $\otimes$ stands for the Kronecker tensor product, $\mathbf{I}_n$ is the identity matrix of size $n$, and $\mathbf{X}$ and $\mathbf{F}_{nl}$ are the vectors containing the Fourier coefficients of the displacement and nonlinear forces, respectively. Similarly, the velocity feedback can be written using the Fourier coefficients of the displacement as:
   \begin{equation}
     \mu\dot{\mathbf{x}}_{k,T}(t-\alpha) =\mu (\mathbf{Q}(t)\otimes \mathbf{I}_n)\left(\mathbf{R_\alpha}\mathbf{T}_T\mathbf{T}_f\boldsymbol{\nabla}(\omega)\otimes \mathbb{I}_{l}\right)\mathbf{X}
     \label{xdot_shift_t}
 \end{equation}
 where
\begin{itemize}
\item $\left(\boldsymbol{\nabla}(\omega)\otimes \mathbb{I}_{l}\right)\mathbf{X}$ contains only the Fourier coefficients of the velocity measured at the $l$-th DOF. $\mathbb{I}_{l}$ is a $(n\times n)$ a null matrix except for the $l$-th diagonal term which is equal to 1. The operator $\boldsymbol{\nabla}(\omega)=$ diag$(0,\boldsymbol{\nabla}_1,...,\boldsymbol{\nabla}_k,...,\boldsymbol{\nabla}_{N_h})$ is a differential operator coming from the time derivative of $\mathbf{Q}(t)$with
\begin{equation}
	\boldsymbol{\nabla}_k=
	\begin{bmatrix}
	0 & -\omega_k \\
	\omega_k  & 0
    \end{bmatrix}
\end{equation}
\item $\mathbf{T}_f$ is a $(2N_h+1)×(2N_h+1)$ null matrix except for 
the two diagonal elements corresponding to the $k$-th harmonic which are equal to 1. $\mathbf{T}_f$ thus filters out all harmonic components different from $k$, as schematized in the time domain in Figure \ref{fig:FORCING_STEP2}.
\item $\mathbf{T}_T$ is a $(2N_h+1)\times(2N_h+1)$ null matrix except for 
the elements whose rows and columns correspond to the Fourier coefficients of $\omega$ and $\omega_k$, respectively. Those elements are equal to 1. $\mathbf{T}_T$ transforms the velocity into a $T$-periodic signal, see Figure \ref{fig:FORCING_STEP3}.
\item $\mathbf{R}_\alpha$ = diag$(0,\mathbf{0},...,\mathbf{R}_\nu,...,\mathbf{0})$ is a rotation matrix which shifts the harmonic component by an angle $\frac{\nu}{k}\alpha\omega_k=\alpha\omega$, as in Figure \ref{fig:FORCING_STEP4}, where 
\begin{equation}
    \mathbf{R}_\nu=
    \begin{pmatrix}
        \cos\alpha\omega & \sin\alpha\omega\\
        -\sin\alpha\omega & \cos\alpha\omega
    \end{pmatrix}.
    \label{eq:R}
\end{equation}
For primary and odd secondary resonances, $\alpha=0$, which means that the rotation matrix $\mathbf{R}_\nu$ is the identity matrix. For the other resonances, $\alpha=\frac{1}{\omega_k}\left(\frac{\pi}{2}-\frac{3\pi}{4\nu}\right)$.
\end{itemize}
 
 By inserting Equations (\ref{xt}), (\ref{fnlt}) and  \eqref{xdot_shift_t} into Equation (\ref{eq:PRNMt}) and removing the time dependency with a Galerkin procedure, the PRNMs can be obtained by solving the system:
\begin{equation}
    \mathbf{A}(\omega)\mathbf{X}+\mathbf{F}_{nl}(\mathbf{X})-\mu\left(\mathbf{R_\alpha}\mathbf{T}_T\mathbf{T}_f\boldsymbol{\nabla}(\omega)\otimes \mathbb{I}_{l}\right)\mathbf{X}=0
    \label{eq:PRNM_FREQ}
\end{equation}
where $\mathbf{A}(\omega)=\boldsymbol{\nabla}^2(\omega)\otimes \mathbf{M}+\boldsymbol{\nabla}(\omega)\otimes \mathbf{C} + \mathbf{I}_{2N_h+1} \otimes \mathbf{K}$ is the dynamic stiffness matrix.

Equation (\ref{eq:PRNM_FREQ}) provides $2N_h+1$ equations for $2N_h+2$ unknowns, namely the vector $\mathbf{X}$ and the gain $\mu$. An additional equation is therefore required to close the system. This equation, termed phase condition, sets the sine coefficient of the $k$-th harmonic component to 0. Eventually, a resonance can be characterized at different amplitudes by taking the frequency $\omega$ as a continuation parameter. The PRNMs and the corresponding resonance frequencies are obtained through vector $\mathbf{X}$ and $\omega$, respectively. To retrieve the phase lag $\phi_k$ defined in Equation \eqref{eq:FourierDispPhase}, the relation $\phi_k=\atantwo(-c_k^x,s_k^x)-\omega_k\alpha - \frac{\omega_k}{\omega}\left(\atantwo(-c_k^x,s_k^x)-\frac{\pi}{2}\right)$ must be considered.

To initiate the continuation process, we require an initial guess $(\mathbf{X}_{(0)},\mu_{(0)})$. This initial guess is taken as one point on the NFRC located in the vicinity of the resonance of interest calculated for a forcing amplitude $f$ and corresponding to a frequency $\omega_{(0)}$. The Fourier coefficients of this point are rotated until the $k$-th harmonic of the displacement is such that $s_k=0$ and $c_k<0$. Since $\mu>0$, this latter condition ensures that the feedback acts as a negative damping term. The initial gain $\mu_{(0)}$ is $-f/c_k$. This strategy is easily implemented for resonances that appear in the direct continuation of the main branch of the NFRC, i.e., for primary and $k:1$ superharmonic resonances. Because the other resonances appear as isolated branches, the computation of basins of attraction is required to provide the initial guess. 

Finally, we note that an alternative implementation could (i) enforce directly the desired phase lag as the phase condition, and (ii) express the velocity feedback as $\mu\omega_kA_k\sin{\omega t}$, as suggested by Equation \eqref{eq:velocity_forcing}.




\section{Numerical demonstration}
The algorithm for the PRNM computation is demonstrated in this section using the one- and two-degree-of-freedom systems introduced in Sections \ref{SECTION:DYNDUFFING} and \ref{motiv}, respectively.

\subsection{The Duffing oscillator}

\begin{figure}[ht] 
  \begin{subfigure}[b]{0.5\linewidth}
    \centering
    \includegraphics[width=1\linewidth]{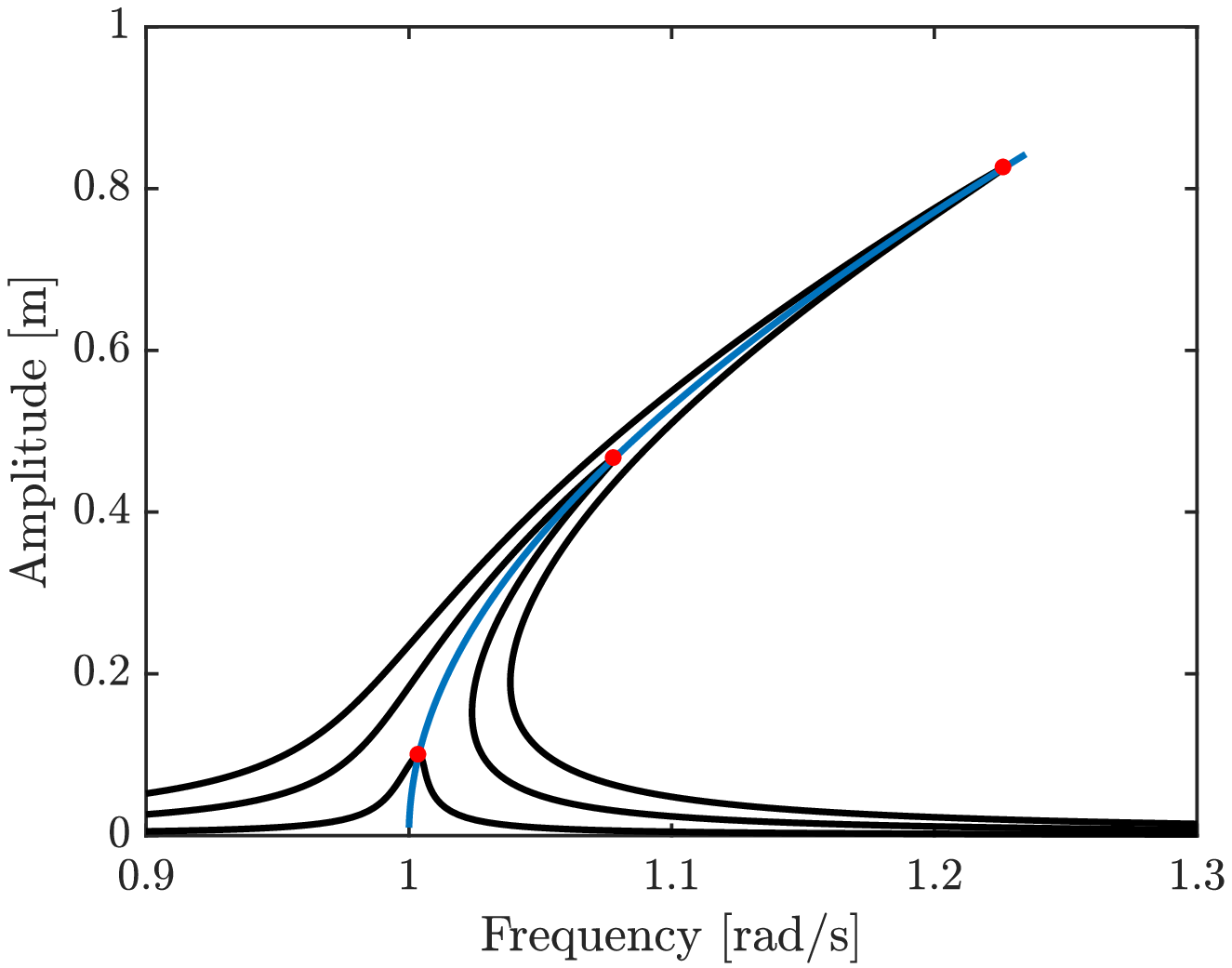} 
    \caption{\label{fig:FUNDRED_NFRC_PRNM}}
  \end{subfigure}
  \begin{subfigure}[b]{0.5\linewidth}
    \centering
    \includegraphics[width=1\linewidth]{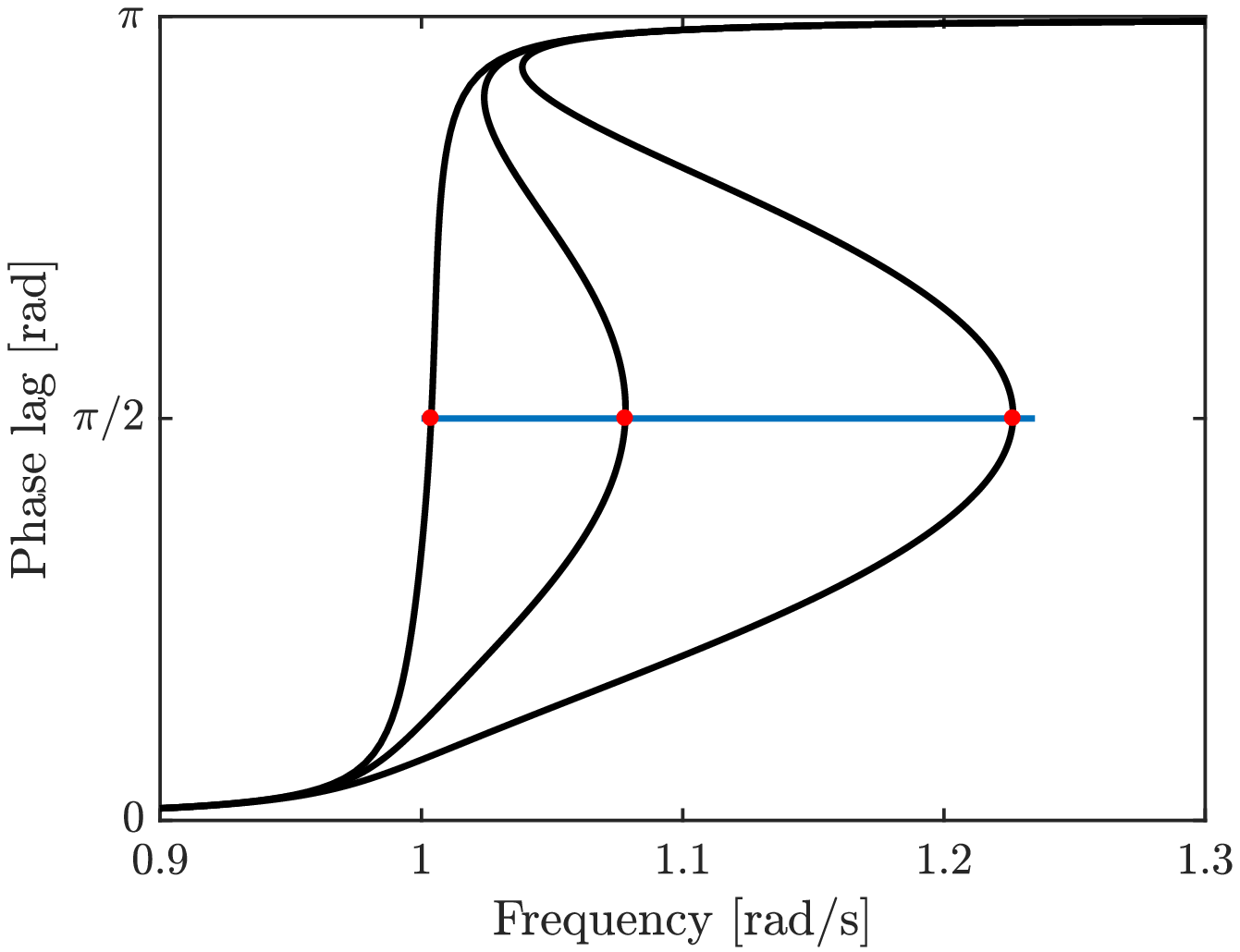}
    \caption{\label{fig:FUNDRES_PHASELAG_PRNM}}
  \end{subfigure} 
   \begin{subfigure}[b]{0.5\linewidth}
    \centering
    \includegraphics[width=1\linewidth]{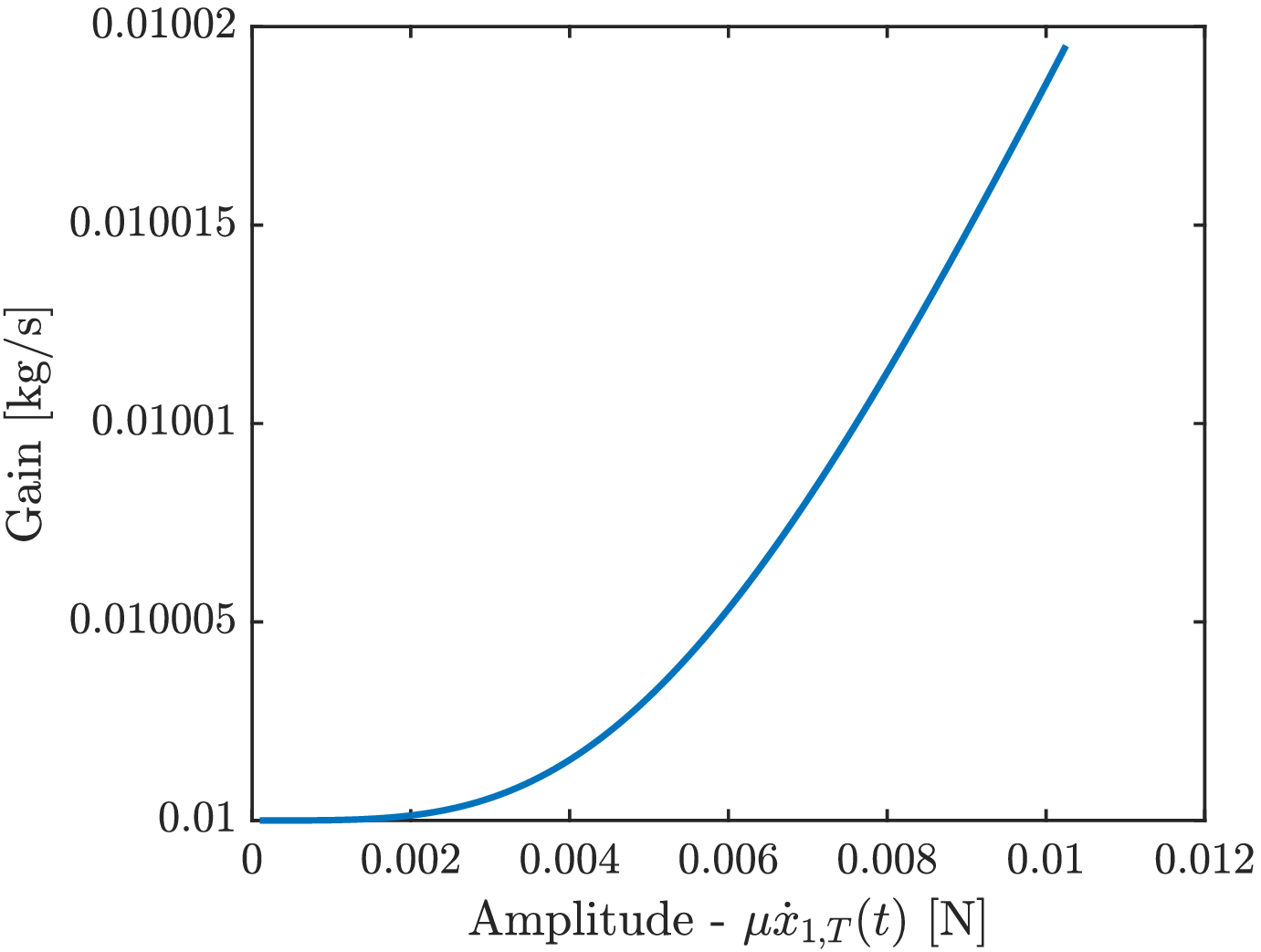}
    \caption{\label{fig:FUNDRES_mu_F}}
  \end{subfigure}
  \begin{subfigure}[b]{0.5\linewidth}
    \centering
    \includegraphics[width=1\linewidth]{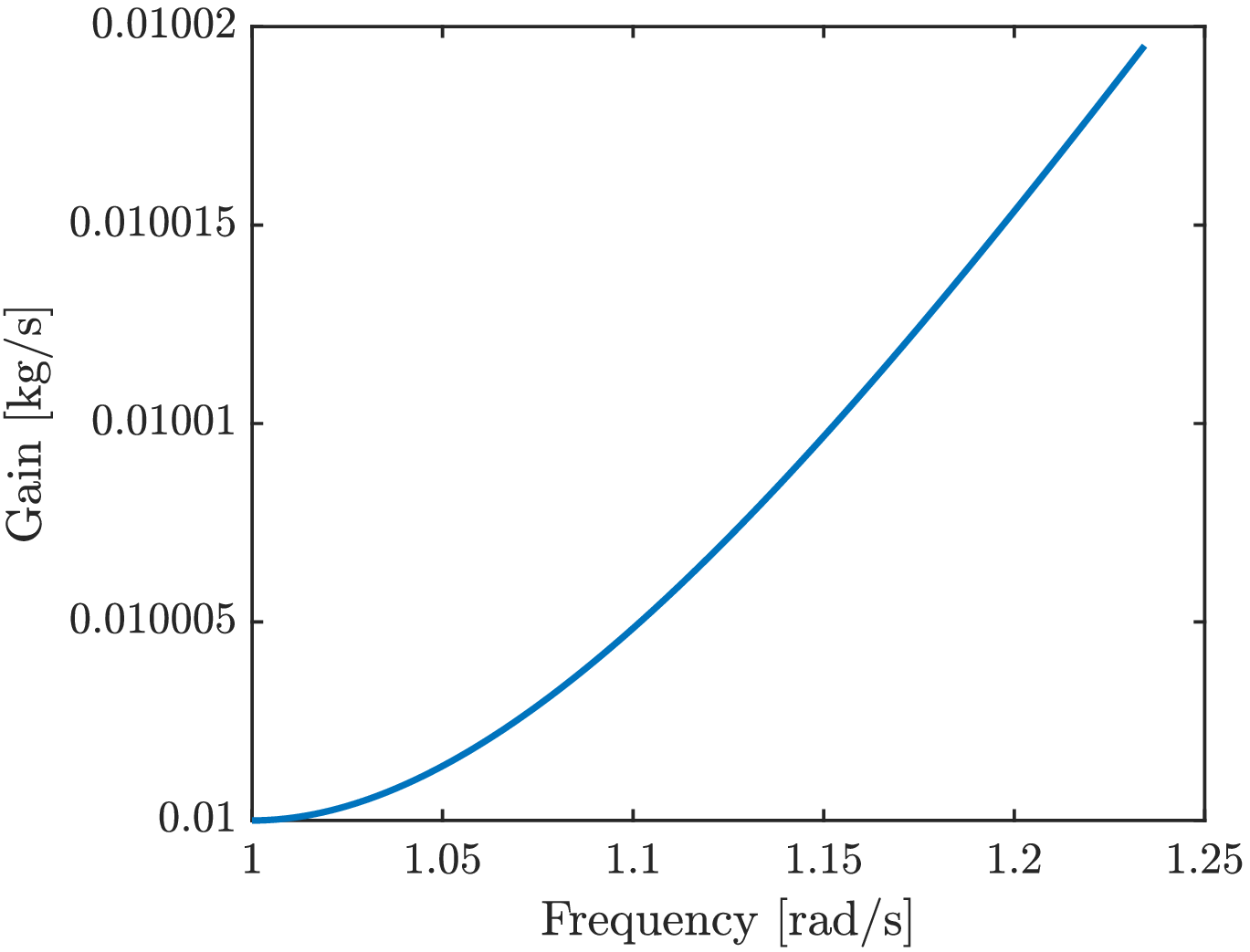}
    \caption{\label{fig:FUNDRES_mu_w0}}
  \end{subfigure}
  \caption{NFRCs (black) and PRNMs (blue) of the primary resonance of the Duffing oscillator for 3 forcing amplitudes: (\subref{fig:FUNDRED_NFRC_NO_PRNM}) amplitude, (\subref{fig:FUNDRES_PHASELAG_NOPRNM}) phase lag, (\subref{fig:FUNDRES_mu_F}) gain vs. forcing, (\subref{fig:FUNDRES_mu_w0}) gain vs. frequency.}
  \label{fig:FUNDRES_PRNM} 
\end{figure}

\begin{figure}[ht] 
  \begin{subfigure}[b]{0.5\linewidth}
    \centering
    \includegraphics[width=1\linewidth]{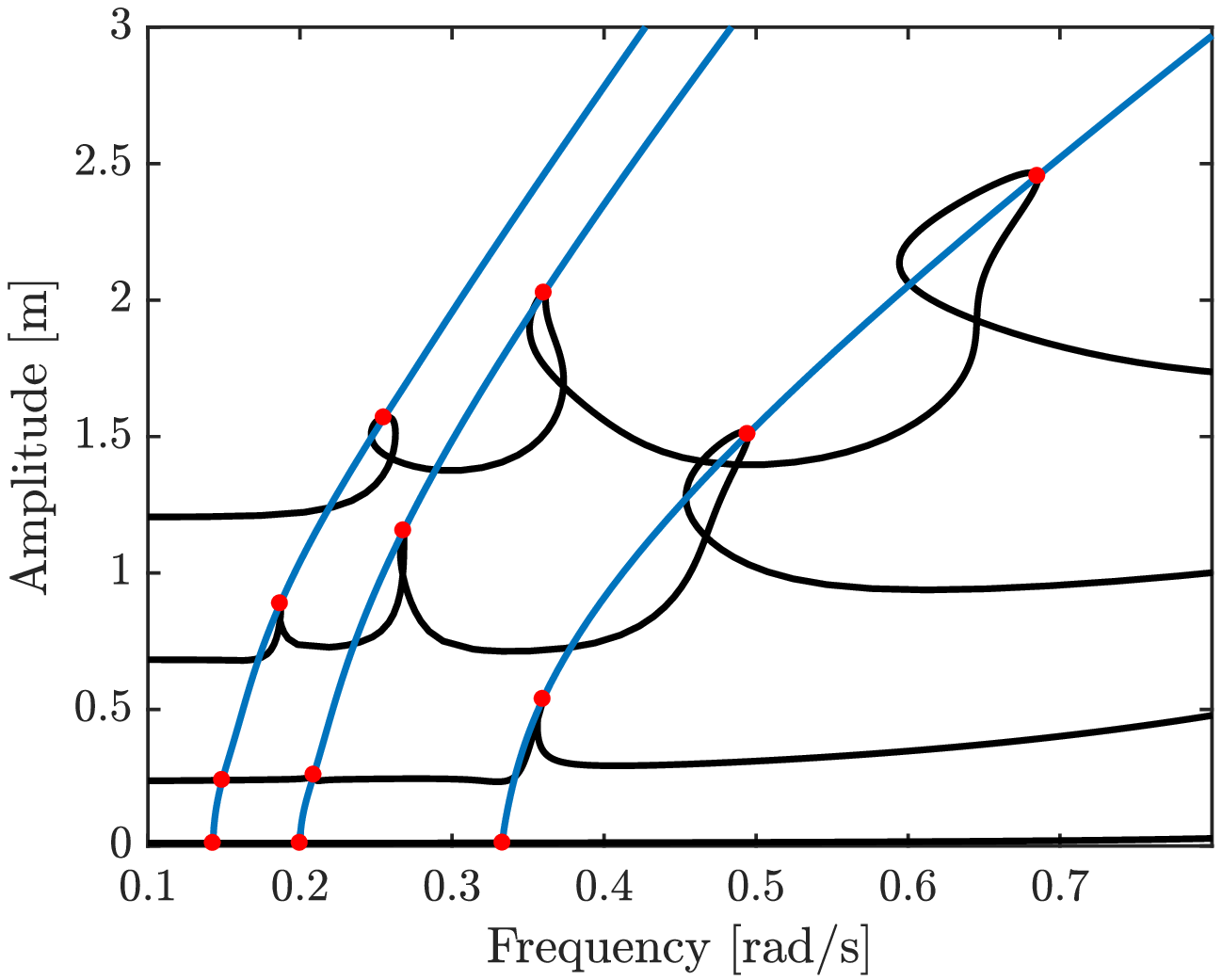} 
    \caption{\label{fig:SUPERHODD_NFRC_PRNM}}
  \end{subfigure}
  \begin{subfigure}[b]{0.5\linewidth}
    \centering
    \includegraphics[width=1\linewidth]{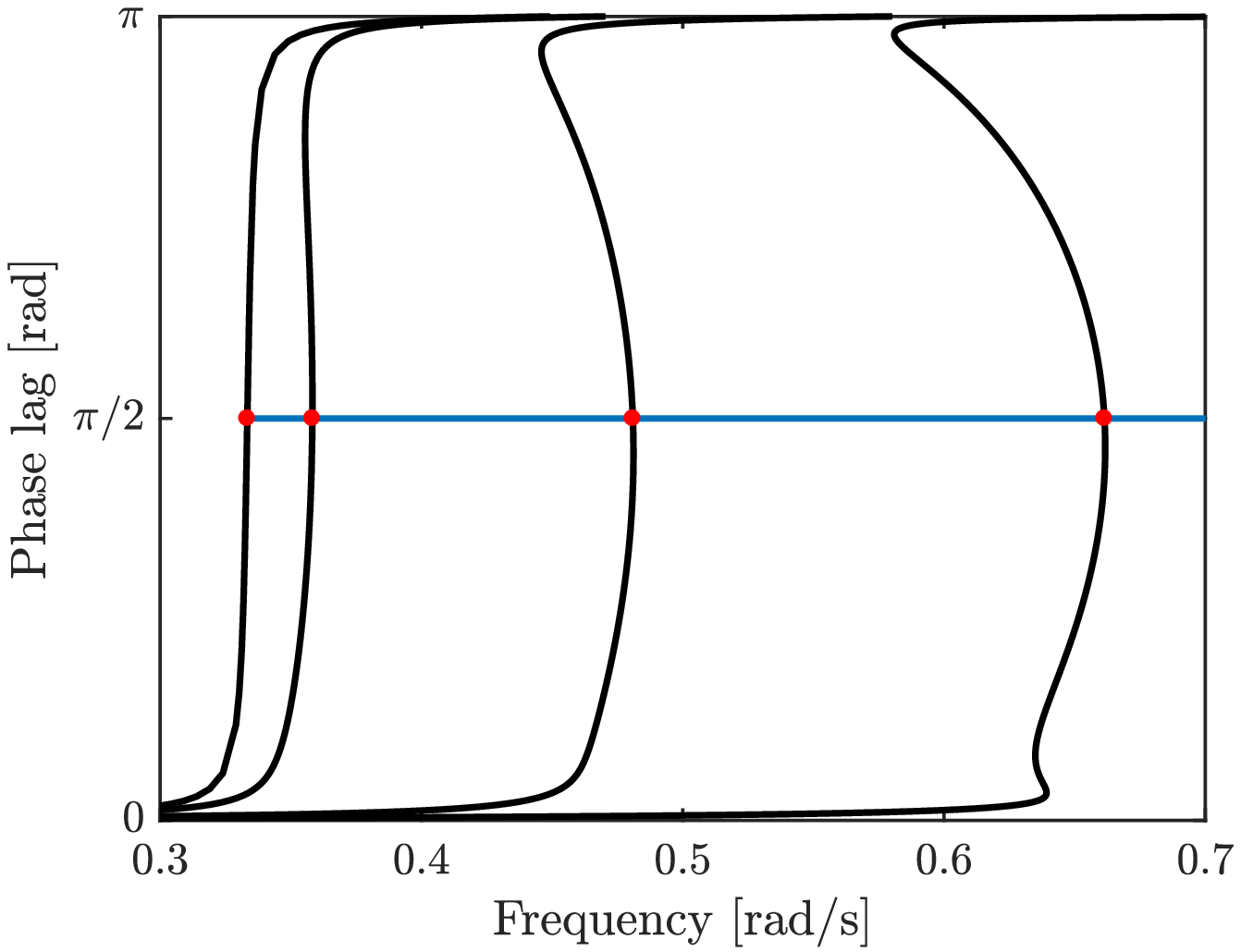}
    \caption{\label{fig:SUPERHODD_PHASELAG_PRNM}}
  \end{subfigure} 
  \caption{NFRCs (black) and PRNMs (blue) of the odd superharmonic resonances for 4 forcing amplitudes: (\subref{fig:SUPERHODD_NFRC_NO_PRNM}) amplitude and (\subref{fig:SUPERHODD_PHASELAG_NOPRNM}) phase lag of the $3^{rd}$ harmonic component of the 3:1 superharmonic resonance.}
  \label{fig:SUPHERHODD_PRNM} 
\end{figure}

\begin{figure}[ht] 
  \begin{subfigure}[b]{0.5\linewidth}
    \centering
    \includegraphics[width=1\linewidth]{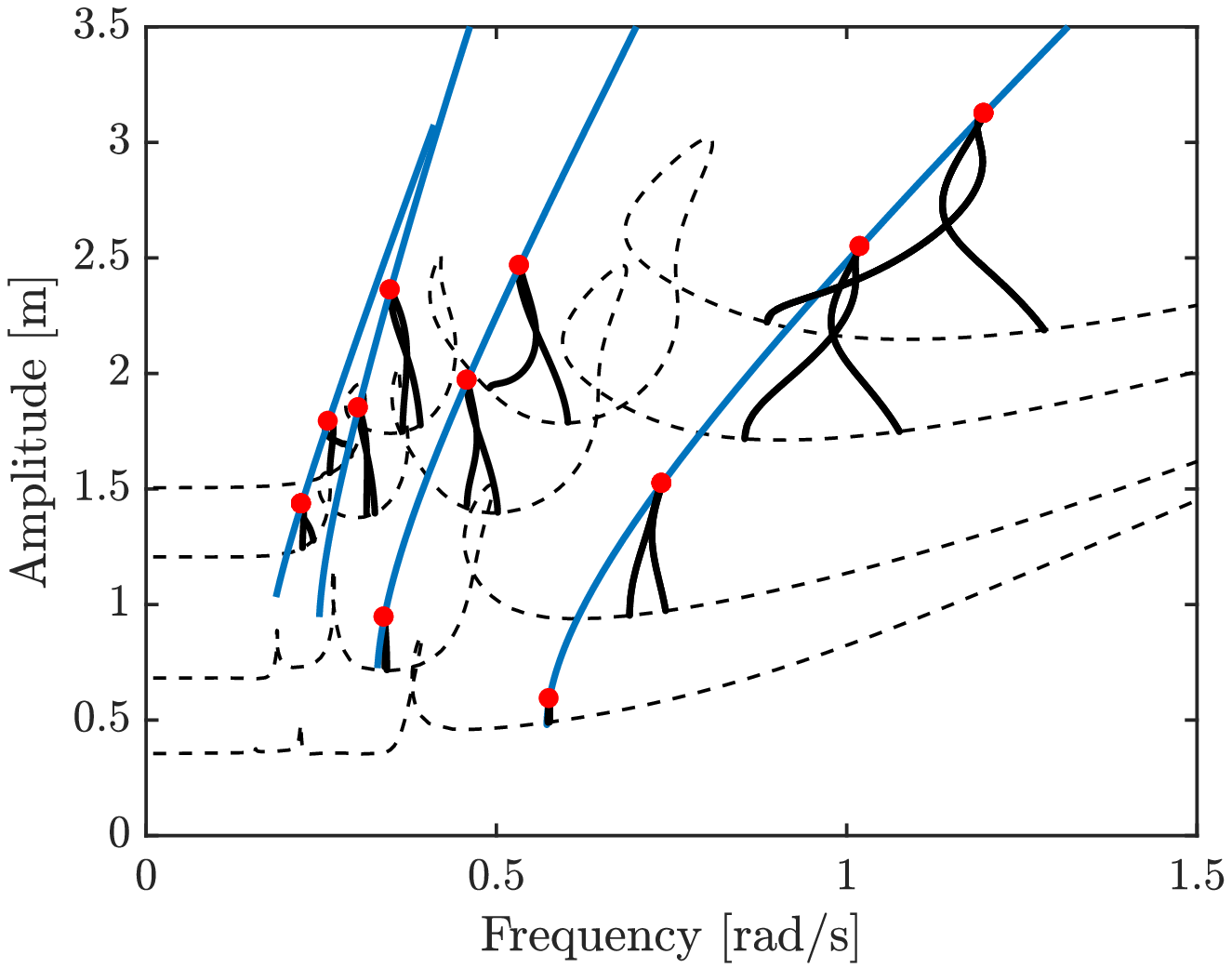}
    \caption{\label{fig:SUPERHEVEN_NFRC_PRNM}}
  \end{subfigure}
  \begin{subfigure}[b]{0.5\linewidth}
    \centering
    \includegraphics[width=1\linewidth]{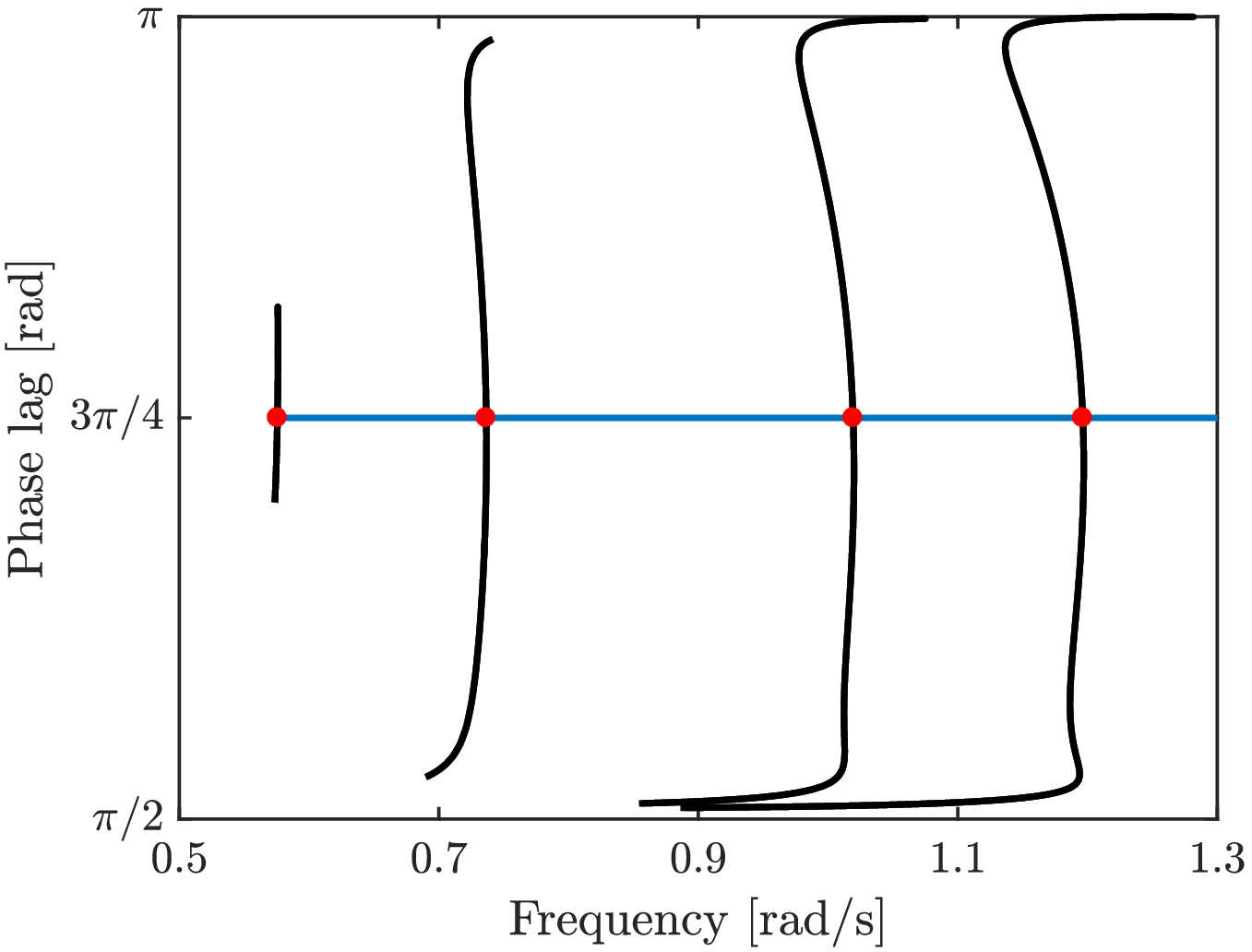}
    \caption{\label{fig:SUPERHEVEN_PHASELAG_PRNM}}
  \end{subfigure} 
  \caption{NFRCs (black) and PRNMs (blue) of the even superharmonic resonances for 4 forcing amplitudes: (\subref{fig:SUPERHEVEN_NFRC_NO_PRNM}) amplitude and (\subref{fig:SUPERHEVEN_PHASELAG_NO_PRNM}) phase lag of the $2^{nd}$ harmonic component for the 2:1 superharmonic resonance.}
  \label{fig:SUPHERHEVEN_PRNM} 
\end{figure}

\begin{figure}[ht] 
  \begin{subfigure}[b]{0.5\linewidth}
    \centering
    \includegraphics[width=1\linewidth]{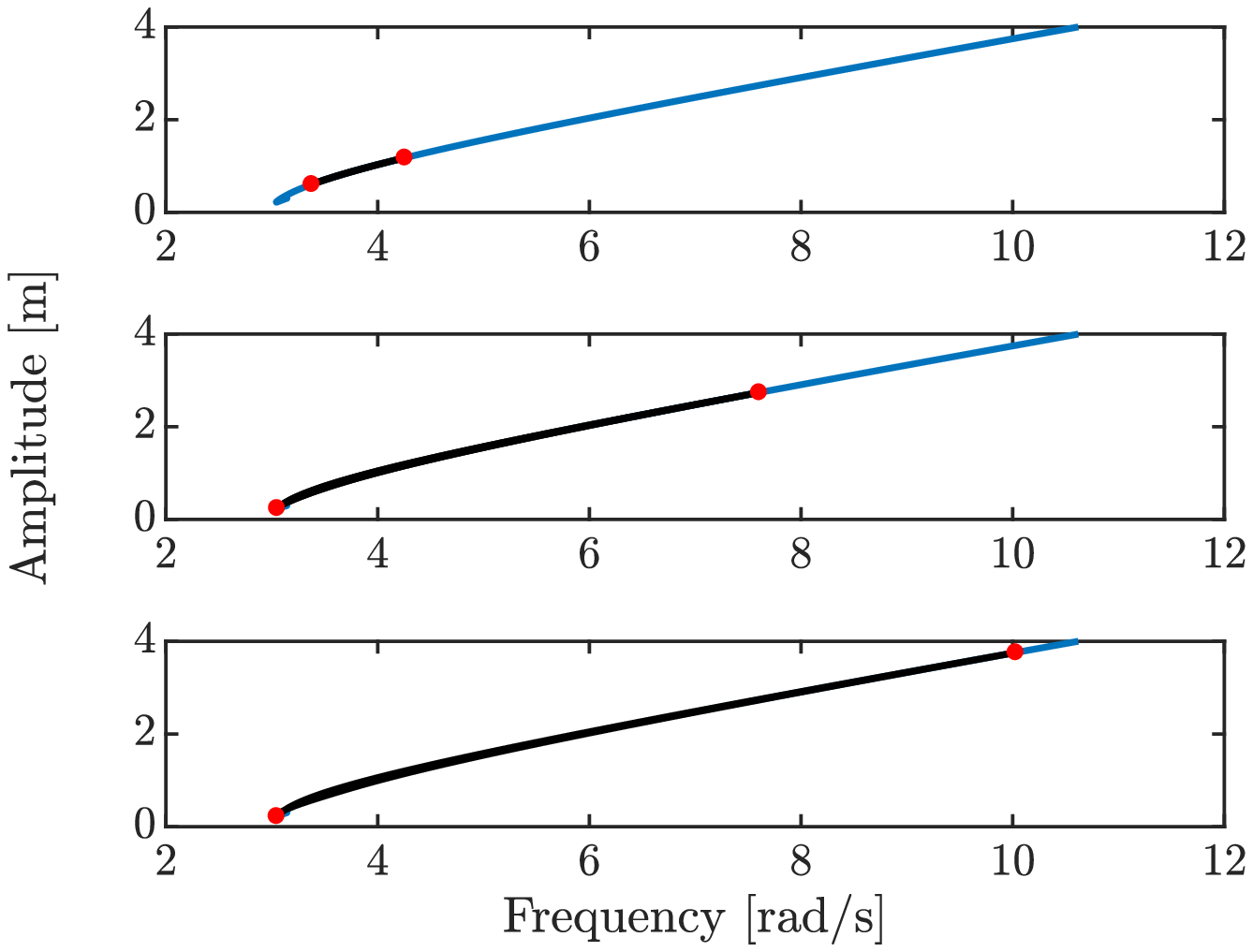}
    \caption{\label{fig:SUBHODD_NFRC_PRNM}}
  \end{subfigure}
  \begin{subfigure}[b]{0.5\linewidth}
    \centering
    \includegraphics[width=1\linewidth]{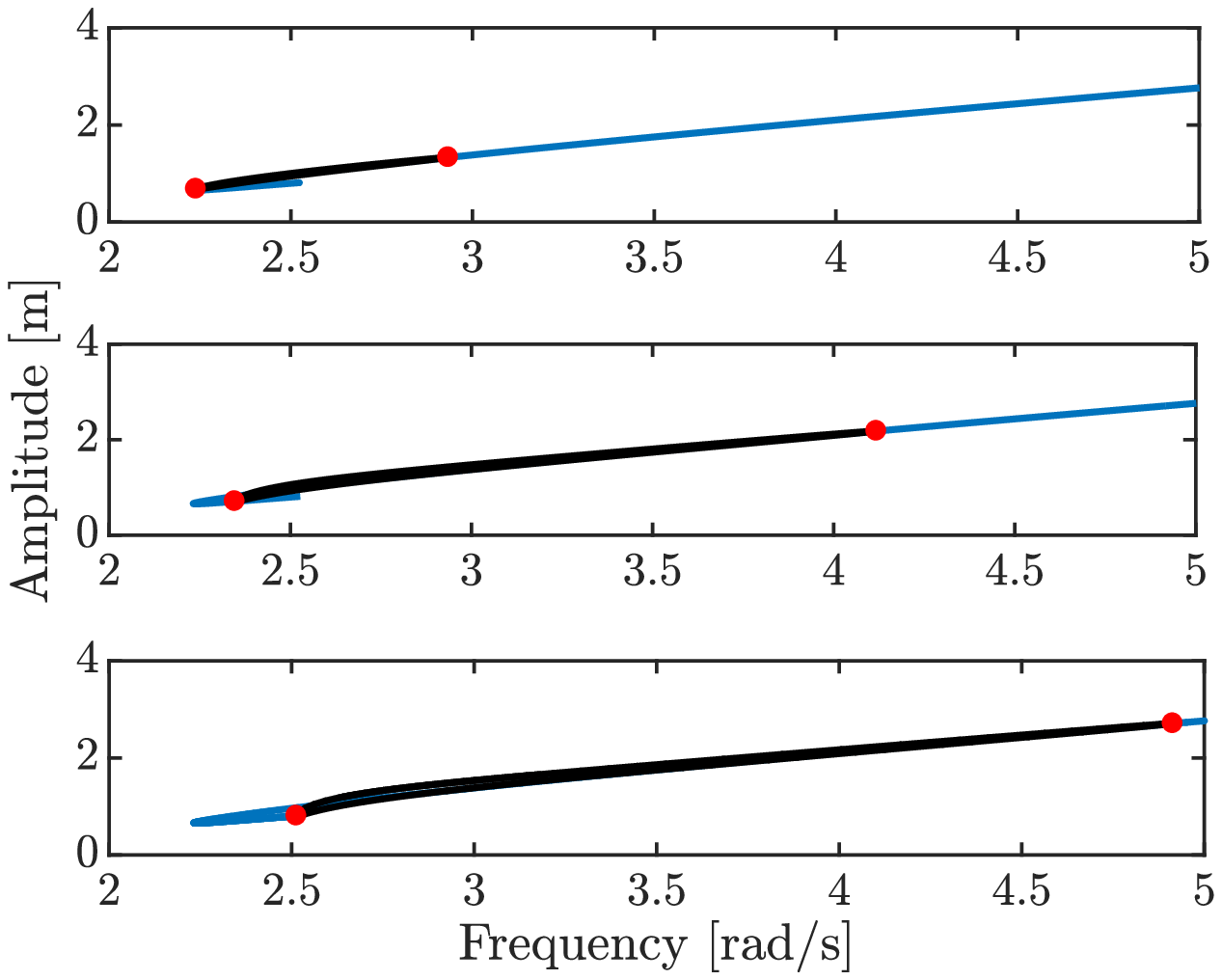}
    \caption{\label{fig:SUBHEVEN_NFRC_PRNM}}
  \end{subfigure} 
  \caption{NFRCs (black) and PRNMs (blue) of the 1:3 and 1:2 subharmonic resonances for 3 forcing amplitudes: (\subref{fig:SUBHODD_NFRC_NO_PRNM}) 1:3 and (\subref{fig:SUBHODD_PHASE_NO_PRNM}) 1:2.}
  \label{fig:SUBH_PRNM} 
\end{figure}

\begin{figure}[ht] 
  \begin{subfigure}[b]{0.5\linewidth}
    \centering
    \includegraphics[width=1\linewidth]{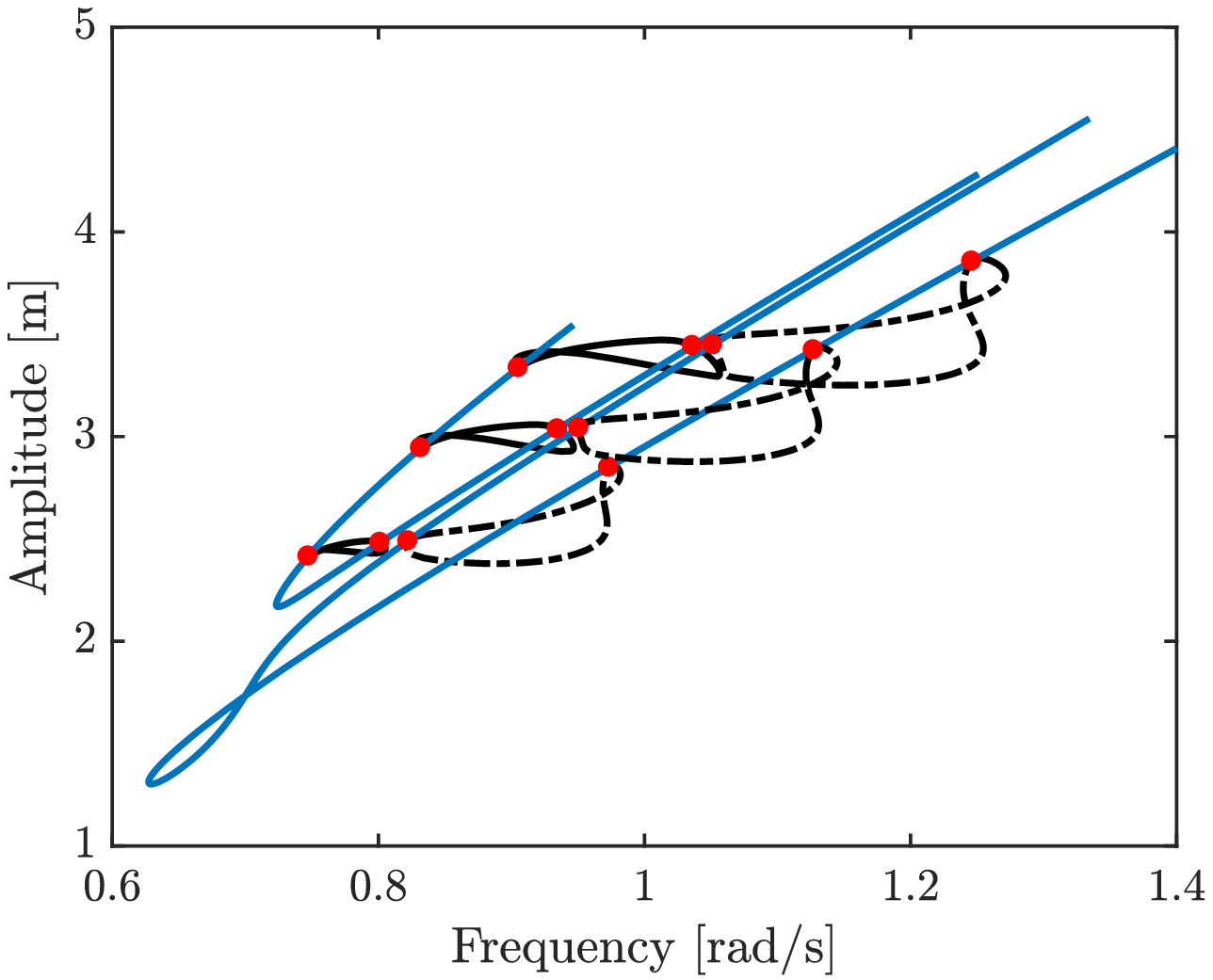}
    \caption{\label{fig:RATIONAL_PRNM_7_3}}
  \end{subfigure}
  \begin{subfigure}[b]{0.5\linewidth}
    \centering
    \includegraphics[width=1\linewidth]{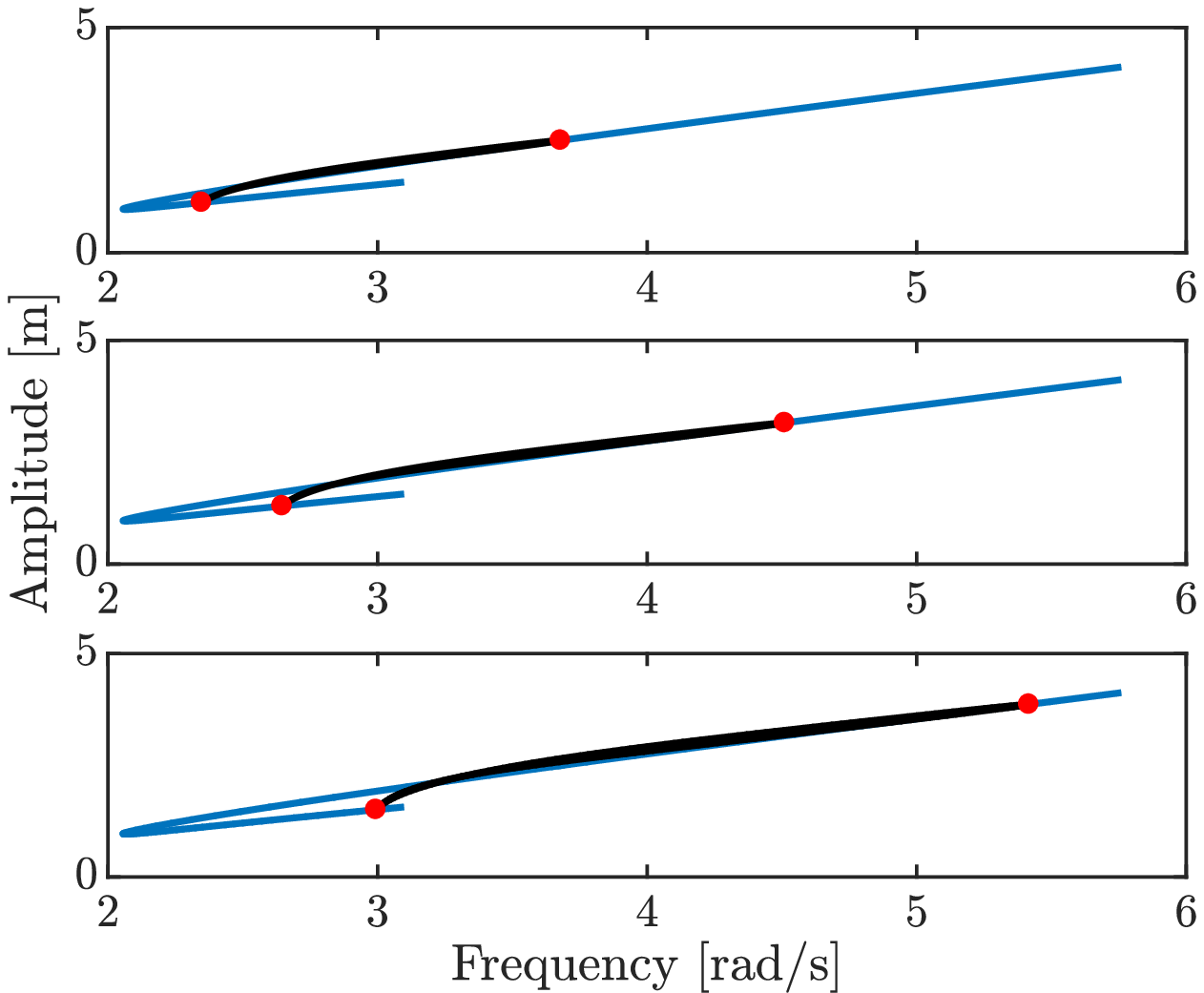}
    \caption{\label{fig:RATIONAL_PRNM_3_5}}
  \end{subfigure}
  
  \begin{subfigure}[b]{0.5\linewidth}
    \centering
    \includegraphics[width=1\linewidth]{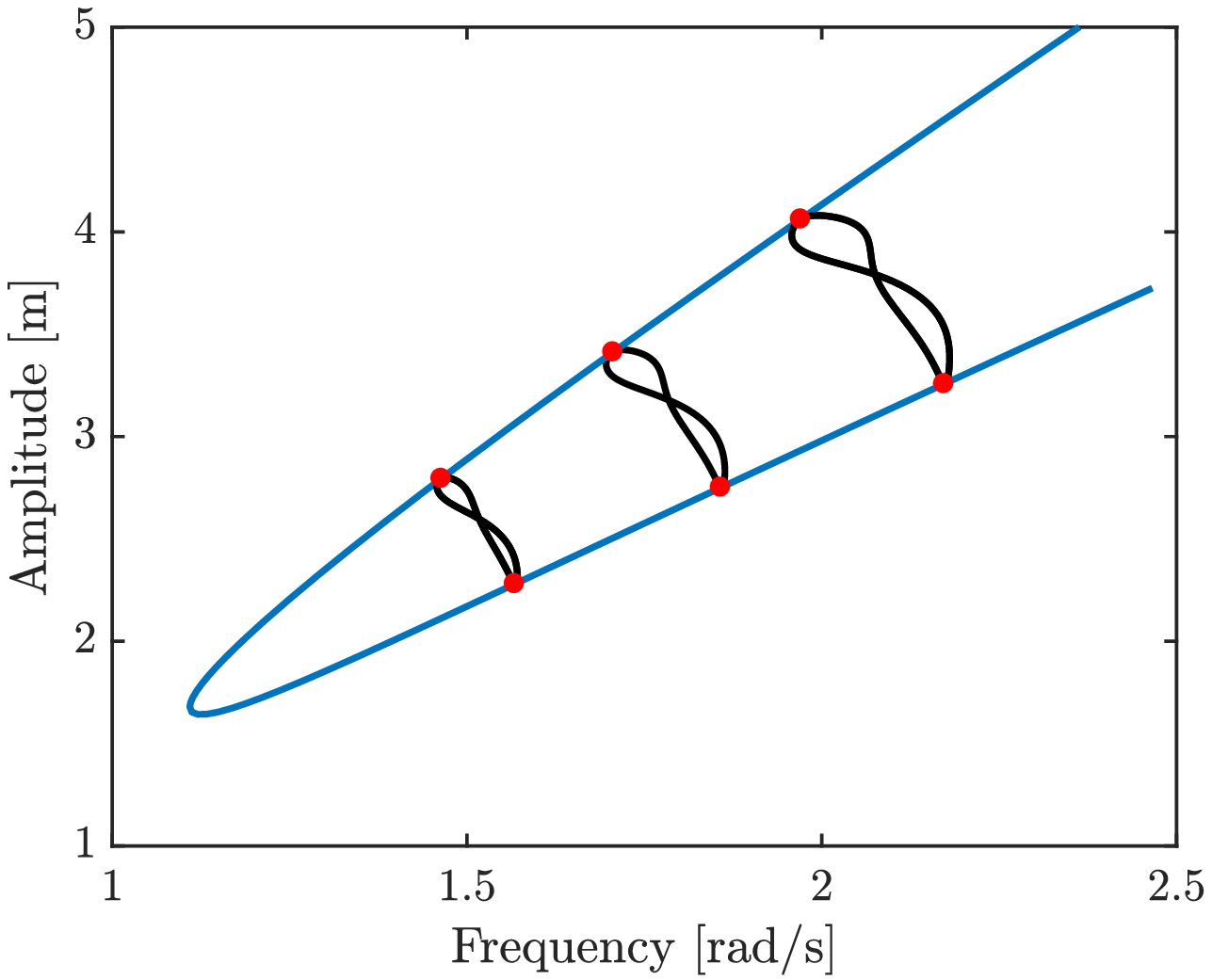}
    \caption{\label{fig:RATIONAL_PRNM_3_2}}
  \end{subfigure}
  \begin{subfigure}[b]{0.5\linewidth}
    \centering
    \includegraphics[width=1\linewidth]{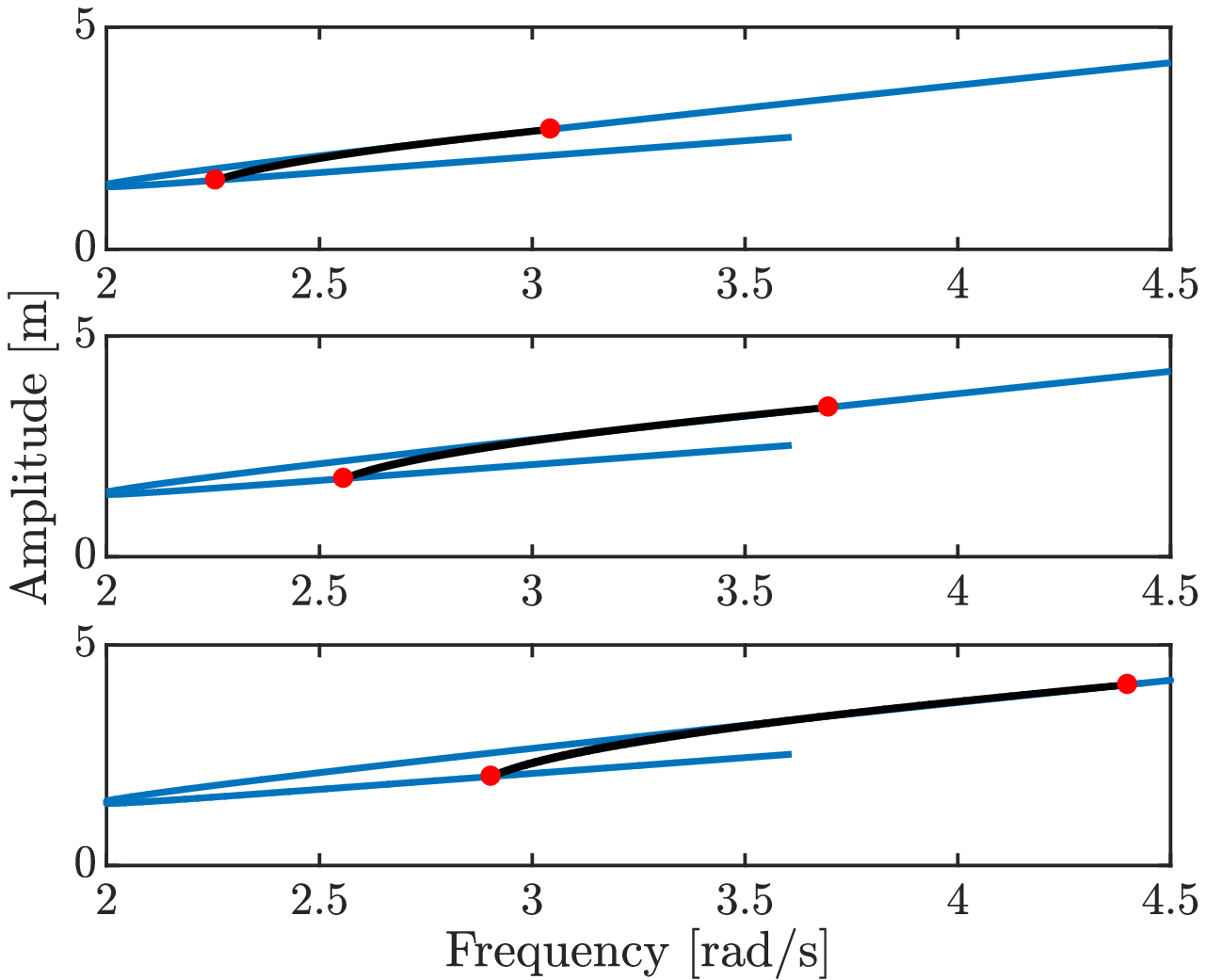}
    \caption{\label{fig:RATIONAL_PRNM_3_4}}
  \end{subfigure}
  
  \begin{subfigure}[b]{0.5\linewidth}
    \centering
    \includegraphics[width=1\linewidth]{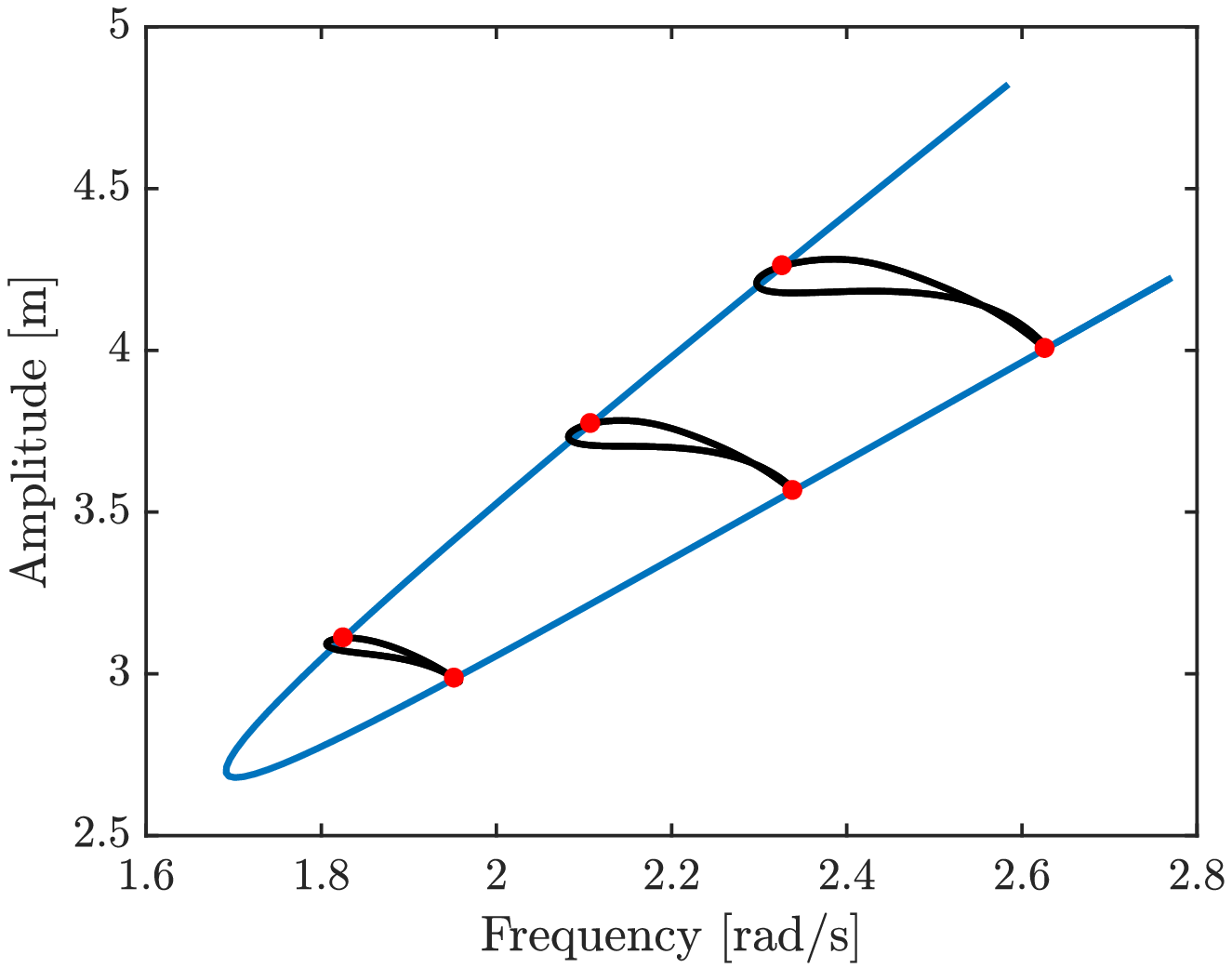}
    \caption{\label{fig:RATIONAL_PRNM_4_3}}
  \end{subfigure}
  \begin{subfigure}[b]{0.5\linewidth}
    \centering
    \includegraphics[width=1\linewidth]{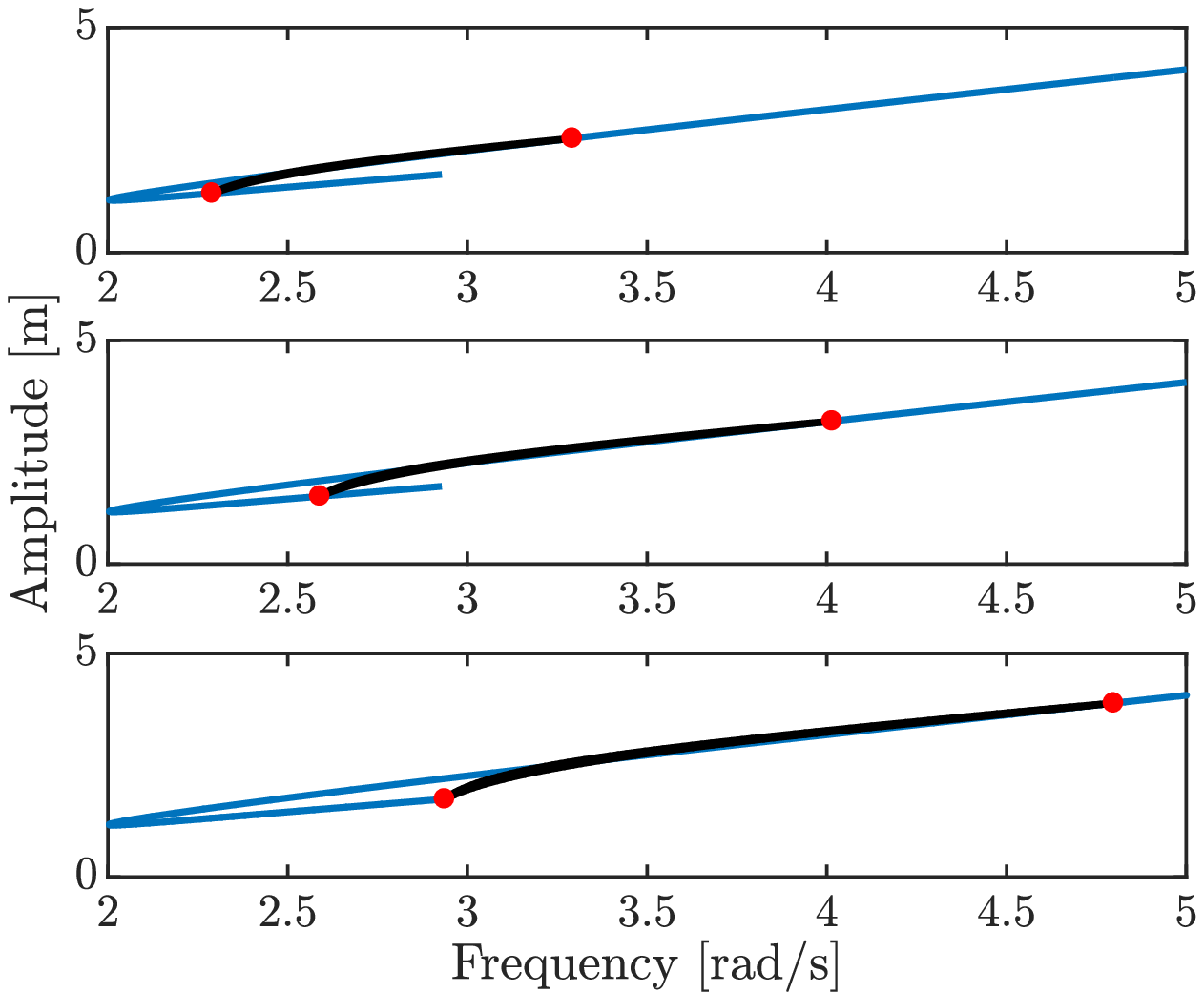}
    \caption{\label{fig:RATIONAL_PRNM_2_3}}
  \end{subfigure}
  \caption{NFRCs (black) and PRNMs (blue) for (\subref{fig:RATIONAL_PRNM_7_3}) 7:3  (\subref{fig:RATIONAL_PRNM_3_5}) 3:5 (\subref{fig:RATIONAL_PRNM_3_2}) 3:2  (\subref{fig:RATIONAL_PRNM_3_4}) 3:4  (\subref{fig:RATIONAL_PRNM_4_3}) 4:3  (\subref{fig:RATIONAL_PRNM_2_3}) 2:3.}
  \label{fig:RATIONAL_PRNM} 
\end{figure}

The PRNMs were calculated for the primary resonance (Figure \ref{fig:FUNDRES_PRNM}), for superharmonic resonances (Figures \ref{fig:SUPHERHODD_PRNM} and \ref{fig:SUPHERHEVEN_PRNM}), for subharmonic resonances (Figure \ref{fig:SUBH_PRNM}) and for ultra-subharmonic resonances (Figure \ref{fig:RATIONAL_PRNM}). The PRNMs are represented using blue lines in these figures whereas the red dots represent the nonlinear phase resonance conditions defined in this paper, i.e., either $\pi/2$ or $3\pi/4\nu$. It can be observed that the PRNMs pass through all red points, corresponding to the nonlinear phase resonances discussed in the paper, which confirms the accuracy of the proposed computational framework.

For the primary resonance, Figures \ref{fig:FUNDRES_mu_F} and \ref{fig:FUNDRES_mu_w0} show the evolution of the gain $\mu$ with respect to the amplitude and frequency of the forcing, respectively. As expected, when $\omega$ is close to $\omega_0$, \textit{i.e.}, the system behaves almost linearly, $\mu$ tends to the value of the damping coefficient $c=0.01kg/s$.


For subharmonic and ultra-subharmonic resonances, it was shown that they exhibit two phase resonance points. Figures \ref{fig:SUBH_PRNM} and \ref{fig:RATIONAL_PRNM} confirm that the continuation of the PRNMs curve is able to capture both of them. In fact, such resonances (and hence the PRNMs) exist only above a critical forcing amplitude. At the forcing amplitude for which a PRNM first appears, the PRNMs curve features a turning point corresponding to the PRNM with the minimum frequency.

\subsection{A two-degree-of-freedom system}

\begin{figure}[ht] 
  \begin{subfigure}[b]{0.5\linewidth}
    \centering
    \includegraphics[width=1\linewidth]{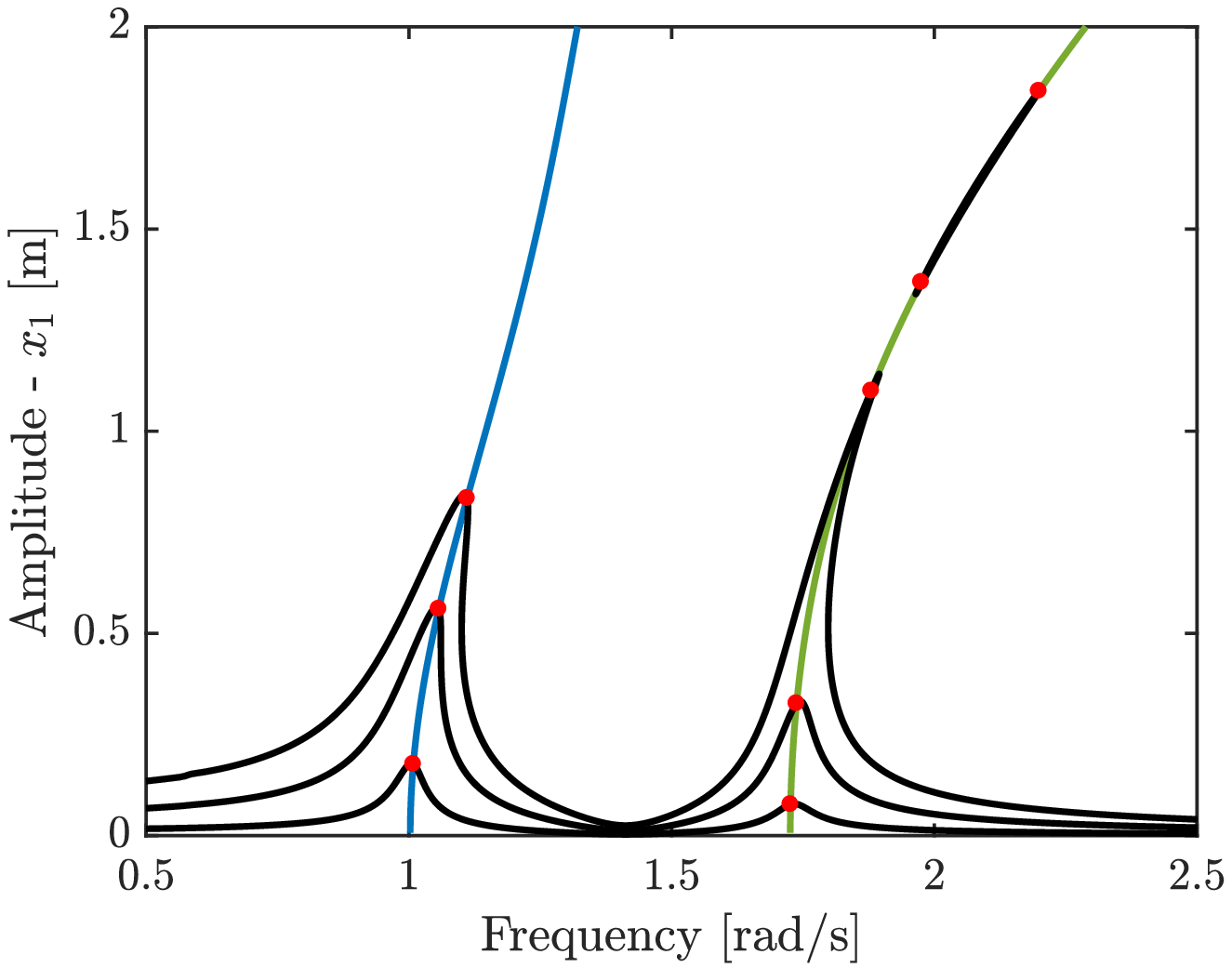}
    \caption{\label{fig:2DOFs_3NFRC_PLNM}}
  \end{subfigure}
  \begin{subfigure}[b]{0.5\linewidth}
    \centering
    \includegraphics[width=1\linewidth]{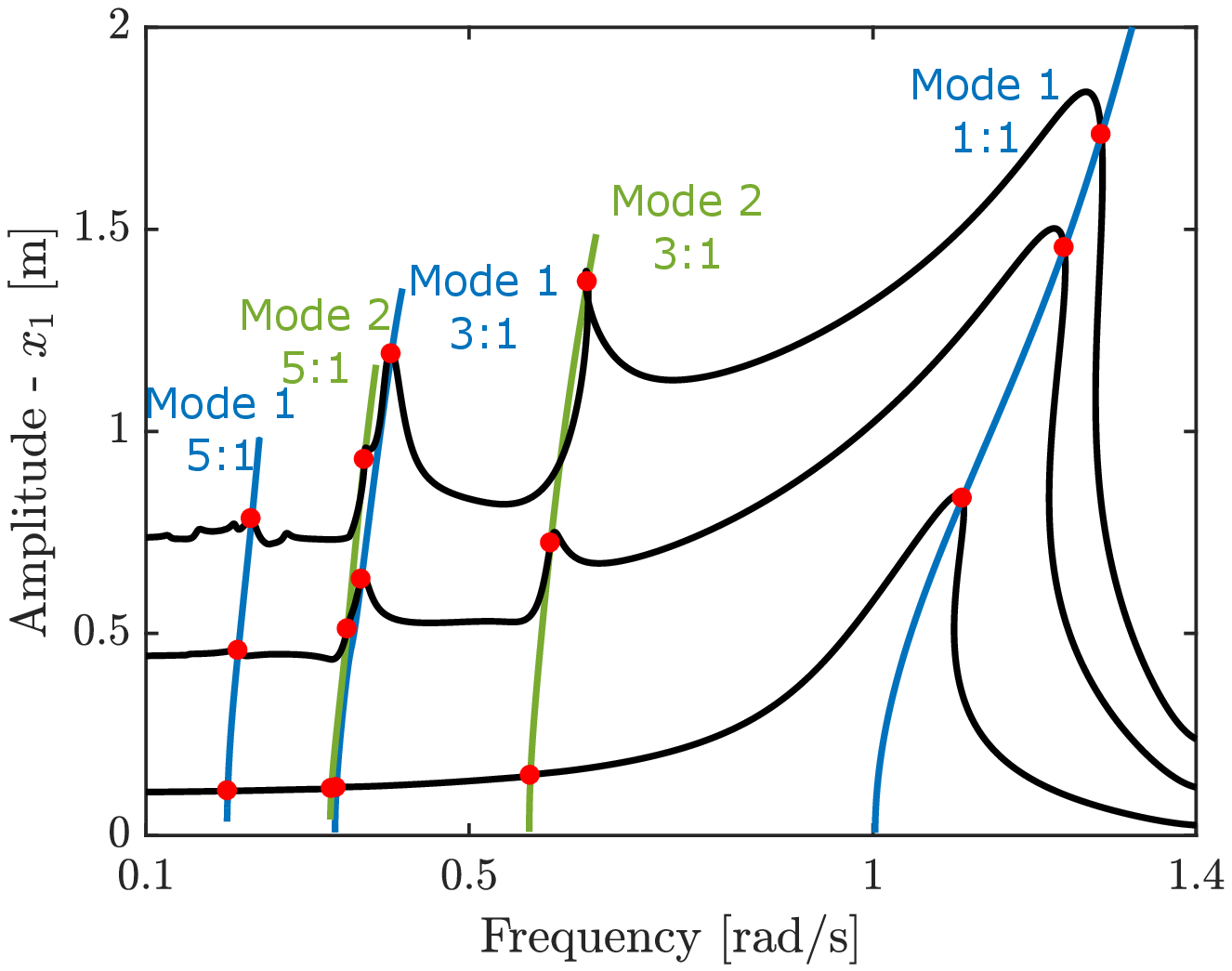}
    \caption{\label{fig:2DOFs_SUPERH}}
  \end{subfigure} 
  \caption{NFRCs (black) and PRNMs (blue: mode 1 and green: mode 2) for the 2DOF system: (\subref{fig:2DOFs_3NFRC_PLNM}) The two vibration modes ($f=0.02$ to $0.161N$) and (\subref{fig:2DOFs_SUPERH}) close-up around the superharmonic resonances ($f=0.161$ to $1.5N$). Black: NFRC; green: PRNMs of mode 1; blue: PRNMs of mode 2.}
  \label{fig:2DOFs} 
\end{figure}

To establish that the usefulness of PRNMs extends beyond  single-degree-of-freedom oscillators, the two-degree-of-freedom introduced in Section \ref{motiv} is reconsidered. The NFRCs for different forcing amplitudes together with the PRNMs and the nonlinear phase resonance conditions are displayed in Figure \ref{fig:2DOFs}. Each mode features 3:1 and 5:1 superharmonic resonances in Figure \ref{fig:2DOFs_SUPERH}. In addition, when $f=0.161$N, a 1:1 isolated resonance branch related to the second mode can be observed in Figure \ref{fig:2DOFs_3NFRC_PLNM}. Many other resonances exist in the system and can be captured by the proposed technique but are not shown for brevity and are thus not further discussed herein. 

As for the Duffing oscillator, the PRNMs pass through all red points and are found to offer a very accurate characterization of the different resonances at hand, including the isolated resonance branch for the second mode.

\begin{figure}[ht] 
  \begin{subfigure}[b]{0.5\linewidth}
    \centering
    \includegraphics[width=1\linewidth]{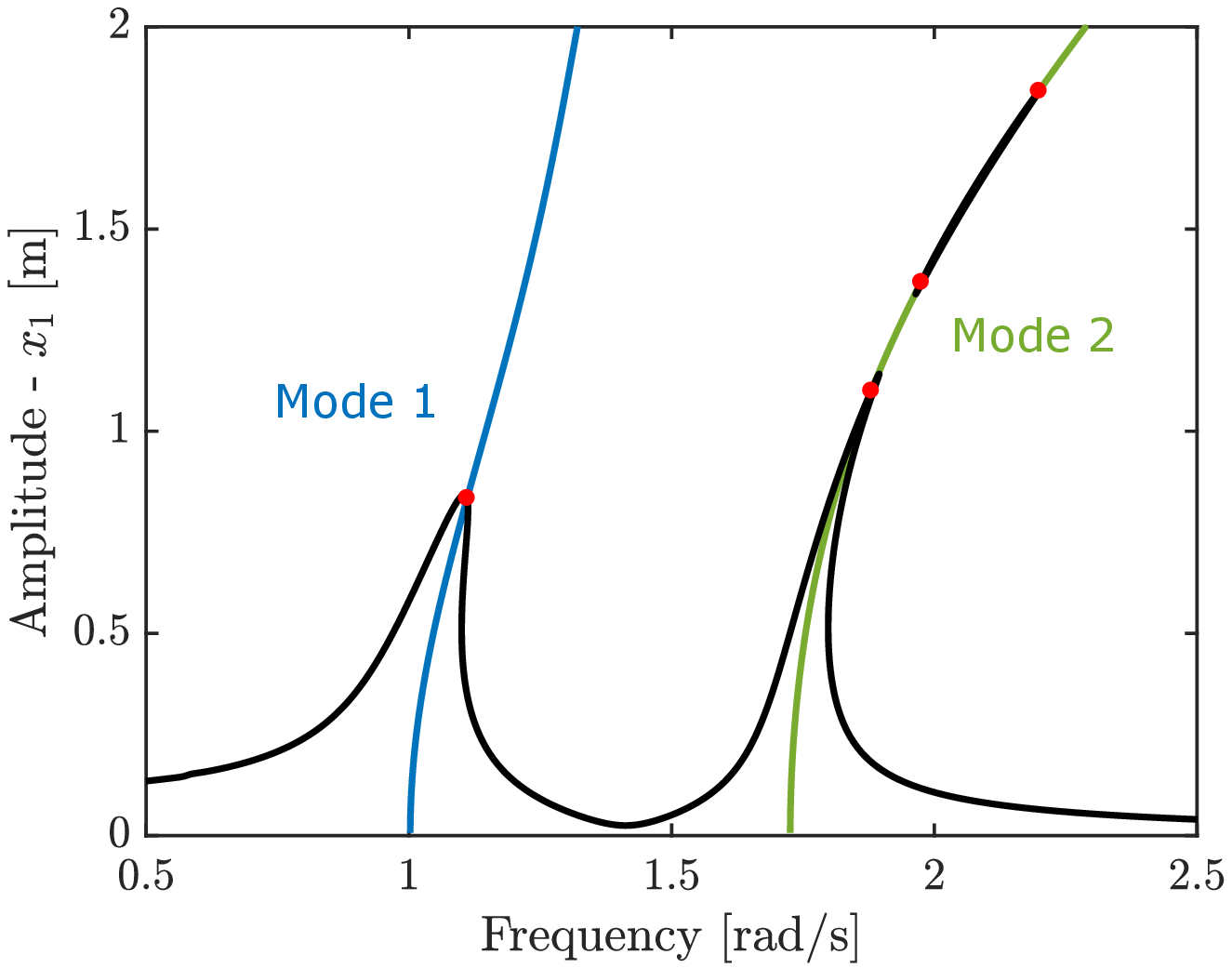}
    \caption{\label{fig:2DOFs_NFRC_PLNM}}
  \end{subfigure}
  \begin{subfigure}[b]{0.5\linewidth}
    \centering
    \includegraphics[width=1\linewidth]{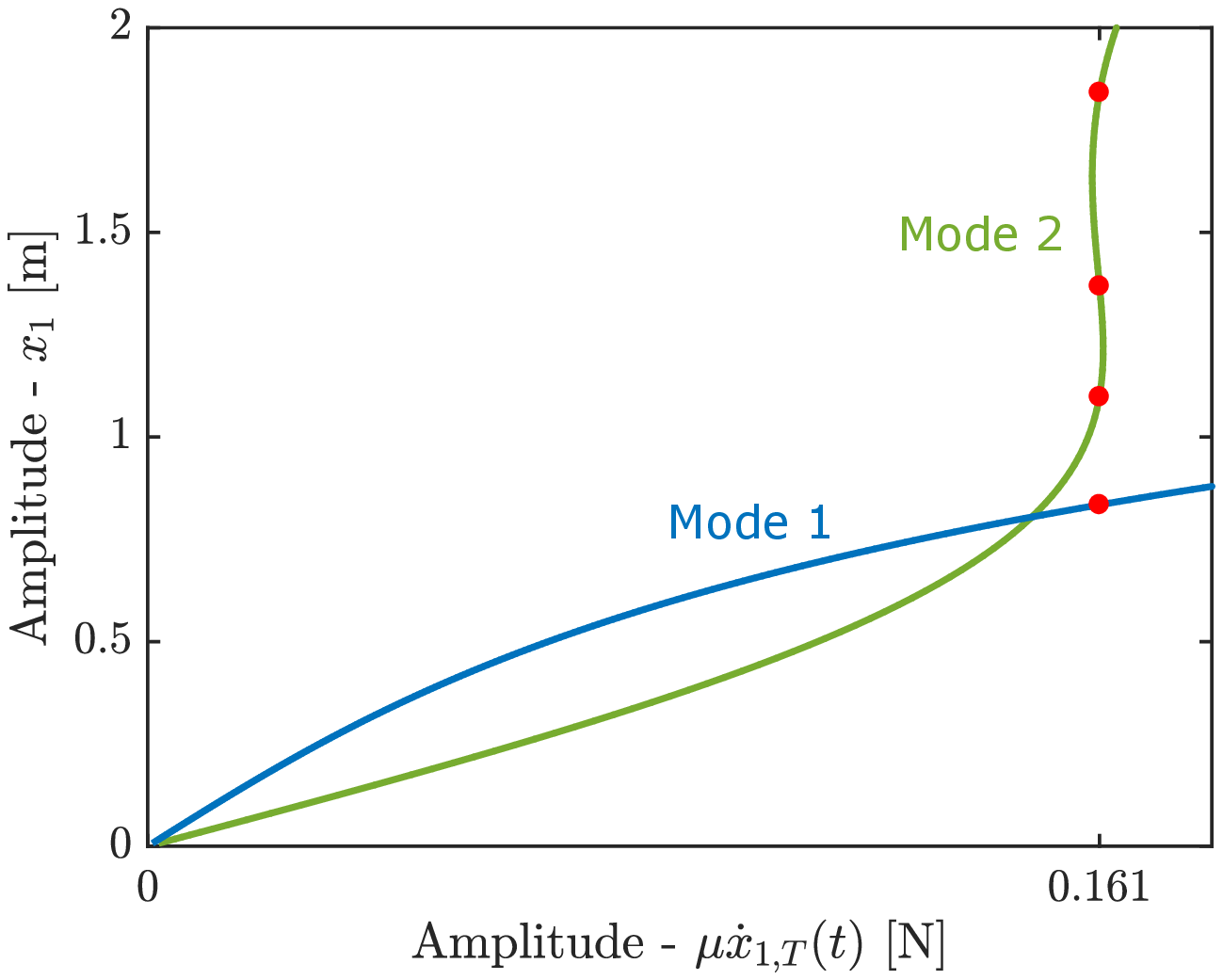}
    \caption{\label{fig:2DOFs_FORCING_AMPLITUDE}}
  \end{subfigure}
  \begin{subfigure}[b]{0.5\linewidth}
    \centering
    \includegraphics[width=1\linewidth]{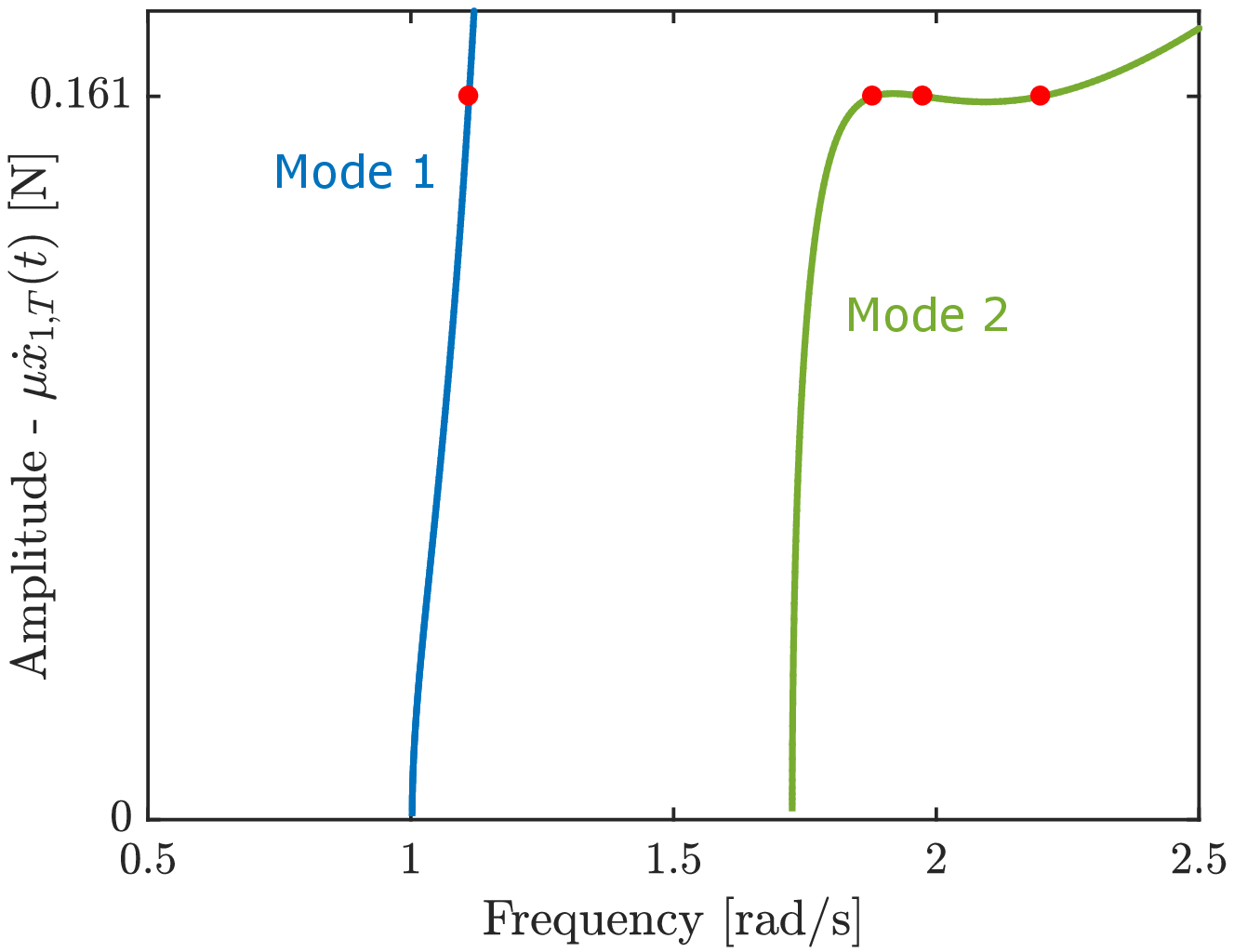}
    \caption{\label{fig:2DOFs_FORCING_FREQUENCY}}
  \end{subfigure}
  \caption{PRNMs (blue: mode 1 and green: mode 2) for the 2DOF system: (\subref{fig:2DOFs_NFRC_PLNM}) NFRC (black) for $f=0.161$N, (\subref{fig:2DOFs_FORCING_AMPLITUDE}) PRNM amplitude at the first degree of freedom and (\subref{fig:2DOFs_FORCING_FREQUENCY}) PRNM frequency.}
  \label{fig:2DOFs_ISOLA} 
\end{figure}

Another insightful graphical depiction exploits the fact that the velocity feedback plays the role of fictitious forcing. It is therefore possible to plot the evolution of the amplitude and frequency of the PRNMs with respect to the external forcing, as achieved in Figures  \ref{fig:2DOFs_FORCING_AMPLITUDE} and \ref{fig:2DOFs_FORCING_FREQUENCY}, respectively. Interestingly, for $f=0.161$N in Figure \ref{fig:2DOFs_FORCING_AMPLITUDE}, three resonance points exist for the second mode (see the red dots). The point with the lowest amplitude corresponds to the usual resonance peak whereas the other two points are associated with the extremities of the isolated resonance. Thus, such a plot nicely reveals the existence of the isolated resonance \cite{Hill22}. The frequencies associated with the isolated branch can be read in Figure \ref{fig:2DOFs_FORCING_FREQUENCY}.


\section{Conclusion}
\label{SECTION:CONCLUSION}

This paper considered the resonant dynamics of nonlinear systems subjected to single-point, single-harmonic forcing. The concept of a nonlinear phase resonance, defined originally for primary resonances, was generalized to $k:\nu$ superharmonic, subharmonic, and ultra-subharmonic resonances. Specifically, it was shown that appropriate resonance conditions correspond to phase quadrature when $k$ and $\nu$ are both odd and to a phase lag equal to $3\pi/4\nu$ otherwise. A second contribution of this paper was to develop a harmonic balance-based computational framework which can accurately predict the mode shapes and oscillation frequencies at phase resonance. To do so, the harmonic forcing was replaced by a $T$-periodic velocity feedback comprising the $k$-th harmonic. 

The resulting nonlinear modes were termed phase resonance nonlinear modes (PRNMs). Because they rely on classical harmonic forcing, the PRNMs exhibit two key properties, namely (i) the PRNMs are actual periodic orbits on the NFRCs and (ii) the PRNMs correspond exactly to the nonlinear modes identified experimentally using phase-locked loops \cite{THOMAS,SCHEEL} and control-based continuation \cite{RENSON}. This paves the way for a rigorous correlation between numerical and experimental nonlinear modes. By contrast, NNMs are not guaranteed to be actual orbits on the NFRCs and may lack significance as discussed in \cite{CENEDESE,HILL2}. In addition, NNMs require an unpractical multi-point, multi-harmonic forcing to be identified experimentally.

Finally, this study considered low-dimensional, weakly to moderately damped systems with cubic stiffness vibrating in weakly and strongly nonlinear regimes of motion. Thus, the findings of this work  should be verified (or adapted) for other types of nonlinearities and for higher-dimensional systems with potentially coupled and highly-damped modes. Our next efforts will also calculate the stability and bifurcations of the PRNMs through the Floquet exponents, similarly to what was achieved for NFRCs in \cite{Detroux}.


    
    \bibliographystyle{unsrt}
    \bibliography{bibli}

\end{document}